\documentclass[10pt,a4paper]{article}

\usepackage{amsmath,amsfonts,amssymb,epsfig}
\usepackage[all]{xy}
\usepackage[latin1]{inputenc}
\newcommand{\atan}{\mbox{atan}}
\newcommand{\R}{{\mathbb{R}}}
\newcommand{\C}{{\mathbb{C}}}
\newcommand{\Z}{{\mathbb{Z}}}

\def\ha{\frac{1}{2}}

\def\pa{\partial}
\def\ra{\rightarrow}
\def\preuve{\begin{proof}}

\def\ga{\alpha}
\def\gb{\beta}

\def\ge{\varepsilon}

\def\gf{\varphi}
\def\gg{\gamma}

\def\gl{\lambda}
\def\go{\omega}

\def\gm{\mu }

\def\gs{\sigma}

\def\aujour{\ifnum\day=1 1\ermini\else\number\day\fi\
\ifcase\month\or janvier\or f\'evrier\or mars\or avril\or mai\or juin\or
juillet\or aout\or septembre\or octobre\or novembre\or d\'ecembre\fi\
\number\year}

\input{macro}

\begin{document}

\title{SINGULAR BOHR-SOMMERFELD RULES FOR 2D INTEGRABLE SYSTEMS}
\author{Yves Colin de Verdi\`ere$^{*,~**}$ and San V\~u Ng\d oc$^{**}$
  }

\date{May 22, 2000}

\maketitle

\noindent
$^{*}$ Institut Universitaire de France\\
$^{**}$ Institut Fourier, Unit{\'e} mixte de recherche CNRS-UJF 5582 \\
 BP 74, 38402-Saint Martin d'H\`eres Cedex (France)\\
yves.colin-de-verdiere@ujf-grenoble.fr\\
san.vu-ngoc@ujf-grenoble.fr\\

\vskip 1cm
\begin{center}
  { \large Prépublication de l'Institut Fourier
    n$^{\textrm{\underline{o}}}$~508 (2000)}
  \texttt{http://www-fourier.ujf-grenoble.fr/prepublications.html}
\end{center}
\vskip 1cm

\begin{abstract}
  In this paper, we describe Bohr-Sommerfeld rules for semi-classical
  completely integrable systems with 2 degrees of freedom with non
  degenerate singularities (Morse-Bott singularities) under the
  assumption that the energy level of the first Hamiltonian is non
  singular.  The more singular case of {\it focus-focus }
  singularities is studied in \cite{san-focus} and \cite{san-these}.
  The case of 1 degree of freedom has been studied in \cite{colin-p3}.
  
  Our theory is applied to some famous examples: the geodesics of the
  ellipsoid, the $1:2$-resonance, and Schrödinger operators on the
  sphere $S^2$. A numerical test shows that the \semicla\ \BS\ rules
  match very accurately the ``purely quantum'' computations.
\end{abstract}

\vfill


\noindent {\bf Mathematical classification}: 34E20 - 34L25 - 81Q20 -
58F07 - 58C40 - 58C27

\section*{Introduction}
In this paper, we describe extensions of results of \cite{colin-p3} to
completely integrable semi-classical systems with 2 degrees of
freedom. If $\hat{H}_1,\hat{H}_2 $ are two commuting
$h$-pseudo-differential operators on a 2-manifold $X$, we introduce
the {\it momentum map} $F=(H_1,H_2):T^\star X \ra \R^2$ where $H_j$
($j=1,2$) is the principal symbol of $\hat{H}_j$; we assume that $F$
is a proper map.  If $o=(0,0)$ is not a critical value of $ H$, the
existence and the construction of solutions of the system $\hat{H}_ju=
O(h^\infty )$ has a long history and existence of solutions is well
known to depend on the {\it Bohr-Sommerfeld} rules involving action
integrals and Maslov indices of loops generating the homology of the
fibre of $F$ which is a 2-torus.

We will assume that $o$ is a critical value and more precisely that
the critical points are of Morse-Bott type. A very simple
classification of such points with the corresponding normal forms is
given in \cite{san-fn}.  The {\it focus-focus} case have been
described in \cite{san-focus}.  We will be interested in the case
where the set of critical points of $H$ is a 2-dimensional manifold
with a transversally hyperbolic (saddle) singularity: it means that
$0$ is not a critical value of $H_1$ and that $H_2$ restricted to
Poincaré sections of the flow of $H_1$ admits critical points of
saddle type.  The set of critical values of $F$ is then a
1-dimensional submanifold of $\R^2$.  The main result of our paper is
a description of {\it singular Bohr-Sommerfeld rules} in this
situation. These rules give necessary and sufficient conditions for
existence of solutions of the system $\hat{H}_j u=O(h^\infty ),~j=1,2$
and approximations of the solutions.  More precisely, we show the
existence near the singular fibre $\Lambda _o=F^{-1}(o)$ of a
Hamiltonian $H_p$ with periodic flow which allows to {\it reduce} the
classical and also the semi-classical study to the 1-dimensional case
on the reduced phase space.  This reduced phase space can be {\it
  singular} because the $S^1$ action induced by the flow of $H_p$ is
not principal in general; non trivial isotropy group isomorphic to
$\Z/2\Z$ may appear.

 We provide a description
of $\Lambda _o$: the topological type can be rather unusual like 
a Klein bottle!

The precise description of the quantisation  rules 
is given in theorem 2.7. There is one rule giving  a quantum number
associated to the periodic orbits of $H_p$
 and  another rule  given in terms
of the graph $G$ which is the quotient of $\Lambda _o$ by the
$S^1$-action associated with $H_p$.
With a viewpoint similar to \cite{san-focus} but somewhat different
from \cite{colin-p3}, we interpret these rules as a universal
\emph{regularisation} of the usual \BS\ rules for tori. This, in
addition, allows us to have a description of the joint spectrum
inside a fixed neighbourhood of $o$, but raises some technical
difficulties due to the possible non-connectedness of the fibres of
$F$, especially in the $\Cinf$ category.

At the end we describe 3 examples for which we provide explicit
calculation and numerical checking:
\begin{enumerate}
\item High energy limit for eigenvalues of the Laplacian on ellipsoids.
\item Semi-excited states for anharmonic oscillators  with a resonance
$2:1$.
\item High energy limit for the Schrödinger spectrum on the 2-sphere.
\end{enumerate}
For the last two examples, numerical computations of eigenvalues of
large matrices are compared with the eigenvalues obtained from the
singular Bohr-Sommerfeld rules. We observe a very good accuracy of the
results even for not very big quantum numbers.


\newcommand{\M}{\mathcal{M}}
\renewcommand{\P}{\mathcal{P}}
\renewcommand\H{\mathfrak{H}}
\section{Classics}

\label{sec:classical}
The goal of this section is to give a semi-global description of the 
Lagrangian fibres of a 2-degrees-of-freedom integrable system having 
non-degenerate rank-one singularities of hyperbolic type.

Let $(M,\omega)$ be a symplectic manifold of dimension 4, and let
$H_1$, $H_2$ be two Poisson commuting Hamiltonian functions in
$C^{\infty}(M,\R)$.  The corresponding momentum map will be denoted by
$F=(H_1,H_2)$; we shall always assume $F$ to be proper, $H_1$ to be
non-singular on the level set $\Lambda_o:= F^{-1}(o)$ (for some point
$o\in\R^2$).  Moreover, we assume that $\Lambda _o$ is connected and
that the critical points of $F$ on $\Lambda_o$ are transversally
non-degenerate, in the following sense:

{\it for any  Poincar{\'e} section $\Sigma $
  of ${\cal X}_1$, the restriction
  of $H_2$ to $\Sigma $ is a Morse function.}

The main results of this section are
\begin{itemize}
\item  The description of the topology of the fibre $\Lambda_o$ 
  (section \ref{ss:lambdao}).

\item The construction of partial action-angle coordinates in a full
  neighbourhood of $\Lambda_o$; we show in particular that there exists
  a Hamiltonian $H_p$ defined in some neighbourhood of $\Lambda_o$
  that Poisson commutes with $H_j,~j=1,2,$ and all orbits of which are
  periodic (theorem \ref{theo:periodicham}).  Moreover, up to finite
  covering, a neighbourhood of $\Lambda_o$ is symplectomorphic to a
  product of $T^*S^1$ by a ``global'' Poincaré section.

\item The construction of normal forms for the system near each
  connected component of the critical set of $F$ (theorem
  \ref{theo:fn-classique}).
\end{itemize}

The core of these results can be obtained from Fomenko's
classification \cite{fomenko}, and the first two points are actually
due to Nguyên Tiên Dung \cite{zung-I}.  However, we felt that these
results deserved a detailed description with independent arguments.

\subsection{Notations}

For convenience of the reader we group here some of the notations: our
integrable system is given by the momentum map
$F=(H_1,H_2):M\rightarrow \R^2$.  Because $F$ is proper, it is a
momentum map for a Hamiltonian action of $\RM^2$ on $M$.
~$o=(a_o,b_o)$ is a critical value of $F$, $U$ small disk around $o$
in $\R^2$.  The smooth energy level is $S_o=H_1^{-1}(a_o)$.  $\Omega $
is $F^{-1}(U)$; $\Lambda _o= F^{-1}(o)$.  $\gamma =\cup _{i=1}^N
\gamma _i $ is the critical set of $F$ in $\Lambda _o$.  $\Gamma =\cup
\Gamma _i$ is the critical set in $\Omega$.  $\Lambda _{o}\setminus
\gamma =\cup \Lambda_{\{i,j\}}^k$ where $\Lambda_{\{i,j\}}^k$ is a
smooth Lagrangian cylinder whose closure is $\Lambda_{\{i,j\}}^k\cup
\gamma_i \cup \gamma_j$.  For any Hamiltonian $H_{something}$ we will
denote by ${\cal X}_{something}$ the associated Hamiltonian vector
field.  $G$ is a graph with $N$ vertices associated to $\Lambda _o$.
The Hessian of $H_2$ restricted to $\Sigma$ is
$\mathcal{H}_\Sigma(H_2)$.  The absolute value of its determinant
(with respect to the density induced by the symplectic for on
$\Sigma$) is independent on $\Sigma$ and denoted by
$|\mathcal{H}_\Sigma(H_2)|$.

\subsection{Topology of $\Lambda _o$}\label{ss:lambdao}

\begin{prop}
  The critical set $\gamma $ of $F$ in $S_o$ is a compact submanifold 
  of dimension 1 of $\Lambda _o$ which is a finite union of disjoint 
  periodic orbits $\gamma _i~(i=1,\cdots ,N)$ of ${\cal X}_1$.  The 
  $\gamma _i$ admits {\rm orbit-cylinders} $\Gamma _i$ which consists 
  of $\gamma _i (a)$, $a$ close to $a_o$.
\end{prop}
\begin{demo} Locally the reduced manifold of $S_o:=\{H_1=a_o\}$ is
  symplectomorphic to a Poincar\'e section $\Sigma $ of ${\cal X}_1$.
  Because $\ham{1}$ is transversal to $\Sigma$ in $S_o$, a
  neighbourhood $V$ of $\Sigma$ in $S_o$ is diffeomorphic to
  $\Sigma\times I$, for a small interval $I$, in such a way that the
  trajectories of $\ham{1}$ are of the form $\{\sigma\}\times I$,
  $\sigma\in\Sigma$.  Because $H_2$ is constant under the
  $\ham{1}$-flow, $H_2(\{\sigma\}\times I)=(H_2)_{\restr
    \Sigma}(\sigma)$.  By hypothesis, $H_2$ restricted to $\Sigma$
  admits a non degenerated critical point at $\sigma:=\gamma \cap
  \Sigma$.  Since the rank of $F$ if invariant by the flow of
  $\ham{1}$, we must have $\gamma\cap V=\{\sigma\}\times I$, which
  says that $\gamma$ is a smooth manifold of dimension 1.  Because, by
  definition, $\gamma$ is closed in the compact set $\Lambda_o$ and
  hence is compact, it must be a union of circles $\gamma=\bigcup
  \gamma_i$.  Then each of these circles must be an $\ham{1}$-orbit,
  and only a finite number $N$ of such critical circles arises due to
  the properness of $F$.

  The non-degeneracy of $(H_2)_{\restr \Sigma}$ implies that the 
  isolated critical point $\sigma=\sigma(a_o)$ depends smoothly on $a$ 
  close to $a_o$.  Therefore the above description extends to any leaf 
  of the foliation $H_1=a$, $a$ close to $a_o$, yielding a smooth family 
  of circles $\gamma_i(a)$.  Since a small neighbourhood of $S_o$ in $M$ 
  is diffeomorphic to $S_o\times\RM$ such that $a$ is a coordinate for 
  the second factor, the union $\Gamma_i=\bigcup_{a \text{ close to } 
    a_o}\gamma_i(a)$ is diffeomorphic to a cylinder $S^1\times(\RM,a_o)$.
\end{demo}
\begin{figure}[htbp]
  \centering
  \input{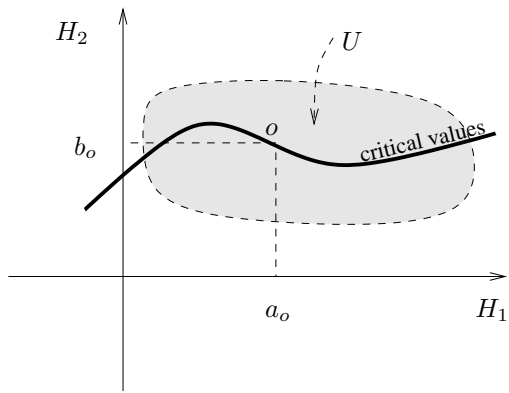}
  \caption{Image of the momentum map. The set of critical values near 
    $o$ is a smooth curve parameterised by $H_1$. If the critical point
    is a saddle point, the regular values are on both sides of the curve
    of critical values. If it is a maximum or a minimum, only one side
    is occupied.}
  \label{fig:image}
\end{figure}

If the critical point $\sigma$ is a local maximum or minimum of $H_2$,
$\Lambda_o$ reduces to one elliptic periodic orbit, a situation which
was studied a long time ago by several people (see
\cite{colinII}).  We will from now assume that {\it this
  critical point is a saddle point}.  Then one can have several critical
circles in $\gamma$; we show now how they can be connected to each
other inside $\Lambda_o$.

\begin{prop} 
  \label{prop:cylinder}
  $\Lambda_o\setminus\gamma $ is a union of $\R^2$-orbits that are 
  cylinders.  Each cylinder contains in its closure 1 or 2 
  $\gamma_i$'s.  We will denote by $\Lambda_{\{i,j\}}=\bigcup_k 
  \{\Lambda_{\{i,j\}}^k\}$ the set of cylinders that connect $\gamma_i$ 
  and $\gamma_j$.
\end{prop}
\begin{demo} Since $\gamma$ is $\RM^2$-invariant,
  $\Lambda_o\setminus\gamma$ is a union of orbits, on all of which the
  action is non-singular.  These orbits are therefore 2-dimensional
  quotients of $\R^2$, hence tori, cylinders or planes.  Because
  $\Lambda_o$ is connected and contains singular points, any orbit
  contains critical points in its closure, which excludes tori.  The
  local structure near each $\gamma_i$ will show independently the
  existence of periodic sub-orbits, leading to the non-triviality of the
  stabilisers of the $\RM^2$-action, which excludes planes.  Hence all
  orbits are cylinders.

  Now, the closure of such a cylinder consists of singular orbits: there
  are 1 or 2 of them.
\end{demo}

\begin{defi} We define the graph $G$ as follows: $G$ has $N$ 
  vertices (where $N$ is the number of critical circles in
  $\Lambda_o$), and there are exactly $|\Lambda_{\{i,j\}}|$ different
  edges connecting the vertices $i$ and $j$.
\end{defi}

Let $\Sigma$ be a Poincaré section at a point $m\in\gamma_i$.  By
hypothesis, $\Sigma\cap\Lambda_o$ is diffeomorphic to a ``hyperbolic
cross'', the union the local stable and unstable manifolds $W^\pm(m)$
for the flow of $\ham{\mbox{$H_2$}_{\restr \Sigma}}$.  Let $\Omega$ be
  a small neighbourhood of $\gamma_i$, and define $W^\pm(\gamma_i)$ as
  the union of the connected components of
  $(\Lambda_o\setminus\gamma_i)\cap\Omega$ intersecting $W^\pm(m)$.
  These manifolds do not depend on the choice of $m$.

\begin{prop}
  \label{prop:un-stable}
  Either $W^+(\gamma_i)$ and $W^-(\gamma_i)$ are diffeomorphic to the
  disjoint union of 2 cylinders or both are diffeomorphic to 1
  cylinder.  In the first case the vertex $i$ of $G$ has degree $4$
  while in the second it has degree $2$.
\end{prop}
\begin{figure}[htbp]
  \centering
  \input{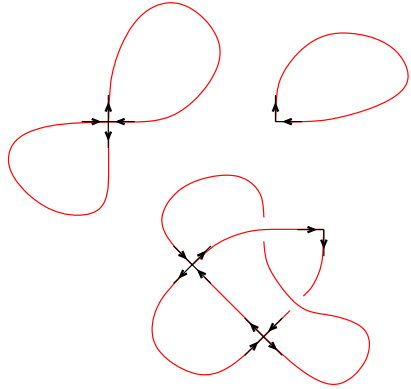}
  \caption{The vertices of $G$ are of degree 2 or 4.}
  \label{fig:graph}
\end{figure}

\begin{demo} The 2-manifold $\tilde{W}^+(\gamma_i) = W^+(\gamma_i)
  \cup \gamma _i$ is a bundle on $\gamma _i$ whose fibre is an interval.
  There are exactly 2 possibilities up to diffeomorphism: the trivial
  and the Moebius bundle.  In the first case, removing $\gamma _i$ gives
  2 cylinders, while in the second it gives only 1 cylinder.  Both
  bundle are isomorphic because the sum of their tangent bundle on
  $\gamma _i$ is a $\R^2-$bundle that is orientable as a symplectic
  bundle.
\end{demo}
\begin{defi} In the first case, $\gamma _i$ is called {\rm direct},
  in the second case, it is called {\rm reverse}.
\end{defi}
\begin{figure}[hbtp]
  \begin{center}
    \leavevmode
    \input{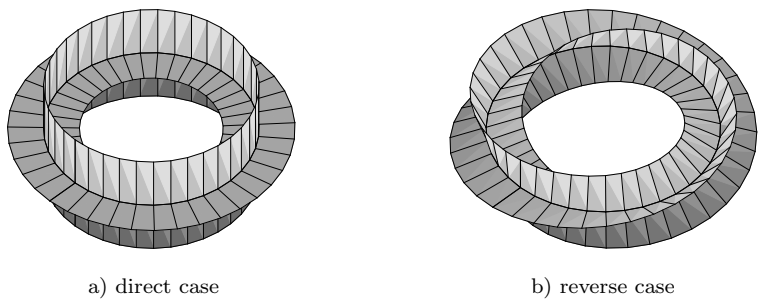}
    \caption{The neighbourhood of a critical circle}
    \label{fig:un-stables}
  \end{center}
\end{figure}

\subsection{The classical commutant}
\begin{lemm}
  \label{lemm:commutant}
  For any smooth function $K$ commuting with $H_1$ and $H_2$, the 
  vector field $\ham{K}$ can be 
  uniquely written (in a neighbourhood of $\Lambda_0$)
  \[ \ham{K} = a\ham{1}+b\ham{2}, \]
  for smooth functions $a$, $b$ commuting with $H_1$ and $H_2$.
\end{lemm}
\begin{demo}
  Near any nonsingular point of $F$, apply the Darboux-Carathéodory 
  theorem which makes $H_1=\xi$ and $H_2=\eta$ in a local symplectic 
  coordinate chart $(x,y,\xi,\eta)$.  Then $K=K(\xi,\eta)$, and the 
  result follows.  
  
  Near a critical point in some $\gamma_i$, we use Theorem
  \ref{theo:fn-local} below that reduces the situation to
  $H_1=\xi+a_o$ and $H_2=\Phi(\xi,y\eta)$, and hence to ($H_1=\xi$,
  $H_2=y\eta$), for which we use the result of \cite[Lemme 2.6]{san-fn}
  (\cite[Lemme 2.2.7]{san-these}).
\end{demo}

\begin{theo}
  \label{theo:fn-local}
  Around any point $m\in\gamma_i$, there exists a canonical chart 
  $(x,y,\xi,\eta)$ in which
  \[\left\{ 
    \begin{array}{rcl}
      H_1-a_o & = & \xi\\
      H_2 & = & \Phi (\xi, y\eta),
    \end{array}
  \right.
  \]
  for some smooth function $\Phi$ defined near the origin, with 
  \begin{equation}
    \label{equ:sign-Phi}
    \partial_2\Phi(0,0) > 0.
  \end{equation}
\end{theo}
\begin{demo}
  First construct a Darboux-Carathéodory chart for $H_1$, i.e. a
  system of symplectic coordinates $(x,y,\xi,\eta)$, with canonical
  form $\omega_0$, in which $\xi=H_1-a_o$.  In these coordinates, the
  plane $\{0\}\times\{0\}\times\RM^2$ is a Poincaré section for $H_1$,
  hence by hypothesis $(y,\eta)\fleche H_2(\xi,y,\eta)$ has, for
  each small $\xi$, a non-degenerate saddle point.  The application of
  the Isochore Morse Lemma \cite{colin-vey} with parameter $\xi$
  yields a local symplectomorphism $\phi_\xi(y,\eta)$ of the
  $(y,\eta)$-space, depending smoothly on $\xi$, such that
  \begin{equation}
    H_2(x,\xi,\phi_\xi(y,\eta)) = \Phi(\xi,y\eta), 
    \label{equ:fn1}
  \end{equation}
  for a function $\Phi$ with $\partial_2\Phi(0,0) \neq 0$.  Applying
  the canonical transformation $(y,\eta)\fleche(-\eta,y)$ if
  necessary, one can assume that $\Phi$ satisfies the condition
  (\ref{equ:sign-Phi}) of the theorem.  The map $\phi :
  (x,\xi,y,\eta)\mapsto(x,\xi,\phi_\xi(y,\eta))$ is a local
  diffeomorphism but need not be symplectic.  A modification of $x$
  shall solve the problem.

  The 2-form $\omega_1:= \phi^*(\omega_0)$ splits as follows~:
  \[ \omega_1 = \omega_0 + d\xi\wedge\beta, \] for some 1-form
  $\beta=\beta(\xi,y,\eta)$.  From $d\omega_j=0$, one gets $d\xi\wedge
  d\beta=0$, which means $d_{(y,\eta)}\beta = 0$.  Hence there exists a
  smooth function $f(\xi,y,\eta)$ such that $d_{(y,\eta)}f=\beta$.  Let
  $\phi_1$ be the diffeomorphism~:
  \begin{equation}
    \phi_1(x,\xi,y,\eta) := (x-f(\xi,y,\eta),\xi,y,\eta), 
    \label{equ:fn2}
  \end{equation}
  so that 
  \begin{equation}
    \phi_1^*\omega_1=\phi_1^*\omega_0 + \phi_1^*(d\xi\wedge\beta) =
    \phi_1^*\omega_0 + d\xi\wedge\beta
    \label{equ:fn3}
  \end{equation} 
  \begin{equation} 
    = \omega_0 - d\xi\wedge df +
    d\xi\wedge\beta = \omega_0.
    \label{equ:fn4}
  \end{equation}
  Thus $\phi\circ\phi_1$ is symplectic, and because it does not change 
  $H_1=\xi-a_o$, it answers the question. 
\end{demo}
%
%

\subsection{A periodic Hamiltonian $H_p$ in $\Omega$}
\label{ss:periodicham}

\begin{theo}[\cite{zung-I}]
  \label{theo:periodicham}
  There exists a unique (up to additive constant) Hamiltonian $H_p$ in
  $\Omega $ that fulfils the following conditions:
  \begin{enumerate}
  \item $H_p$ Poisson commutes with $H_1$ and $H_2$.
    
  \item The flow of $H_p$ is $2\pi $ periodic with minimal period
    $2\pi $ outside $\Gamma$.
    
  \item On $\Gamma_i$,
    \begin{equation}
      \label{equ:transv-Hp}
      {\cal X}_p=\alpha_i(H_1) {\cal X}_1 ,
    \end{equation}
    with $\alpha_i$ a positive function.
  \end{enumerate}
  Then the flow of $H_p$ is $2\pi$-periodic on $\Gamma_i$ if the
  vertex $\{i\}$ is direct, and $\pi$-periodic in the reverse case.
\end{theo}

\begin{lemm}
  \label{lemm:critic}
  Near each $\gamma _i$, there exists a unique (up to some additive 
  constant if $a_o\ne 0$) Hamiltonian $H_{2,i}=H_2-\lambda _i(H_1)H_1$ 
  which is critical on $\Gamma _i$.
\end{lemm}
\begin{demo} 
  If $dH_2=f_i(H_1)dH_1$ on $\Gamma _i$, one gets the differential 
  equation \[ t\lambda_i '+\lambda_i =f_i(t) \] which always admits a 
  local solution.
\end{demo}

\begin{defi}
  A class of path $z\rightarrow [\gamma _z]$ (where $[~]$ means a
  homotopy or homology class) is called {\rm smooth} if there exists
  locally a representative that smoothly depends on $z$.
\end{defi}

\begin{lemm}
  \label{lemm:class}
  Let us denote by $L_z$ the orbit of $z$ by the $\R^2-$action.  There
  exists a unique mapping $z\rightarrow [\gamma _z]$ that is a smooth
  map from $\Omega $ into $H_1(L_z,\Z)$ and such that, if $z\in \gamma
  _i$, $[\gamma _z]=\nu [\gamma _i]$ where $\gamma _i$ is oriented
  according ${\cal X}_1$ and $\nu=1$ (resp.  $2$) if $\gamma _i$ is
  direct (resp.  reverse).
\end{lemm} 
\begin{demo} 
  a) The main point for this proof is to construct such a smooth family 
  $[\gamma_z]$ near $\gamma_i$.  By the Morse-Bott lemma (see Appendix 
  \ref{app:morse-bott}) applied to the Hamiltonian $H_{2,i}$ of Lemma 
  \ref{lemm:critic}, there are a coordinate $x_1$ on the circle $\gamma 
  _i=(\R/\Z)$ and a fibre bundle $F_\pm$ of dimension 2 on $\gamma_i$ 
  defined as quotient of the trivial bundle on $\R$ by identifying 
  $(x_1,w)\in \gamma_i \times \R^2$ with $(x_1+1,\pm w)$ such that 
  $H_{2,i}=x_3x_4$.  We can assume that ${\cal X}_1$ has the same 
  orientation as $\partial / \partial x_1$ on $\gamma _i$.

  The above Morse-Bott lemma can be applied with the parameter $H_1$.  
  Therefore, one gets coordinates $(x_1,x_2,x_3,x_4)$ on a full 
  neighbourhood of $\gamma_i$ in $M$ by letting $x_2=H_1$.

  We choose then $\gamma _z$ to be the path given in these coordinates 
  by $t\rightarrow (t,x_2,x_3,x_4)$ with $t\in \R/\Z$ (direct case) and 
  $t\in \R /2\Z$ (reverse case).  These path are drawn on the leaf 
  $(H_1,H_{2,i})=const$, and hence on a leaf of $F=(H_1,H_2)$.  The last 
  assumption of lemma \ref{lemm:class} follows from proposition 
  \ref{prop:un-stable}.

  b) Far from $\gamma_i$, we construct $\gamma_z $ on $\Lambda_o$ by
  taking the cylinder equator with the orientation which in the affine
  structure of $\Lambda_{\{i,j\}}^k$ is given by projecting ${\cal X }_1$. 
  Then we extend to the nearby Lagrangian leaves by local triviality of
  the foliation by orbits.  Because these $\gamma_z$ are homotopic to
  the ones constructed in a) and lying on the same Lagrangian leaf, it
  is easy to realise this homotopy as an isotopy, thus yielding a smooth
  family of loops in $\Omega$.
\end{demo}

c) It remains now to use these $\gamma_z$ in order to define the action
variable.
\begin{lemm} 
  \label{lemm:alpha}
  The symplectic form $\omega $ is exact in $\Omega $, i.e. there
  exists $\alpha $, 1-form in $\Omega$ with $d\alpha =\omega$.
\end{lemm}
\begin{demo} 
  The set $\Lambda_o$ is Lagrangian and any 2-cycle can be deformed
  inside $\Lambda_o$.
\end{demo}

\begin{demo}[(of Theorem \ref{theo:periodicham})]
  Put $H_p(z)=(1/2\pi)\int _{\gamma _z}\alpha$.
  $H_p$ is smooth and commutes with $H_1$ and $H_2$ (because it is
  constant on $\R^2-$orbits).  Moreover on $\R^2-$orbits that are tori
  the orbits of ${\cal X}_p$ are $2\pi -$periodic with orbits
  homotopic to $\gamma _z$ (by usual action-angle coordinates).  On
  the $\gamma _i$'s, the period is $2\pi $ in the direct case and $\pi
  $ in the reverse case.
  
  The affine structure on the Lagrangian cylinders
  $\Lambda_{\{i,j\}}^k$ and the condition (3.)  implies the uniqueness
  of $\ham{p}$ on $\Lambda_o$.  Now suppose $H_p'$ is another
  Hamiltonian with the same properties, and let $z$ be a point in
  $\Lambda_o\setminus\gamma$.  Since the orbits under $\ham{p}$ and
  $\ham{p}'$ of $z$ are equal, the orbits of points near $z$ in a same
  level set of $F$ (different from $\Lambda_o$) are homotopic.  But
  these level sets are Liouville tori for which we know that $\ham{p}$
  and $\ham{p}'$ must be equal.
\end{demo}

\begin{rema}
  Step a) of the proof does not use the nature of $\Lambda_{\{i,j\}}^k$. 
  Therefore, a)+c) gives a Hamiltonian $H_p$ verifying the conditions 
  of Theorem \ref{theo:periodicham} but in a neighbourhood of 
  $\gamma_i$ only. This suffices to prove that the $\RM^2$-action has 
  non-trivial stabilisers, whence $\Lambda_{\{i,j\}}^k$ must be cylinders, 
  thus finishing the proof of proposition \ref{prop:cylinder}.
\end{rema}

\subsection{$S^1$ reduction}
\label{sec:reduction}
The flow of $\ham{p}$ yields a locally free Hamiltonian action of
$S^1$ on $\Omega$, which is free outside $\Gamma$.  Let
$c_o=H_p(m)=H_p(\Lambda_o)$. We denote by $W$
the reduced space 
\[ W = H_p^{-1}(c_o)\cap\Omega/S^1.  
\]
$W$ is a symplectic orbifold (see eg.  \cite{cushman-book}).  It is a
smooth manifold if and only if the action is free, that is if and only
if no vertex of reverse type are present in the graph $G$.  Otherwise,
it has singularities at the critical orbits $\gamma_i$.  Since these
critical orbits come in families depending on the value of $H_1$,
yielding local orbit cylinders, and because $dH_p(m)=\lambda dH_1(m)$,
$\lambda\neq 0$, only one orbit of each local cylinder meets
$H_p^{-1}(c_o)$, ensuring that the critical orbits give isolated
singularities in $W$. 

Let
\begin{equation}
  H_q=-b(H_1-a_o)+a(H_2-b_o),
  \label{equ:Hq}
\end{equation}
where $a>0$ and $b$ are the real constants such that, on $\Lambda_o$, 
$\ham{p}=a\ham{1}+b\ham{2}$ (cf.  Lemma \ref{lemm:commutant}).  Then 
$\Lambda_o=H_p^{-1}(c_o)\cap H_q^{-1}(0)$.  (This still holds of course 
for a generic choice of $(a,b)$.)  Since $H_q$ is $S^1$-invariant, it 
defines a smooth Hamiltonian function $\tilde{H}_q$ on $W$.  The graph 
$G$ can be viewed as the quotient of $\Lambda_o$ by $S^1$, and thus is 
identified to the level set $\tilde{H}_q^{-1}(0)$.

\begin{prop}
  \label{prop:poincare-global}
  If the $S^1$-action is free (\emph{i.e.} all vertices of $G$ of
  degree 4), then a neighbourhood $\Omega$ of $\Lambda_o$ in $M$ is
  diffeomorphic to the direct product $S^1\times\RM\times W$ (hence
  $\Lambda_o$ is diffeomorphic to the direct product $S^1\times G$ --
  these diffeomorphisms are equivariant with respect to the natural
  action of $S^1$ on itself).
\end{prop}
\begin{rema}
  In this case, $W$ can be regarded as a ``global'' Poincaré section 
  for $\ham{p}$.
\end{rema}
\begin{demo} 
  We choose now $\Omega$ to be of the form $\Omega_o\times I$ where
  $I$ is some open interval around $0$, $\Omega_o$ is a small
  invariant neighbourhood of $\Lambda_o$ in $H_p^{-1}(c_o)$, and
  $H_p(\Omega_o\times\{\xi\})-c_o=\xi$.  If the action is free,
  $\Omega_o$ is a principal $S^1$-bundle over $W$.  It is
  topologically classified by its holonomy class in
  $H^1(W,S^1_{\mathfrak{d}})$, where $S^1_{\mathfrak{d}}$ is the sheaf
  of germs of smooth functions on $W$ with values in $S^1$ (see
  \cite{hirzebruch}).  Using the short exact sequence $0\fleche \ZM
  \fleche \RM \fleche S^1 \fleche 0$, and the fact that the sheaf
  $\RM_{\mathfrak{d}}$ is fine, one gets an isomorphism
  \[
  H^1(W,S^1_{\mathfrak{d}}) \simeq H^2(W,\ZM),
  \] 
  yielding the so-called Chern class of the bundle.  But $W$ retracts
  onto $G$ and $G$ is 1-dimensional, so $H^2(W,\ZM)=0$, and $\Omega_o$
  is a trivial bundle.
\end{demo}

\begin{theo}[\cite{zung-I}]
  \label{theo:action-angle}
  If the $S^1$-action is free, then $(\Omega,\omega)$ is
  symplectomorphic to a neighbourhood of $S^1\times\{0\}\times W$ in
  \[ 
  (T^*S^1 \times W,\quad d\xi\wedge dx + \pi^*\omega_W), 
  \] 
  with $H_p-c_o=\xi$.  Here $\pi$ is the projection onto $W$ and $\omega_W$
  is the symplectic form of $W$.
\end{theo}
\begin{demo}
  First apply proposition \ref{prop:poincare-global} and let $\xi=H_p-c_o$
  be a coordinate for the $\RM$ factor.  Then choose the conjugate
  angle variable $x$ (pick up some coordinate $\theta$ in $S^1$, an
  origin $\theta_0$, and let $x(\theta)$ be the time required to go
  from $\theta_0$ to $\theta$ under the Hamiltonian action of $\xi$),
  so that $\{\xi,x\}=1$.  Because $\omega$ is $S^1$-invariant, it does
  not depend on $x$; using the equivariant Darboux-Weinstein theorem
  \cite{weinstein}, one can assume that $\omega=\omega_{\restr
    \xi=0}$, and so does not depend on $\xi$ either.
  
  Because for any $\xi$, $\ham{p}$ is $\omega$-orthogonal to
  $S^1\times\{\xi\}\times W$, one easily checks by taking local
  coordinates on $W$ that
  \[ 
  \omega = d\xi\wedge dx + d\xi\wedge\pi^*\beta + \pi^*\omega_W, 
  \] 
  where $\beta$ is a one-form on $W$, and $\omega_W$ a non-degenerate
  2-form on $W$.  The closedness of $\omega$ (and its independence on
  $\xi$) implies $d\omega_W=0$ and $d\beta=0$, the latter yielding
  $d(\xi\pi^*\beta)=d\xi\wedge\pi^*\beta$.  Let us now apply Moser's
  path method to get rid of this term.  We let
  \[ 
  \omega_t := d\xi\wedge dx + \pi^*\omega_W + td(\xi\pi^*\beta),
  \]
  and wish to construct an isotopy $\phy_t$ of diffeomorphisms of
  $\Omega$ such that $\phy_t^*\omega_t = \omega_0$.  $\phy_t$ is then
  given as the flow of the vector field $X_t$ defined by
  $i_{X_t}\omega_t + \xi\pi^*\beta=0$.  It is easy to check that
  $\omega_t$ is non-degenerate for all $t$ so that $X_t$ is uniquely
  defined.  Moreover, because of its defining equation, $X_t$ is of
  the form $\xi\iota_*Y_t$ ($\iota$ is the inclusion $W\inject
  T^*S^1\times W$), where $Y_t$ is a vector field on $W$ satisfying
  $i_{Y_t}(\omega_t)_{\restr W} + \beta=0$.  Therefore $\phy_t$ is of
  the form
  \[ 
  (x,\xi,w) \fleche (x,\xi,\phi_{\xi t}(w)),
  \]
  (where $\phi_t$ is the flow of $Y_t$) and hence preserves $x$ and
  $\xi$.  If $Y_t$ can be integrated up to the time $t_0>0$, then
  $X_t$ can be integrated up to the time $1$ for $\xi\leq t_0$.  The
  diffeomorphism $\phy_1$ then answers the question.
\end{demo}

\begin{rema}
  The formula $\omega = d\xi\wedge dx + \pi^*\omega_W$ ensures that 
  $\omega_W$ is the natural symplectic form on $W$ obtained by the 
  reduction process.
\end{rema}

\begin{theo}
  \label{theo:revetement}
  In the general case, there exists a smooth double covering
  $\Omega^*$ of $\Omega$ in which the action is free.  The reduced
  manifold $W^*$ is a covering of $W$ that is ramified of degree 2 at
  the critical orbits $\gamma_i$.
\end{theo}
\begin{demo}
  Choose $\Omega$ to be a relatively compact invariant neighbourhood
  of $\Lambda_o$ in $M$.  Then $H_p^{-1}(c_o)\cap\Omega$ has a smooth
  non-empty invariant boundary, and because this boundary does not
  meet any critical orbit, the closure
  $\overline{W}=\overline{H_p^{-1}(c_o)\cap\Omega}/S^1$ is a relatively
  compact surface with a non-empty smooth boundary.
  Let $p_i, i=1\dots \ell$ be the images under reduction of the 
  critical circles $\gamma_i$ and let $\check{W}$ be the surface 
  $\overline{W}$ after removal of small disks $D_i$ around each $p_i$.  
  It is still a smooth surface with boundary, whose fundamental group 
  is free \cite{ahlfors-sario}, and generated by some 
  $\mu_1,\dots,\mu_k,\delta_1,\dots,\delta_\ell$, with 
  $\delta_i=\partial D_i$.  Let $\mathcal{D}\subset \pi_1(\check{W})$ 
  be the free subgroup generated by 
  $\mu_1,\dots,\mu_k,\delta_1^2,\dots,\delta_\ell^2$, and let 
  $\check{W}^*$ be the corresponding smooth covering of $\check{W}$.  
  Gluing back the disks $D_i$ defines a covering space $W^*$ of $W$ 
  that is ramified of degree two at each $p_i$.  Define now 
  $\Omega^*\subset \Omega\times W^*$ such that the following diagram 
  commutes:
  \[ 
  \begin{array}{ccc}
    \Omega^* & \flechedroite\ & W^*  \\
    \flechebas\ &  & \flechebas\  \\
    \Omega & \flechedroite\ & W
  \end{array}
  \]
  Since the local structure near $\gamma_i$ (see Theorem 
  \ref{theo:fn-classique} in the next section) gives a model for the
  covering $\Omega^*\fleche\Omega$, $\Omega^*$ is 
  naturally endowed with a smooth structure compatible with that of 
  $\Omega$. The lifts of $\gamma_i$ are critical circles in $\Omega^*$ 
  that become of direct type. As a result, all critical circles in 
  $\Omega^*$ are of direct type, which means that the lifted $S^1$ 
  action is free.
\end{demo}

\subsection{Normal forms near $\gamma_i$}
\label{sec:fn}
Choose any Hamiltonian function $H_q$ near $\gamma_i$ that commutes
with $H_1$ and $H_2$ and such that, \emph{for any} $o$ in the local
curve of critical values of $F$,
\begin{itemize}
\item $H_q(\Lambda_o)=0$;
  
\item on $\Lambda_o\setminus\gamma_i$, $\ham{p}$ and $\ham{q}$ are
  linearly independent, and the automorphism:
  $(\ham{1},\ham{2})\fleche (\ham{p},\ham{q})$ is orientation
  preserving.
\end{itemize}
For instance, the previously defined $H_q$ (eq.  (\ref{equ:Hq})) is a
good choice, which is independent on $i$, but any generic linear
combination $H_q=\alpha (H_1-a_o)+ \beta (H_2-b_o)$, $\beta>0$ would
also do.

Then the following theorem states a simple normal form near $\gamma_i$
for the new system $(H_p,H_q)$.  It will be the main tool for the
semi-classical analysis, for it reduces the situation to the case of a
cotangent bundle.

\begin{theo} 
  \label{theo:fn-classique}
  There exists coordinates $(x,y)$ on $R=\R \times ]A,B[$ that give
  local coordinates on $\tilde{W}^+(\gamma_i)$ by taking the quotient
  by
  \[
  (x,y)\rightarrow (x+2\pi,y)
  \]
  in the direct case and
  \[ 
  (x,y)\rightarrow (x+\pi,-y)
  \]
  in the reverse case and a canonical diffeomorphism of a
  neighbourhood of $\gamma _i$ into a neighbourhood of the
  ``$\xi=c_o$-section'' of $T^*(\tilde{W}^+(\gamma_i))$ (recall that
  $c_o=H_p(\Lambda_o)$) such that with respect to canonical
  coordinates we get
  \[ \left\{
    \begin{array}{rcl}
      H_p & = & \xi \\
      H_q & = &\Phi (\xi, y\eta),
    \end{array}
  \right.
  \]
  for some smooth function $\Phi$ defined near the origin, with
  $\partial_2\Phi(0,0) > 0$.
\end{theo}
\begin{demo} 
  1. We first wish to prove that the restriction of $H_q$ to the
  locally reduced manifold $W:=H_p^{-1}(c_o)/\exp{t\ham{p}}$ has a
  non-degenerate saddle point. Let $m\in\gamma_i$ and let
  $(s,\sigma,u,v)$ be local coordinates near $m$ such that
  $\sigma=H_1-a_o$ and the flow of $\ham{1}$ is just translation on
  the $s$ variable. Then because $dH_p(m)=\lambda dH_1(m)$ with
  $\lambda\neq 0$, the map $(s,\sigma,u,v)\fleche(s,H_p-c_o,u,v)$ is a
  local diffeomorphism of $(\RM^4,0)$ that sends $(s,0,0,0)$ to
  $(s,0,0,0)$. Therefore one can take $(u,v)$ as local coordinates for
  $W$, and we wish to prove that $(H_q)_{\restr W}(u,v)$ has a
  non-degenerate saddle point at $(0,0)$.

  Let $\M=
  \left(
    \begin{array}{cc}
      a & b  \\
      c & d
    \end{array}
  \right)$ 
  be the matrix of smooth functions such that
  \[ (dH_1,dH_2) = \M\cdot(dH_p,dH_q). \]
  Since $H_2$ is of the form $H_2=F(H_1,u,v)$, one has
  \[ dH_2 = K.dH_1 + A, \]
  where $K=\partial_1 F$ and $A=A(H_1,u,v)$ is a one-form on
  $\{0\}\times\{0\}\times\RM^2$ depending smoothly on $H_1$, that
  vanishes at $m$ and whose differential at $m$ is a non-degenerate
  quadratic form of hyperbolic type.  Using $\M$ one gets
  \[ (c-Ka)dH_p = (Kb-d)dH_q + A. \]
  The claim is that $(Kb-d)(m)\neq 0$.  Indeed, $(c-Ka)$ and $(Kb-d)$
  cannot simultaneously vanish because $\M$ is invertible.  It
  suffices then to see that $A(m)=0$ whereas by hypothesis
  $dH_p(m)\neq 0$.
  
  Therefore, we have, on $TW\subset\ker dH_p$,
  \[ d(H_q)_{\restr TW} = -\frac{1}{Kb-d}A, \]
  which implies that $d(H_q)_{\restr TW}$ possesses, as $A$ does, a 
  non-degenerate differential of hyperbolic type.
  
  2.  Now, we consider a neighbourhood of the whole critical circle
  $\gamma_i$ and use Weinstein's theorem with $S^1$ action to reduce to
  the cotangent space of $\tilde{W}^+(\gamma_i)$ with $(H_p-c_o=\xi,H_q=
  H_q(\xi,y,\eta))$.  In these coordinates (we consider first the
  direct case), $W$ is naturally identified with
  $\{0\}\times\{0\}\times\RM^2$, so the previous point shows that
  $(y,\eta)\fleche H_q(\xi,y,\eta)$ has, for each small $\xi$, a
  non-degenerate saddle point.  Then the proof goes exactly as that of
  Theorem~\ref{theo:fn-local}.  Equations (\ref{equ:fn1}) through
  (\ref{equ:fn4}) are valid if $H_1$ and $H_2$ are replaced by $H_p$
  and $H_q$, and $\omega_0$ is the canonical 2-form of
  $T^*(\tilde{W}^+(\gamma_i))$.
  
%

  In the reverse case, the proof is the same but we need the Isochore 
  Morse lemma for functions that are invariant by the involution $\sigma 
  (x)=-x $: it it then possible to choose the diffeomorphism $F$ 
  commuting with $\sigma$.  This fact follows easily from the proof 
  given in \cite{colin-vey}: in the {\it lemme principal} p.~283, we 
  choose $\eta$ such that $\sigma ^\star (\eta)=-\eta$.  It implies that 
  $X_t$ commutes with the involution.
\end{demo}

\begin{rema}
  \label{rema:orientation}
  We decide to give to the graph $G$ the orientation of the flow of 
  $H_q$ defined by (\ref{equ:Hq}).  Near a vertex $\gamma_j$, it is 
  also given by the flow of the normal form $y\eta$.
\end{rema}

\subsection{Appendix: Morse-Bott lemma}
\label{app:morse-bott}

\begin{defi}
  Let $f:X\ra \R$ be a smooth function.  A submanifold $W$ of $X$ is
  called a {\rm Morse-Bott critical manifold} if every point $w\in W$
  is a critical point of $f$ and if the restriction of $f"(w)$ to the
  normal bundle $T_wX/T_w W$ is non degenerate.
\end{defi}

Morse-Bott critical manifold arises in many situations especially
when $f$ is invariant by a Lie group action.  An extension of the
Morse lemma is available in that case.  In some situations, there is
global topological problem with the subbundles $N_\pm$ of the normal
bundle generated by eigenspaces of $f"$ associated with $>0$ (resp.
$<0$) eigenvalues.

\begin{lemm}[Morse-Bott lemma]
  Assume we have a Morse-Bott connected critical manifold $W$ for a
  function $f:X\rightarrow \R$.  Let $N$ be the normal bundle of $W$
  and $F$ the Hessian of $f$ which is a non degenerate quadratic form
  on $N$.  Then there exists a diffeomorphism of a neighbourhood of $W
  $ on a neighbourhood of the 0-section in $N$ which conjugates $f$
  to $F+c$.
\end{lemm} 
If $W$ is connected, a complete set of invariant of $f$, up to
smooth conjugacy near $W$, is given by
the pair $(N_+,N_-)$ of bundles on $W$ up to isomorphism.

\renewcommand{\H}{\hat{H}}
\renewcommand{\F}{\hat{F}}
\renewcommand{\L}{\mathfrak{L}}
\renewcommand{\C}{\mathfrak{C}}
\newcommand{\bgamma}{\bar{\gamma}}
\newcommand{\Hun}{\mathop{\mathrm{H}_1}}
\newcommand{\HUN}{\mathop{\mathrm{H}^1}}
\newcommand{\cHUN}{\mathop{\mathrm{\check H}^1}}
\newcommand{\hol}{\mathop{\textbf{\textup{hol}}}}
\newcommand{\ex}{\mathrm{e}}

\newcommand{\e}{\mathfrak{e}}
\newcommand{\K}{\hat{K}}
\renewcommand{\k}{\overline{k}}
\newcommand{\CC}{\mathcal{C}}
\renewcommand{\O}{\mathcal{O}}
\section{Semi-classics}

The aim of this section is to express the singular \BS\ quantisation 
rules for quantum integrable systems whose classical counterpart 
fulfils the hypothesis of the previous section.

Let $X$ be a 2-dimensional differential manifold, and let $\H_1(h)$,
$\H_2(h)$ be commuting $h$-\pdo s, with real principal symbols $H_1$,
$H_2$. Assume that the momentum map $F=(H_1,H_2)$ satisfies the
hypothesis of section \ref{sec:classical}.
In all of this section, the 1-form $\alpha$ of Lemma \ref{lemm:alpha}
is taken to be the canonical Liouville 1-form of the cotangent bundle
$T^*X$.  Then $H_p$ (Theorem \ref{theo:periodicham}) is uniquely
defined as the action integral with respect to $\alpha$.  For any $E$
in the image of $F$, the sub-principal form $\kappa_E$ is the closed
differentiable 1-form on $\Lambda_E:=F^{-1}(E)$ defined at its regular
points by $\kappa_E(\ham{j})=-r_j$, where $r_j$ is the sub-principal
symbol of $\H_j$.

\subsection{The \mi\ normal form}
We will prove here a \semicla\ analogue of Theorem
\ref{theo:fn-classique}, which was particularly fit for this purpose
since it reduced the situation to that of a cotangent bundle, for
which the usual pseudo-differential quantisation can be used.  A
\semicla\ analogue of Theorem \ref{theo:action-angle} should also be
interesting, but would involve symplectically reduced cotangent
bundles, for which Toeplitz quantisation is needed, a theory that we
don't want to enter here.

In this section a critical circle $\gamma_j\subset\Lambda_o$ is fixed. 
Theorem \ref{theo:fn-classique} identifies, via a symplectomorphism
$\psi$, a neighbourhood of $\gamma_j$ in $T^*X$ with a neighbourhood
of the zero section of a cotangent bundle of the form
$T^*(\RM^2/\sigma)$, where $\sigma(x,y)=(x+2\pi,y)$ in the direct
case, and $\sigma(x,y)=(x+\pi,-y)$ in the reverse case.

It is easy to check that Weyl quantisation satisfies, for a symbol
$a\in\Cinf_0(\RM^2)$~:
\[
\sigma^* Op_h^W(a) = Op_h^W(a\circ T^*\sigma)\sigma^*,
\]
where $\sigma^*$ is the adjoint operator $u\mapsto u\circ\sigma$, and
$T^*\sigma$ is the cotangent lift of $\sigma$.  Therefore, if
$a=a\circ T^*\sigma$, then $Op_h^W(a)$ acts on the space of functions
$u$ that are invariant under $\sigma$~: $u\circ\sigma=u$, which is the
space of functions defined on a cylinder in the direct case, and on
the M\oe bius strip in the reverse case.  In particular,
$Q_1(h)=Op_h^W(\xi)$ and $Q_2(h)=Op_h^W(y\eta)$ are well-defined
differential operators on $\RM^2/\sigma$~:
\begin{equation}
Q_1(h)= \frac{h}{i}\deriv{}{x}, \quad 
Q_2(h)= \frac{h}{i} \left(y\deriv{}{y}+\frac{1}{2}\right). 
\label{equ:opQ}
\end{equation} 

Let $\Psi^0$ be the algebra of operators of the form $Op_h^W(a(h))$ 
for classical symbols $a(h)$ on $T^*(\RM^2/\sigma)$, modulo those 
whose symbol is $\oh$.  For more details about the \mi\ calculus and 
$\oh$ remainders, see eg.  \cite{colin-p}, \cite{san-these}.  Before 
stating the result of this section, we introduce the following 
spaces~:
\begin{defi}
  \label{defi:commutant}
  The classical and \semicla\ commutants $\C_{cl}(\gamma_j)$ and
  $\C_h(\gamma_j)$ are defined as follows~:
  \[ 
  \C_{cl}(\gamma_j)=\{f\in\Cinf(T^*(\RM^2/\sigma)), \quad
  \{f,\xi\}=\{f,y\eta\}=0 \textrm{ near } \gamma_j\};
  \]
  \[ 
  \C_h(\gamma_j)=\{P(h)\in\Psi^0, \quad [P,Q_1] \textrm{ and } [P,Q_2]
  \textrm{ are } \oh \textrm{ near } \gamma_j\};
  \]
\end{defi}
(Recall that in our coordinates,
$\gamma_j=((\RM\times\{0\})/\sigma)\times\{c_o\}\times\{0\}$, where
$c_o=H_p(\Lambda_o)$.)  Because the symbols of $Q_1$ and $Q_2$ are
polynomials of degree $\leq 2$, the operators in $\C_h(\gamma_j)$ are
exactly the Weyl quantisations of symbols of the form $\sum h^ka_k$
with $a_k\in\C_{cl}(\gamma_j)$.
\begin{theo}[Microlocal normal form]
  \label{theo:fn-semicla}
  There exists an elliptic \fio\ $U(h)$ associated to the canonical
  transformation $\psi$ of Theorem \ref{theo:fn-classique}, an
  invertible $2\times 2$ matrix $\M(h)$ of \pdo s in $\C_h(\gamma_j)$,
  and complex-valued functions of $h$~: $\epsilon_1(h)$ and
  $\epsilon_2(h)$ admitting an asymptotic expansion in
  $\CM\formel{h}$:
\[
\epsilon_1(h)\sim \sum_{\ell=0}^\infty \epsilon_1^{(\ell)} \qquad
\epsilon_2(h)\sim \sum_{\ell=0}^\infty \epsilon_2^{(\ell)}
\]
 such that, microlocally near $\gamma_j$~:
  \begin{equation}
    \label{equ:fn}
    U^{-1}(\H_1-a_o,\H_2-b_o)U = \M.(Q_1-\epsilon_1,Q_2-h\epsilon_2) + \oh.
  \end{equation}
  If $\H_1$ and $\H_2$ are formally self-adjoint, then $U(h)$ can be
  chosen to be \mi ly unitary, and the functions $\epsilon_j$ are
  real-valued.
\begin{itemize}
  \item  The first terms of $\epsilon_1(h)$ (of order respectively
  $h^{0}$ and $h^1$) are given by the formul\ae\:
\begin{eqnarray}
  \label{equ:epsilon10}
  \epsilon_1^{(0)}  = c_o = \frac{1}{2\pi}\int_\delta \alpha; \\
  \label{equ:epsilon11}
  \epsilon_1^{(1)} = \frac{1}{2\pi}\int_\delta \kappa_o +
  \mu(\delta)/4,
\end{eqnarray}
where $\mu$ is the Maslov index of any regular part of $\Lambda_o$, 
and $\delta$ is any cycle associated to an $S^1$-orbit on 
$\Lambda_o\setminus\gamma$ (and recall that $\kappa_o$ is the 
sub-principal form of the system).
  \item  The first term of $\epsilon_2(h)$ is given by the formula:
\begin{equation}
  \label{equ:epsilon2}
  \epsilon_2^{(0)} = \left(\frac{\lambda 
   r_1-r_2} {|\mathcal{H}_{\Sigma}(H_2)|^{1/2}}\right)_{\restr \gamma_j}, 
\end{equation}
where $\lambda$ is defined in Lemma \ref{lemm:critic} (recall that
$r_i$ is the sub-principal symbol of $\H_i$ and
$\mathcal{H}_{\Sigma}(H_2)$ is the transversal Hessian of $H_2$). Note
that $\mathcal{H}_{\Sigma}(H_2)$ is also equal to the
$(y,\eta)$-Hessian of $H_{2,j}$ (the latter was defined along with
$\lambda$ in Lemma \ref{lemm:critic}).
\end{itemize}
\end{theo}
\begin{rema}
  Recall that there is a choice of sign in the canonical chart $\psi$
  of Theorem \ref{theo:fn-classique}. If the other sign is chosen,
  then $\epsilon_2$ becomes $-\epsilon_2$.
\end{rema}
\begin{demo}
  Consider the direct case first.  The proof goes along the same lines
  as in \cite{san-fn}, so we skip the details and indicate only the
  crucial points. First take $U$ as any \fio\ associated to $\psi$
  (note that by construction this symplectomorphism is exact in the
  sense that it preserve the action integral).  Since $H_1$ and $H_2$
  commute with $H_p$ and $H_q$, Theorem \ref{theo:fn-classique}
  implies that the principal symbols of $U^{-1}\H_1U$ and
  $U^{-1}\H_2U$ are in the classical commutant $\C_{cl}(\gamma_j)$.
  The following division lemma is easily proved as in \cite{san-fn}~:
  \begin{lemm}
    \label{lemm:division}
    Any function $K\in\C_{cl}(\gamma_j)$ that vanish on $\gamma_j$ can
    be written (in a neighbourhood of $\gamma_j$)~:
    \[ 
    K(x,\xi,y,\eta) = K(\xi,y,\eta) = a(\xi-c_o) + by\eta,
    \]
    for some smooth functions $a$ and $b$ in $\C_{cl}(\gamma_j)$.
  \end{lemm}  
  Applying this lemma to $H_1\circ\psi$ and $H_2\circ\psi$ solves the
  principal part of equation (\ref{equ:fn}).
  
  The next steps are obtained by conjugating $U$ by elliptic \pdo s,
  yielding transport equations of the form~:
  \begin{lemm}
    \label{lemm:poincare}
    Given any functions $(r_1,r_2)$ such that 
    \[ 
    \{r_1,y\eta\} = \{r_2,\xi\},
    \] 
    there exists $K_1$, $K_2$ $\in\C_{cl}(\gamma_j)$ and a function
    $f$ such that~:
    \[ 
    \{\xi,f\} = K_1 - r_1 \quad \textrm{and} \quad \{y\eta,f\} = K_2
    - r_2.
    \]
  \end{lemm}
  \begin{demo}
    Let \[ K_1(\xi,y,\eta) = \frac{1}{2\pi}\int_0^{2\pi}
    r_1(x,\xi,y,\eta)dx. \] Of course, $\{\xi,K_1\}=0$, and using the
    hypothesis of the lemma, one has
    \[ \{K_1,y\eta\} =  \frac{1}{2\pi}\int_0^{2\pi} \{r_2,\xi\}dx = \]
    \[  =  \frac{1}{2\pi}\int_0^{2\pi} -\deriv{r_2}{x}dx = 0. \]
    Now, let
    \[ f_1(x,\xi,y,\eta)=\int_0^xK_1(\xi,y,\eta) - r_1(x',\xi,y,\eta) dx'. \]
    Then $f_1$ is a smooth function on $T^*(\RM^2/\sigma)$ that
    satisfies -- using the same kind of calculation as above ~:
    \[ \{f_1,y\eta\}(x,\xi,y,\eta) = r_2(x,\xi,y,\eta) - r_2(0,\xi,y,\eta). 
    \]
    Then the wanted function $f$ is sought under the form $f=f_1+f_2$,
    which leads to the system
    \[ \{\xi,f_2\}=0, \quad \text{ and } \]
    \[ \{y\eta,f_2\}= K_2 -  r_2(0,\xi,y,\eta).\]
    It suffices to see $\xi$ as a parameter and apply a known lemma in
    the $(y,\eta)$-variables (see eg.~\cite[Theorem
    2]{guillemin-schaeffer}).
  \end{demo}
  In the reverse case, the proof of the theorem is the same provided
  we deal with functions that are invariant under $T^*\sigma$.  But if
  $K(\xi,-y,-\eta)=K(\xi,y,\eta)$, lemma \ref{lemm:division} still
  applies, yielding functions $a$ and $b$ with the same properties.
  The same is true for the transport equation.  Then each step of the
  proof can be quantised via Weyl's formula to yield well-defined \pdo
  s on $\RM^2 / \sigma$.  Thus the result still holds for the reverse
  case.
  
  The proofs for formul\ae\ (\ref{equ:epsilon10}),
  (\ref{equ:epsilon11}) and (\ref{equ:epsilon2}) are given in section
  \ref{sec:problem}, but the formula for $\epsilon_1$ is apparent in
  the proof of Theorem \ref{theo:globalQN} below, and the formula for
  $\epsilon_2$ can be directly checked using the fact that the
  subprincipal symbol is preserved under conjugation by an elliptic
  \fio\ at a critical point of the principal symbol.
\end{demo}

\subsection{Microlocal solutions}

We investigate here the solutions of the system
\begin{equation}
 \label{equ:system}
 (\H_1(h)-a_o) u = \oh, \quad (\H_2(h)-b_o) u = \oh,
\end{equation}
microlocally on a neighbourhood of the critical Lagrangian
$\Lambda_o$.  If the operators $\H_j$ depend smoothly on some
additional parameter $E\in\RM^d$ that leaves the principal symbols
intact, then all the results presented here depend smoothly and
locally uniformly on $E$.  This applies in particular to the
investigation of the joint spectrum in a window of size $O(h)$ around
$(a_o,b_o)$, where $\H_j$ is to be replaced by $\H_j-hE_j$.

Theorem \ref{theo:fn-semicla}, applied to all critical circles 
$\gamma_j$, yields a finite set of semi-classical invariants 
$(\epsilon_{1,j}(h),\epsilon_{2,j}(h))$. We show here how these 
quantities are related to the solutions of (\ref{equ:system}).

\begin{theo}[The global quantum number]~\\
  \label{theo:globalQN}
  \begin{itemize}
  \item The asymptotic series $\epsilon_1(h)=\epsilon_{1,j}(h)$,
    modulo $h\ZM$, depend neither on $j$, nor on the particular way to
    achieve the normal form of Theorem \ref{theo:fn-semicla}.
    
  \item The system (\ref{equ:system}) admits a \mi\ solution near
    any (and then all) $S^1$-orbit (including critical circles) if
    and only if the following condition holds~:
    \begin{equation}
      \epsilon_1(h) \in h\ZM + \oh.
      \label{equ:condition1}
    \end{equation}
  \end{itemize}
\end{theo}
\begin{rema}
  \label{rema:quantif-h}
  Since $\epsilon_1(h)$ is determined by $\H_1(h)$ and $\H_2(h)$, the
  fulfilment of equation (\ref{equ:condition1}) seems to impose a
  quantisation condition on $h$. While we can stick here to this
  interpretation, another possibility would be to recall that
  everything (and in particular $\epsilon_1(h)$) smoothly depends on
  the point $o$ in the curve of critical values of $F$. Then equation
  (\ref{equ:condition1}) can be interpreted as a quantisation
  condition on $o$, which leaves $h$ free to vary in a full
  neighbourhood of $0$. This viewpoint is made clear in section
  \ref{sec:problem} (cf. Corollary \ref{coro:cond1E}).
\end{rema}
\begin{demo}[of Theorem \ref{theo:globalQN}]
  Following \cite{san-these}, we introduce the sheaf $(\L,\Lambda_o)$
  of germs of \mi\ solutions on $\Lambda_o$, as a sheaf of
  $\CM_h$-modules, where $\CM_h$ is the ring of all complex functions
  of $h$, $c(h)$, such that \[ |c(h)| \leq C.h^{-N}, \] for some
  constants $C,N$, modulo those functions that are $\oh$.  Note that
  the vector operator $\F=(\H_1,\H_2)$ acts on the huge
  sheaf over $\Lambda_o$ of germs of all admissible distributions
  modulo \mi\ equivalence, and $(\L,\Lambda_o)$ can be seen as the
  kernel of $\F$.  The question is to find out how local
  germs can be glued together to form a nontrivial global section of
  $(\L,\Lambda_o)$, i.e. a solution of (\ref{equ:system}) near
  $\Lambda_o$.
  
  It was shown in \cite{san-these} that the restriction
  $(\L,\Lambda_o\setminus\gamma)$ to the non-singular points of $F$ is
  a locally constant sheaf, and the germs $\L(p)$ at any non-singular
  point $p$ form a free module of rank 1, generated by a standard WKB
  solution.  The existence of nontrivial global sections of
  $(\L,\Lambda_o\setminus\gamma)$ is then characterised by the nullity
  (mod.  $\oh$) of the associated holonomy (or ``\BS\ cocycle'')~:
  \[
  \lambda(h)\in \HUN(\Lambda_o\setminus\gamma,\RM/2\pi\ZM).
  \] 
  Since $\Lambda_o\setminus\gamma$ is a disjoint union of cylinders
  $\Lambda_{\{i,j\}}^k$, whose homology $H_1(\Lambda_{\{i,j\}}^k)$ is
  generated by the cycle represented by any oriented $S^1$-orbit, we
  get a finite set of holonomies $\lambda_{\{i,j\}}^k(h)$ characterising
  $(\L,\Lambda_o\setminus\gamma)$.
  
  Apply now Theorem \ref{theo:fn-semicla}.  The system
  (\ref{equ:system}) is then, on a neighbourhood $\Omega$ of
  $\gamma_j$, equivalent to the following standard system~:
 \begin{equation}
   Q_1u=\epsilon_{1,j}u, \quad Q_2u=h\epsilon_{2,j}u.
   \label{equ:systemQ}
 \end{equation}
 At any non-singular point $p\in\Omega\setminus\gamma_j$, the standard
 WKB $u$ solution generating $\L(p)$ is therefore of the form 
 \[ u(p) =
 \ex^{i\frac{\epsilon_{1,j}}{h}x}v(y), \quad (x,y)\sim p\in
 \RM^2/\sigma.
 \] 
 This implies that
 \[ 
 \frac{1}{2\pi}\int_{\delta}\lambda(h) \equiv
 \frac{\epsilon_{1,j}(h)}{h} + \oh \quad (\textrm{mod} \ZM),
 \]
 where $\delta$ is the cycle on $\Lambda_o\setminus\gamma$ associated
 with the orbit $S^1(p)$.  This proves
 \begin{enumerate}
 \item that $\epsilon_{1,j}(h)$ does not depend on the particular
   way to achieve the normal form;
   
 \item that $\int_{\delta}\lambda(h)$ does not depend on the choice of
   the Lagrangian cylinder containing $p$~ -- i.e. for $j$ fixed, all
   $\lambda_{\{i,j\}}^k(h)$ are equal -- which in turn proves
   
 \item that $\epsilon_1=\epsilon_{1,j}$ does not depend on the
   choice of the critical circle $\gamma_j$ (since $\Lambda_o$ is
   connected).
 \end{enumerate}
 The first part of the theorem is now proved.  Moreover, the condition
 (\ref{equ:condition1}) is necessary and sufficient for the existence
 of a (non-trivial) solution near a regular orbit.  This condition
 remains therefore necessary for the existence of a solution near a
 critical circle.  We are thus left with the proof of the sufficiency
 of this condition for critical circles, which is achieved by the next
 proposition.
\end{demo}

\begin{prop}
  \label{prop:dimension}
  Let $\gamma_j\subset\gamma$ be a critical circle. Let $d=2$ or $4$ 
  be its degree in $G(\Lambda_o)$.
  If the condition (\ref{equ:condition1}) is fulfilled, then the set 
  $\L(\gamma_j)$ of germs of \mi\ solutions on $\gamma_j$ is a free 
  $\CM_h$-module of rank $\frac{d}{2}$.
\end{prop}
\begin{demo}
  Let $n=n(h)\in\ZM$ be such that $\epsilon_1=hn + \oh$, and let
  $p\in\gamma_j$.  We know from \cite[proposition 17]{colin-p} that the
  module of \mi\ solutions of (\ref{equ:systemQ}) at $p$ is free of rank
  2, generated by
    \begin{equation}
      u_\pm \egdef \ex^{inx} \left(1_{\pm y>0}
    \frac{1}{\sqrt{|y|}}\ex^{i\epsilon_2\ln|y|}\right). 
  \label{equ:solu}
\end{equation}
  If $\gamma_j$ is direct, this immediately implies that the module
  $\L(\gamma_j)$ of \mi\ solutions of (\ref{equ:system}) on the whole
  circle $\gamma_j$ is also free and of rank 2.
  
  In the reverse case, the distribution $C_+u_++C_-u_-$ on $\T\times\RM$ 
  is invariant under the involution $\sigma$ if and only if it has the 
  parity of $n$ in the variable $y$, which reads here
  \[ C_- = \ex^{in\pi}C_+. \]
  $\L(\gamma_j)$ is in this case a free module of rank 1, and its 
  generator depends on the parity of $n$. 
\end{demo}

\subsection{The abstract \BS\ rules}
\label{sec:abstract}
We assume here that the first condition (\ref{equ:condition1}) is
fulfilled, and show that the existence of global solutions of
(\ref{equ:systemQ}) can be read on the graph $G=G(\Lambda_o)$.  As
before, let $n=n(h)\in\ZM$ be such that $\epsilon_1=n + \oh$,

Because of Theorem \ref{theo:globalQN}, for any point $p\in\Lambda_o$,
there exists a \mi\ solution on a neighbourhood of the orbit $S^1(p)$.
We shall use this fact to construct from the sheaf $(\L,\Lambda_o)$ a
new sheaf $(\bar{\L},G)$ on $G\subset W$ (recall that $W$ is the
symplectic orbifold of section \ref{sec:reduction}) that will in some
way encode whether $(\L,\Lambda_o)$ has a global section.  To each
point $p\in G$ we associate the free module $\bar{\L}(p)$ generated by
the germs of \textbf{standard basis} at $p$, which will be of rank 1,
as follows. For details of the construction when all vertices are of
direct type, see \cite{colin-p3}.

Denote by $\bgamma_j$ the vertex of $G$ corresponding to the orbit 
$\gamma_j$, and let $\bgamma=\bigcup\bgamma_j$.

\begin{itemize}
\item At a regular orbit in $\Lambda_o\setminus\gamma$, a standard
  basis is just any basis of the space of solutions near $\gamma_j$,
  so we let $\bar{\L}(p)=\L(S^1(p))$.
  
\item At a vertex $\bgamma_j$ of degree 4, we use the definition of a
  standard basis from \cite{colin-p3}, which we recall here.
 
  The edges connecting $\bgamma_j$ are oriented according to the flow
  of $y\eta$.  Moreover, near $\bgamma_j$, $W$ is a smooth oriented
  surface (the orientation is given by minus the symplectic form). It
  is shown in \cite{colin-p} that Proposition \ref{prop:dimension} in
  the direct case can be restated as follows~: let $I_1I_2$ (resp.
  $I_3I_4$) be the disjoint union of the two local edges leaving
  $\bgamma_j$ (resp.  arriving at $\bgamma_j$) with cyclic order
  $(1,3,2,4)$ -- with respect to the orientation of $W$ near
  $\bgamma_j$.  $\bar{\L}(I_1I_2)$ and $\bar{\L}(I_3I_4)$ are free
  modules of rank 2.  Then there exists a linear map
  $T_j:\bar{\L}(I_3I_4)\fleche\bar{\L}(I_1I_2)$ such that $u$ is a
  solution in a neighbourhood of $\gamma_j$ if and only if its
  restrictions satisfy $u_{\bar{\L}(I_1I_2)}= T_j
  u_{\bar{\L}(I_3I_4)}$.  In other words, if we ``feed'' the system
  with two functions on the entering edges $I_3$ and $I_4$, then these
  functions are propagated on the leaving edges $I_1$ and $I_2$ in a
  unique way (see fig.\ref{fig:1234}).
  \begin{figure}[htbp]
    \begin{center}
      \leavevmode \input{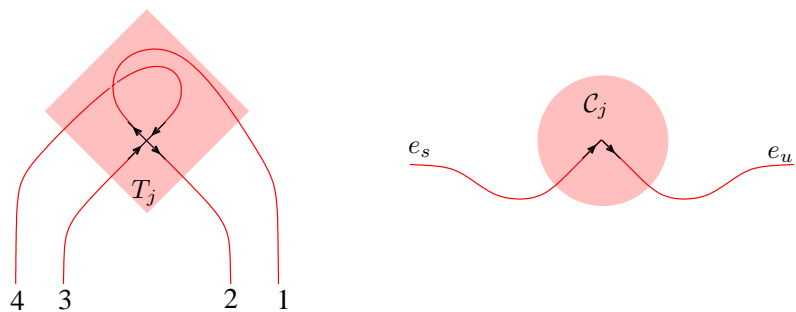}
      \caption{Propagation of solutions at vertices of degree 4 and 2}
      \label{fig:1234}
    \end{center}
  \end{figure}
  
  One can choose a basis element for each $\bar{\L}(I_i)$,
  $i=1,\dots,4$, and express $T_j$ as a $2\times 2$ matrix $\left(
    \begin{array}{cc}
      a & b \\
      c & d
    \end{array}
  \right)$ (defined modulo $\oh$).  Moreover, one can show that the
  entries are all non-vanishing.  It is then easy to check that a new
  choice for the basis elements does not change the cross-ratio
  $\rho_j=\frac{ad}{bc}$.  In our situation, $\rho_j$ can be explicitly
  calculated; one finds (\cite{colin-p})~:
  \begin{equation}
    \rho_j = \rho_j(h) = -\ex^{2\pi\epsilon_{2,j}(h)}.
    \label{equ:rho}
  \end{equation}
  The choice of a matrix $T_j$ fixes the basis elements up to their
  multiplication by a same factor.  We shall call the choice of the
  basis elements of $\bar{\L}(I_i)$, $i=1,\dots,4$, a \textbf{standard
    basis} whenever $T_j$ has the following expression~:
  \begin{equation}
    \label{equ:basis4}
    T_j = \frac{1}{\sqrt{2\pi h}}
    \Gamma(\beta)\ex^{\beta\ln h}\left(
      \begin{array}{cc}
        \ex^{-i\beta\frac{\pi}{2}} & \ex^{i\beta\frac{\pi}{2}} \\
        \ex^{i\beta\frac{\pi}{2}} & \ex^{-i\beta\frac{\pi}{2}}
      \end{array}
    \right),
  \end{equation}
  with $\beta=\frac{1}{2}+i\epsilon_{2,j}$; or with the notations of
  \cite{colin-p3}~:
  \begin{equation}
    \label{equ:basis4bis}
    T_j = T(\epsilon_{2,j}(h)) = \ex^{-i\frac{\pi}{4}}\mathcal{E}_j \left(
      \begin{array}{cc}
        1 & i\ex^{-\epsilon_{2,j}\pi} \\
        i\ex^{-\epsilon_{2,j}\pi} & 1
      \end{array}
    \right),
  \end{equation}
  with
  \begin{equation}
    \label{equ:constE}
    \mathcal{E}_j =  \mathcal{E}(\epsilon_{2,j}(h)) =
    \frac{1}{\sqrt{2\pi}}\Gamma(\frac{1}{2}+i\epsilon_{2,j}) 
    \ex^{\epsilon_{2,j}(\frac{\pi}{2}+i\ln h)}.
  \end{equation}
  \begin{rema}
    In \cite{colin-p3} the factor $\ex^{-i\frac{\pi}{4}}$ was absent in
    the definition of $T_j$ (\ref{equ:basis4bis}). Its introduction
    here will greatly simplify the treatment of Maslov indices (see
    also \cite{child-book}).
  \end{rema}
  \begin{rema}
    Equation (\ref{equ:rho}) proves that $\epsilon_{2,j}$ is a
    \semicla\ invariant (modulo $i\ZM$) of the critical circle
    $\gamma_j$~: it does not depend on the particular way the normal
    form is achieved.
  \end{rema}
  \begin{rema}
    As it is presented here, the notion of a standard basis seems to
    be attached to the graph $G$ endowed with a specific labelling at
    vertices of degree 4.  The form of the matrix $T_j$ shows that the
    different possible labellings of the four hyperbolic branches
    $I_i$, $i=1,2,3,4$ yield the same set of standard basis, provided
    $I_1$ and $I_2$ are the local unstable manifolds (for the flow of
    $y\eta$ -- which means for the flow of $H_q$), $I_3$ and $I_4$ are
    the local stable manifolds, and on the oriented manifold $W$
    (which is smooth at vertices of direct type) the branches appear
    in cyclic order $(1,3,2,4)$. Furthermore, it can be easily checked
    using the standard basis (\ref{equ:base-std}) given below and the
    fact that $\fourier^2f=\check f$ (where $\check f(y)=f(-y)$) that
    $T(\epsilon)T(-\epsilon)=\left(
      \begin{array}{cc}
        0 & 1 \\ 1 & 0
      \end{array}
    \right)$. Therefore, exchanging the local un/stable manifolds just
    amounts to changing the sign of $\epsilon_{2,j}$.
  \end{rema}
  
\item At a vertex $\bgamma_j$ of degree 2, the space of solutions has
  dimension 1, so we could just, as in the regular case, call any
  solution a standard basis.  However, in order to isolate the
  ``singular'' components of the holonomy, we prefer the following
  convention which is more in accordance with the previous case
  (degree 4).
  
  Let $I_u$ and $I_s$ be the local unstable and stable manifolds of
  $\bgamma_j$.  A choice of basis elements $(e_u,e_s)$ for
  $\bar{\L}(I_u)$ and $\bar{\L}(I_s)$ will be called a standard basis
  if
  \begin{equation}
    \label{equ:basis2}
    e_u = \mathcal{C}_j e_s,
  \end{equation}
  with
  \begin{eqnarray}
    \label{equ:constC}
    \mathcal{C}_j & = & \mathcal{C}(n,\epsilon_{2,j}(h)) =
    \ex^{-i\frac{\pi}{4}} \ex^{-in\frac{\pi}{2}} \mathcal{E}_j (1+i(-1)^n 
    \ex^{-\epsilon_{2,j}\pi}).\\
    & = & \sqrt{\frac{2}{\pi h}} \Gamma(\beta)\ex^{\beta\ln  h}
    \cos(\frac{\pi}{2}(\beta+n)).  \nonumber  
  \end{eqnarray}
  Notice that $\mathcal{C}_j$ depends on $n \mod 4$.
\end{itemize}
\begin{rema}
  If $\epsilon_{2,j}\in \RM$, one has
  \begin{equation}
    \mathcal{E}_j =
    \frac{1}{\sqrt{1+\ex^{-2\pi\epsilon_{2,j}}}}\ex^{i\arg
      \Gamma(1/2+i\epsilon_{2,j}) + i\epsilon_{2,j}\ln(h)}.
  \end{equation}
  Therefore $T_j$ is unitary and $|\mathcal{C}_j|=1$.
\end{rema}

$(\bar{\L},G)$ is a locally flat sheaf of rank-one modules, and hence 
is characterised by its holonomy: 
\[ 
\hol : \Hun(G)\fleche \CM_h.
\] 
In terms of \v Cech cohomology, if $\gamma$ is a loop in $G$, and
$\Omega_1$, $\Omega_2$, \dots, $\Omega_\ell=\Omega_1$ is an ordered
sequence of open sets covering the image of $\gamma$, each $\Omega_i$
being equipped with a standard basis $u_i$, then
\[ 
\hol(\gamma) \egdef x_{1,2}x_{2,3}\dots x_{\ell-1,\ell}, 
\] 
where $x_{i,j}$ is defined in $\CM_h$ by $u_i=x_{i,j}u_j$ on 
$\Omega_i\cap\Omega_j$.
\begin{defi}[The singular invariants]
  \label{defi:invariants}
  \begin{enumerate}
  \item The ``principal value'' $\tilde{\kappa}_o$ of the
    sub-principal form $\kappa_o$ is the cocycle on $\Lambda_o$
    defined as follows (see also \cite{colin-p3}):
    \begin{itemize}
    \item if $[A,B]\subset\Lambda_o\setminus\gamma$ is a non singular
      path, then
      \[
      \int_{[A,B]}\tilde{\kappa_o} := \int_{[A,B]} \kappa_o.
      \]
    \item if $[A,B]\subset\Lambda_o$ is a path intersecting once and
      transversally a unique critical circle $\gamma_j$, and oriented
      according to the flow of $H_q$ (ie. $A$ is on the local stable
      manifold and $B$ is on the local unstable manifold) then
      \[
      \int_{[A,B]}\tilde{\kappa_o} := \lim_{a,b\fleche m} \left(
        \int_{[A,a]}\kappa_o + \int_{[b,B]}\kappa_o +
        \epsilon_{2,j}^{(0)}\ln\left|\int_{R_{a,b}}\omega\right|
      \right),
      \]
      where $R_{a,b}$ is the parallelogram (defined in any coordinate
      system) built on the vectors $\overrightarrow{ma}$ and
      $\overrightarrow{mb}$ ($m=[A,B]\cap\gamma_j$ -- see
      fig.\ref{fig:Rab}).
    \end{itemize}
    \begin{figure}[htbp]
      \begin{center}
        \input{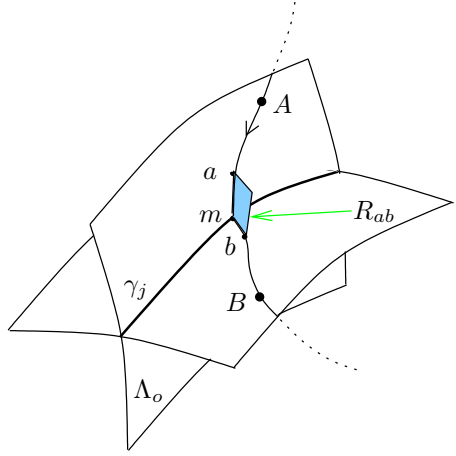}
        \caption{Regularization of $\kappa_o$}
        \label{fig:Rab}
      \end{center}
    \end{figure}
  \item The ``regularized'' Maslov index $\tilde{\mu}$ on $\Lambda_o$
    is defined as follows:
    \begin{itemize}
    \item The contribution from a regular path in
      $\Lambda_o\setminus\gamma$ is the usual Maslov index of the
      path.
    \item Let $\delta=[A,B]\subset\Lambda_o$ is a small path
      intersecting once and transversally a unique critical circle
      $\gamma_j$, such that $A$ belongs to one of the hyperbolic
      branches $(1,3,2,4)$ and $B$ to an adjacent branch (ie. $\delta$
      makes a turn of angle $\pm\frac{\pi}{2}$).  $\delta$ can be
      continuously deformed into a path $\delta_t$ drawn on a regular
      leaf of $F$. Then the usual Maslov index for this path is
      constant for $t$ small enough ($\delta_0=\delta$), and we define
      \begin{equation}
        \label{equ:maslov}
        \tilde{\mu}(\delta) := \mu(\delta_t) \pm \left(\frac{1}{2} +
          \chi_{\{d_j=2\}}n\right),
      \end{equation}
      where $\pm=$``$+$'' if $\delta$ turns in the direct sense (with
      respect to the cyclic order of the branches) and ``$-$''
      otherwise, and $\chi_{\{d_j=2\}}=1$ if $\gamma_j$ is of degree 2
      and $\chi_{\{d_j=2\}}=0$ otherwise.
      \begin{figure}[htbp]
        \begin{center}
          \input{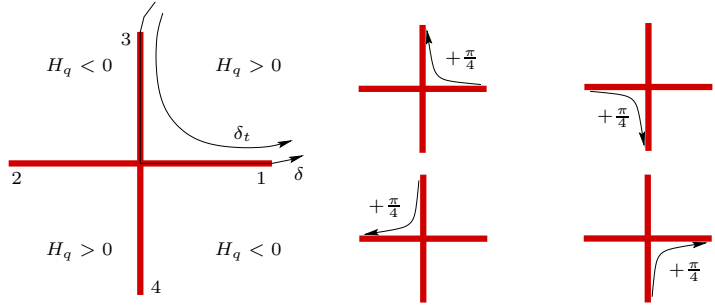}
          \caption{The local Maslov correction (for a vertex of degree
            2, replace $\frac{\pi}{4}$ by
            $\frac{\pi}{4}+n\frac{\pi}{2}$)}
          \label{fig:maslov}
        \end{center}
      \end{figure}
    \end{itemize}
  \end{enumerate}
\end{defi}
\begin{rema}
  The sign in (\ref{equ:maslov}) is negative if $\delta$ is oriented
  according to the flow of $H_q$ and $\delta_t$ belongs to a region in
  phase space where $H_q>0$, and changes whenever one of these two
  conditions changes. Of course the rule of Fig. \ref{fig:maslov} is
  simpler to use, but this correspondence will be used for the proofs
  of section \ref{sec:problem}.
\end{rema}
\begin{prop}
  \label{prop:holo}
  \begin{enumerate}
  \item The holonomy of the sheaf $(\bar{\L},G)$ has the form
    $\hol=\ex^{i[\theta(h)]/h}$, where $[\pi^*\theta(h)]\in
    \cHUN(\Lambda_o,(\CM,+))$ admits an asymptotic expansion in
    non-negative powers of $h$. ($\pi:\Lambda_o\fleche G$ is the
    projection associated to the $S^1$-reduction.)
  \item Let $\sum_{\ell\geq 0} [\tilde{\theta}_\ell] h^\ell$ be the
    asymptotic expansion of $[\tilde{\theta}(h)]:=\pi^*[\theta(h)]$.
    Then the first two terms are given by the following formul\ae\ :
    \begin{itemize}
    \item $[\tilde{\theta}_0]=[\alpha]$ (the Liouville 1-form on
      $\Lambda_o$);
    \item $[\tilde{\theta}_1] =
      [\tilde{\kappa}_o]+\tilde{\mu}\frac{\pi}{2}$.
    \end{itemize}
\end{enumerate}
\end{prop}
\begin{demo}
  We just prove here the existence of the claimed asymptotic
  expansion. The formul\ae\ for $\tilde{\theta}$ will follow from our
  refined analysis in Section \ref{sec:problem} (Corollary \ref{coro:coeff}).
  
  For this purpose, it suffices to show that one can choose local
  sections $u_\alpha$ of $(\bar{\L},G)$ for which the transition
  constants $c_{\alpha,\beta}$ have the required form. On the edges of
  $G$, this follows from the regular theory of WKB solutions (see
  eg.\cite{san-focus}). At a critical circle, we apply the normal form
  (Theorem \ref{theo:fn-semicla}), and choose the following standard
  basis ($u_\pm^{\epsilon_2}$ is defined in Eq.(\ref{equ:solu})):
  \begin{itemize}
  \item in the direct case,
    \begin{equation}
      \label{equ:base-std}
      \left\{    
        \begin{array}{rcl}
          e_1 & = & u_+^{\epsilon_2}  \\
          e_2 & = & u_-^{\epsilon_2}  \\
          e_3 & = & \fourier^{-1}(u_+^{-\epsilon_2}) \\
          e_4 & = & \fourier^{-1}(u_-^{-\epsilon_2});
      \end{array}
    \right.
  \end{equation}
  
  \item in the reverse case,
    \[
    \left\{    
      \begin{array}{rcl}
        e_u & = & e_1+(-1)^ne_2  \\
        e_s & = & i^ne_3+(-i)^ne_4.
      \end{array}
    \right.
    \]
  \end{itemize}
  We see then that the restrictions of these solutions to any edge are
  standard WKB solutions, whose phases admit an asymptotic expansion
  in powers of $h$.
\end{demo}

The dimension $b_1=\dim \Hun(G)$ is given by Euler-Poincaré formula : 
\[ 
b_1 = \#\{\textrm{edges of }G\} - N+1 
\]
(recall that $N$ is the number of vertices of $G$).  Moreover, if we 
write $N=p+q$ with $p$ the number of tetravalent vertices and $q$ the 
number of divalent vertices, then it is easy to see that
\[ 
\#\{\textrm{edges of }G\} = 2p+q,
\]
so that $b_1=p+1$.

We can now cut $b_1$ edges of $G$, each one corresponding to a cycle 
$\delta_i$ in a basis $(\delta_1,\dots,\delta_{b_1})$ of $\Hun(G)$, in 
such a way that the remaining graph is a tree $T$ ($\Hun(T)=0$).  Then 
the sheaf $(\bar{\L},T)$ has a nontrivial global section, i.e. there 
exists a standard basis $u_\alpha$ on each edge $e_\alpha$ such that 
they extend to a standard basis at each vertex.
\begin{theo}
  \label{theo:BSformel}
  $(\L,\Lambda_o)$ has a nontrivial global section if and only if the
  following linear system of $3p+q+1$ equations with the $3p+q+1$
  unknowns $(x_\alpha\in\CM_h)_{\alpha\in\{\textrm{edges of } T\}}$
  has a nontrivial solution~:
  \begin{enumerate}
  \item if the edges $(\alpha_1,\alpha_3,\alpha_2,\alpha_4)$ connect
    at a tetravalent vertex $\gamma_j$ (with the prescribed
    orientation), then
     \[ 
     (x_{\alpha_3},x_{\alpha_4}) = T_j(x_{\alpha_1},x_{\alpha_2});
     \] 
   \item if the edges $(\alpha_u,\alpha_s)$ connect at a divalent
     vertex $\gamma_j$ (with the prescribed orientation), then
     \[ 
     x_{\alpha_s} = \mathcal{C}_jx_{\alpha_u};
     \] 
  
   \item if $\alpha_0$ and $\alpha_1$ are the extremities of a cut
     cycle $\delta_i$, then
    \[ x_{\alpha_0} = \hol(\delta_i)x_{\alpha_1}. \]
    Here we assume the following orientation: $\delta_i$ can be
    represented by a closed path starting on the edge $\alpha_0$ and
    ending on $\alpha_1$.
  \end{enumerate}
\end{theo}
\begin{rema}
  When solving the system, it is immediate (if $p\neq 0$) to replace
  the equations of type 2. and 3. into those of type 1., in order to
  finally obtain a linear system of size $(2p)\times(2p)$. If $p=0$
  then equations of type 2. and 3. combine together to yield a unique
  equation in one variable.
\end{rema}
\begin{demo}
  Any global section $u$ of $(\L,\Lambda_o)$ can be characterised by
  the set of constants $x_\alpha\in\CM_h$ such that $u_{\restr
    \bar{\L}(e_\alpha)}=x_\alpha u_\alpha$.  By definition of the
  standard basis $(u_\alpha)$, conditions (1.)  and (2.)  are
  necessary and sufficient for $(x_\alpha u_\alpha)$ to extend to a
  solution near every critical circle $\gamma_j$.  It remains to check
  that the solutions at the extremities $\alpha_0,\alpha_1$ of a cut
  cycle $\delta_i$ can be consistently glued back together.  Since
  $(u_\alpha)$ is a global section of $(\bar{\L},T)$, $u_{\alpha_1}$
  is the parallel transport of $u_{\alpha_0}$ along $\delta_i$, which
  means that, as local sections of $(\bar{\L},G)$ (or
  $(\L,\Lambda_o)$), they satisfy
  $u_{\alpha_1}=\hol(\delta_i)u_{\alpha_0}$.  Therefore the solutions
  $x_{\alpha_0}u_{\alpha_0}$ and $x_{\alpha_1}u_{\alpha_1}$ can be
  glued back if and only if condition (3.)  holds.
\end{demo}

\subsection{The spectral problem}
\label{sec:problem}
The goal of this section is to investigate uniform estimates for our
system when it depends on spectral parameters. Specifically, we look
now at the system
\begin{equation}
 \label{equ:systemE}
 (\H_1(h)-E_1)u = \oh, \quad (\H_2(h)-E_2)u = \oh,
\end{equation}
where $E_1$ and $E_2$ are real numbers.  Here we shall assume that
$\H_1$ and $\H_2$ are formally self-adjoint.  If we are only
interested in studying the spectrum in a window of size $O(h)$ around
the origin, we can let $E_i=h\varepsilon_i$ and there is nothing more
to be done: Theorem \ref{theo:fn-semicla} holds uniformly for
$(\varepsilon_1,\varepsilon_2)$ varying in a compact set of $\RM^2$,
so that all the results of the previous sections apply.  However,
Theorem \ref{theo:fn-semicla} does \emph{not} apply to the system
(\ref{equ:systemE}) with uniform estimates for $E=(E_1,E_2)$ in a
compact.  Indeed, it would imply that $H_p$ has a unique value on the
local level set $(H_1,H_2)=(E_1,E_2)$ near $\gamma_j$; if $E$ is a
\emph{regular} value of $F$, this level set may fail to be connected
and it is easy to construct a situation where $H_p$ has different
values on each component.  Actually, $H_p$ by definition is a function
on the set of leaves of the Lagrangian foliation defined by $F$; and
the following diagram is in general non-commutative:
\begin{equation}
  \label{equ:factor}
  \xymatrix{ {\Omega} \ar[r]^{F} \ar @/_/ [dr]^{H_p} & {U} \ar @{.>}[d]^{??} \\
    & \RM }
\end{equation}
Instead, we need to work with the space of leaves
$\bar{\Omega}$, equipped with the momentum map $\bar{F}$.

For any $E\in U\subset (\RM^2,0)$, let $\Lambda_E=F^{-1}(E)$. If $U$
is a sufficiently small ball around some critical value, the curve
$\CC_c\subset U$ of critical values of the momentum map $F$ separates
the set of regular values in $U$ into two simply connected open sets
$U^+$ and $U^-$. Using $H_q$ defined in section~\ref{sec:fn}, we take
the following convention~: $U^\pm:=\{\pm H_q>0\}$.  Let
$D^\pm:=U^\pm\cup \CC_c$. Let $N^+$ and $N^-$ be the sets of connected
components of the open sets $F^{-1}(U^+)$ and $F^{-1}(U^-)$
respectively.  In each of $U^+$ and $U^-$, the levels sets of $F$ have
a unique topological type, namely they are unions of a finite number
$\tilde{N}^\pm$ of Liouville tori: for any $E\in U^\pm$,
$F^{-1}(E)=:\bigsqcup_{k^\pm\in N^\pm}T_{k^\pm}(E)$. Of course,
$\tilde{N}^\pm = |N^\pm|$.
\begin{prop}
 \label{prop:reeb}
 A smooth function $K$ commuting with $H_1$ and $H_2$ in $\Omega$ is
 characterised by the data of $|N^+|$ functions $f_{k^+}\in\Cinf(D^+)$
 and $|N^-|$ functions $f_{k^-}\in\Cinf(D^-)$ such that
  \begin{enumerate}
  \item For all $k^\pm\in N^\pm$, $K_{\restr k^\pm}= f_{k^\pm}\circ
    F_{\restr k^\pm}$;
  \item For all $k^+\in N^+$ and $k^-\in N^-$, the function equal to
    $f_{k^+}$ in $D^+$ and to $f_{k^-}$ in $D^-$ is smooth on $U$.
  \end{enumerate}
\end{prop}
\begin{defi}
  The space of smooth functions in $\Omega$ commuting with $H_1$ and
  $H_2$ will be denoted by $\mbox{$C^\infty$}^F(\Omega)$. The space of
  leaves together with the smooth structure described in the above
  proposition will be called the \textbf{Reeb graph} of $F$.
\end{defi}
\begin{figure}[htbp]
  \begin{center}
    \input{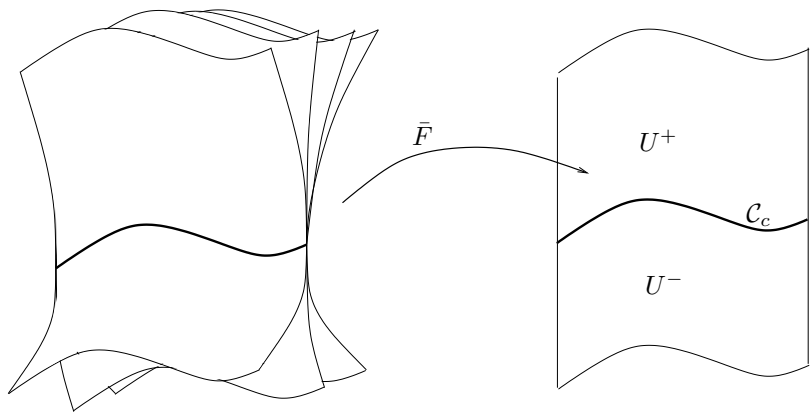}
    \caption{The Reeb graph of $F$}
    \label{fig:leaves}
  \end{center}
\end{figure}
\begin{demo}
  a) Given $K$, the condition 1.) uniquely defines the functions
  $f_{k^\pm}$. Their smoothness is given by the Darboux-Carathéodory
  theorem that states that near every non-singular point of
  $\Lambda_E$, there is a canonical chart $(x,y,\xi,\eta)$ in which
  $H_1=\xi$ and $H_2=\eta$. The fact that such non-singular point
  exist even on a singular leaf $\Lambda_{E_c}$, $E_c\in\CC_c$, shows
  that the smoothness extends to $D^\pm$. The same argument shows that
  the condition 2.) holds whenever $\k^+\cap\k^-\neq\emptyset$.  Then
  condition 2.) without this restriction holds because $\Lambda_{E_c}$
  is connected.
  
  b) Conversely, the data of all the functions $f_{k^\pm}$ defines a
  unique function $K$. The smoothness of $K$ outside of the critical
  points comes from the same argument as above; its smoothness at
  critical points comes from Theorem \ref{theo:fn-local} and the
  following lemma \ref{lemm:commutantE} (proved as in \cite{san-fn}).
\end{demo}
\begin{lemm}
  \label{lemm:commutantE}
  Each function $K$ commuting with $\xi$ and $y\eta$ near $(x,0,0,0)$
  is characterised by two functions $f_+,f_-\in\Cinf(\RM^2,0)$
  satisfying
  \begin{equation}
    \label{equ:plat}
    f_-(\xi,t)-f_+(\xi,t) = O(t^\infty),
  \end{equation}
  locally uniformly in $\xi$, such that
  \[ K(x,y,\xi,\eta) = \left\{
    \begin{array}{ll}
      f_+(\xi,y\eta) & \textrm{ if } y\geq 0 \\
      f_-(\xi,y\eta) & \textrm{ if } y < 0
    \end{array} \right. \]
\end{lemm}
\begin{rema}
  The distinction between the functions $f_+$ and $f_-$ is of course
  irrelevant when all the data is analytic. Neither has it any impact
  for \emph{semi-excited} regions, ie. for $E=O(h^\gamma)$,
  $0<\gamma<1$ (see \cite[chapter 5]{san-these}). In these cases, no
  further modification of the results of the previous sections are
  required (the diagram (\ref{equ:factor}) is always commutative), and
  the following section \ref{sec:globalQNE} becomes rather
  straightforward. We have laid down all the details to cope with the
  $\Cinf$ case, which makes the statements and proofs more technical.
  On the other hand, the statements in section
  \ref{sec:regularisation} are non-trivial even in the analytic case,
  and represent some of the most crucial results of this article.
\end{rema}


\subsubsection{The ``global'' quantum number}
\label{sec:globalQNE}
In this section, the issue is to generalise the ``global'' quantum
number $\epsilon_1(h)$ of Theorem \ref{theo:globalQN}. In the smooth,
non-analytic category, this leads to a subtle repartition property for
the \semicla\ spectrum.

For any $E\in U\setminus\CC_c$, denote by $(\L,\Lambda_E)$ the sheaf
of germs of \mi\ solutions of (\ref{equ:systemE}) on $\Lambda_E$. We
know from the regular theory (see \cite{san-these}) that
$(\L,\Lambda_E)$ is just a flat bundle characterised by its holonomy,
which was called in \cite{san-these} the \BS\ cocycle, and denoted by
$\lambda_h\in \HUN(\Lambda_E,\RM/2\pi\ZM)$.  When $E$ is restricted to
any compact subset $K\subset U^\pm$, and $\lambda_h$ is restricted to some
connected component $k^\pm$, $h\lambda_h$ has a uniform asymptotic
expansion in $\Cinf(K)\formel{h}$. This is \emph{no} longer true
on $D^\pm$. However, the following statement holds~:
\begin{theo}
  \label{theo:globalQNE}
  The function that assigns to a leaf $T_{k^\pm}(E)$ the integral
  \[
  \frac{h}{2\pi}\int_\delta \lambda_h,
  \]
  where $\delta$ is any $S^1$-orbit in $T_{k^\pm}(E)$, defines an element
  $\e_1(h)\in{\Cinf}^F(\Omega)$ that admits an asymptotic expansion of
  the form~:
  \[
  \e_1(h) = \sum_{\ell=0}^\infty h^\ell \e_1^{(\ell)} \quad \in
  {\Cinf}^F(\Omega)\formel{h},
  \]
  with $\e_1^{(0)}=H_p$, and
  $\e_1^{(1)}=-\frac{a}{2\pi}\int_\delta r_1 -
  \frac{b}{2\pi}\int_\delta r_2 + \mu(\delta)/4$, where $r_j$ is the
  sub-principal symbol of $\H_j$, and $a$, $b$ in ${\Cinf}^F(\Omega)$
  are the functions defined by $\ham{p}=a\ham{1}+b\ham{2}$, and $\mu$
  is the Maslov cocycle of $T_{k^\pm}(E)$.
\end{theo}
\begin{rema}
  \label{rema:epsilon1}
  Fix $k^\pm\in N^\pm$ and realise $(\e_1)_{\restr k^\pm}$ as a smooth
  function $E\fleche \e_1(E;h)$ in $\Cinf(D^\pm)$.  Now let
  $E=(a_o,b_o)+h(e_1,e_2)$, where $e_j$ vary in a compact (recall that
  $(a_o,b_o)$ is a critical value of $F$); then $\e_1(E;h)$ admits
  an asymptotic expansion in powers of $h$ whose coefficients are
  smooth functions of $(e_1,e_2)$. But these coefficients are made out
  of the Taylor series of $E\fleche \e_1(E;h)$, and therefore, in view
  of the definition of ${\Cinf}^F$, they do not depend on the
  component $k^\pm$. We obtain this way an element
  $\epsilon_1(h)\in\Cinf\formel{h}$ that is nothing else but the
  global quantum number of Theorem \ref{theo:globalQN}. 
  
  Notice that, in order to compute the principal term in
  $\epsilon_1(h)$, it is not interesting to use the formula
  $\epsilon_1(h)=\e_1(E;h)$, since it would involve the derivative of
  $\e_1^{(0)}$.  Instead, apply the formula of the theorem above to
  a system whose principal symbol is independent of $E$, and whose
  sub-principal symbol is $(r_1-e_1,r_2-e_2)$. We obtain this way the
  claim (\ref{equ:epsilon10})-(\ref{equ:epsilon11}) in Theorem
  \ref{theo:fn-semicla}.
\end{rema}
\begin{demo}[of the theorem]
  The fact that $\e_1(h)\in{\Cinf}^F(\Omega)$ is obvious form the
  construction. To prove the existence of the claimed asymptotic
  expansion, it suffices to \mi ize near a critical circle $\gamma_j$.
  
  Using Lemma \ref{lemm:commutantE} in the coordinates of Theorem
  \ref{theo:fn-classique}, one sees that the functions $f_+$ and $f_-$
  of Lemma \ref{lemm:commutantE} are the same if the degree $d=2$. In
  this case, Theorem \ref{theo:fn-semicla} generalises to
  \[
  U^{-1}(\H_1-E_1,\H_2-E_2)U = \M.(Q_1-\e_1,Q_2-\e_2) + \oh,
  \]
  where $\M$, $\e_1$, and $\e_2$ depend smoothly on $E$, which gives
  the result.
  
  The case $d=4$ is more intricate, and follows from Proposition
  \ref{prop:ansatzE} below.
\end{demo}
We shall need the following slightly weaker version of Theorem
\ref{theo:fn-semicla}~:
\begin{prop}
  \label{prop:fn-weaker}
  There exists an elliptic \fio\ $U(h)$ associated to the canonical
  transformation $\psi$ of Theorem \ref{theo:fn-classique} such that,
  \mi ly near $\gamma_j$~:
  \[U^{-1}\H_j U = \K_j,\]
  where $\K_j\in \C_h(\gamma_j)$.
\end{prop}
\begin{demo}
  The same proof scheme as that of Theorem \ref{theo:fn-semicla}
  applies, using Lemma \ref{lemm:commutantE}  instead of Lemma
  \ref{lemm:division}.
\end{demo}
\begin{prop}
  \label{prop:ansatzE}
  Let $\gamma_j$ be a critical circle of degree 4.  For each $E$ close
  to zero, the set of microlocal solutions of (\ref{equ:systemE}) on a
  small neighbourhood of any point of $\gamma_j$ is a free
  $\CM_h$-module of rank 2. In the coordinates of Theorem
  \ref{theo:fn-classique}, it has a basis of the form
  \[u^\pm_E = \ex^{i\e_1^\pm x/h}\left(1_{\pm
      y>0}\frac{1}{\sqrt{|y|}}\ex^{i\e_2^\pm\ln |y|/h} \right),\] where
  $\e^\pm_j=\e^\pm_j(E;h)$ admits an asymptotic expansion in
  non-negative powers of $h$ whose coefficients are smooth functions
  of $E$.  For this basis, the system (\ref{equ:systemE}) is solved
  locally uniformly with respect to $(E_1,E_2)$ near $(0,0)$.
  Moreover, the functions $(E_1,E_2)\fleche\e_j^+-\e_j^-$ are flat on
  the set $\CC_c$ of critical points of $F$.
\end{prop}
The proof of this proposition relies on the following lemma~:
\begin{lemm}
  \label{lemm:ansatzE}
  Let $p\in\Cinf_0(T^*S^1\times T^*\RM)$ be a Hamiltonian satisfying
  \[
      p(x,\xi,y,\eta) = 0  \textrm{ for } y\geq 0
  \]
  (i.e. $p=1_{y\leq 0}p$). Then
  \[ 
  Op^W_h(p) = Op^W_h(p)\circ 1_{y\leq 0} = 1_{y\leq
    0}.Op^W_h(p)  \quad (\textup{mod } \oh).
  \] 
\end{lemm}
Of course the symmetric result (with respect to $y=0$) holds.

\begin{demo}
  Recall that Weyl quantisation of $p$ is defined by~:
  \[ v(x,y)=Op_h^W(p)u(x,y) = \]
  \[ \frac{1}{(2\pi h)^2} \int \ex^{\frac{i}{h}((x-x')\xi+(y-y')\eta)}
  p(\frac{x+x'}{2},\xi,\frac{y+y'}{2},\eta) u(x',y')dx'dy'd\xi d\eta.
  \] 
  We prove the first estimate by showing that $\|v\|=\oh$ whenever
  $u=1_{y\geq 0}u$. For this, we consider the two regions $|y|\geq
  h^{\delta_1}$ and $|y|\leq h^{\delta_2}$, with $0<\delta_2\leq
  \delta_1<1$. If $|y|\leq h^{\delta_2}$, then only the domain
  $|\frac{y+y'}{2}|\leq h^{\delta_2}$ contributes to the
  integral; and the result follows from the flatness of $p$ with
  respect to its third variable~: for all $N\in\NM,
  |v|=O(h^{N\delta_2})$. Let us now look at the region $|y|\geq
  h^{\delta_1}$. Since $v(\cdot,\cdot,y\geq 0,\cdot)=0$, one can
  assume that $y\leq -h^{\delta_1}$, which implies $|y-y'|\geq
  h^{\delta_1}$. Now the usual trick applies~: a repeated integration
  by parts with respect to the operator
  $\frac{h}{i(y-y')}\deriv{}{\eta}$ (or standard estimates for the
  Fourier transform) gives $|v|=O(h^{N(1-\delta_1)})$
  for any integer $N$.
  
  The same method can be applied to show that $\| 1_{y\geq 0}v\| =
  \oh$ whenever $u= 1_{y\leq 0} u$, thus proving the second estimate.
\end{demo}
\begin{demo}[of Proposition \ref{prop:ansatzE}]
  The fact that the set of solutions is a free module of rank 2 is
  due, for $E\not\in\CC_c$, to the regular theory (the local
  Lagrangian manifold has two connected components, on each of which
  the set of solutions is a free module of rank 1), and, for
  $E\in\CC_c$, to Proposition \ref{prop:dimension}.
  
  We prove the rest of the proposition for $u_E^+$; the same argument
  applies to $u_E^-$.  First apply proposition \ref{prop:fn-weaker} to
  assume in what follows that $\H_j\in \C_h(\gamma_j)$. Since $y\eta$
  is a quadratic function, every element of $\K\in\C_h(\gamma_j)$ can
  be written $\K=Op^W_h(K_h)$, with $K_h\sim\sum_{\ell\geq 0}h^\ell
  K^{(\ell)}$ and $K^{(\ell)}\in\C_{cl}(\gamma_j)$.  Because of Lemma
  \ref{lemm:commutantE}, each $K^{(\ell)}$ is defined by two functions
  $f^{(\ell)}_\pm$.  Let
  \[ 
  F_h(x,\xi,y,\eta)\sim \sum_{\ell\geq 0}h^\ell 
  f^{(\ell)}_+(\xi,y\eta),
  \] 
  and $R_h=K_h-F_h$.  Let us prove now that there is a unique symbol 
  $g_h(\xi,t)\sim\sum_{\ell\geq 0}h^\ell g^{(\ell)}(\xi,t)$ such that
  \[ 
  Op^W_h(F_h) = g_h(Q_1,Q_2) + \oh,
  \]
  ($Q_j$ is defined in Eq.(\ref{equ:opQ})).  Indeed, $g^{(0)}$ is 
  necessarily equal to $f^{(0)}_+$; therefore,
  \[
  Op^W_h(F_h) = g^{(0)}(Q_1,Q_2) + h\hat{S},
  \]
  where $\hat{S}\in\C_h(\gamma_j)$ and is of order $0$.  Then 
  $\hat{S}$ is similarly decomposed -- and the claim is proved by 
  induction -- provided we show that its Weyl symbol is, as for $F_h$, 
  a function of $(\xi,y\eta)$.  This is achieved by applying Lemma 
  \ref{lemm:commutantE} and remarking that $Op^W_h(F_h)$, as well as 
  $g^{(0)}(Q_1,Q_2)$, commute with the involution $y\fleche -y$ (and 
  thus their Weyl symbols are invariant under 
  $(y,\eta)\fleche(-y,-\eta)$).  Summing up, we have proved so far 
  that any operator $\K\in\C_h(\gamma_j)$ can be written~: 
  \[ 
  \K = 
  g_h(Q_1,Q_2) + Op^W_h(R_h),
  \]
  where all the coefficients in the expansion of $R_h$ are smooth 
  functions verifying the hypothesis of Lemma \ref{lemm:ansatzE}.

  Applying this to $\H_j$, we obtain the existence of two symbols
  $g_{1,h}$ and $g_{2,h}$ such that
  \[
  \H_j u_E^+ = g_{j,h}(\e^+_1,\e^+_2)u_E^+ + \oh.
  \]
  The independence of $H_1$ and $H_2$ implies that the principal
  term~:
  \[
  (\xi,t)\fleche(g_1^{(0)}(\xi,t),g_2^{(0)}(\xi,t))
  \]
  is local diffeomorphism; therefore the symbol $(g_{1,h},g_{2,h})$ is
  invertible, and the proposition is proved with
  \[
  (\e^+_1,\e^+_2) \sim (g_{1,h},g_{2,h})^{-1}(E_1,E_2).
  \]
\end{demo}

\begin{figure}[htbp]
  \begin{center}
    \input{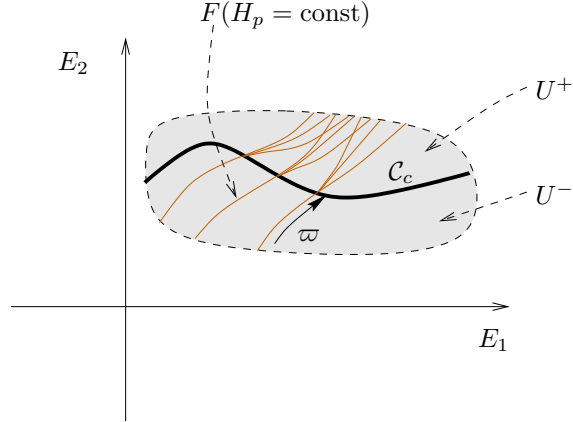}
    \caption{Level sets of $H_p$}
    \label{fig:feuillesbadE}
  \end{center}
\end{figure}

\begin{coro}
  \label{coro:cond1E}
  Fix $k^\pm\in N^\pm$ and realise $(\e_1)_{\restr k^\pm}$ as a smooth
  function $E\fleche \e_1(E;h)$ in $\Cinf(D^\pm)$. Then the condition
  \begin{equation}
    \e_1(E;h) \in h\ZM + \oh
    \label{equ:cond1E}
  \end{equation}
  is necessary and sufficient for the existence of a uniform solution
  of (\ref{equ:systemE}) \mi ised in a neighbourhood (in
  $\overline{k}^\pm$) of any $S^1$-orbit in $\overline{k}^\pm$.
\end{coro}
\begin{defi}
  For any natural integers $\ell^+$ and $\ell^-$, we call an
  \textbf{$(\ell^+,\ell^-)$-curve} the union of $\ell^+$ smooth curves in $D^+$
  and $\ell^-$ smooth curves in $D^-$ that are transversal to $\CC_c$
  and infinitely tangent to each other on $\CC_c$.
\end{defi}
These curves are just the image by $F$ of a level set of a smooth
function $K$ on the Reeb graph of $F$, if $\ham{K}\neq 0$ on the
critical leaves $\Lambda_{E_c}$ (this is a consequence of Proposition
\ref{prop:reeb}). This holds for instance for $K=H_p$
(Fig.~\ref{fig:feuillesbadE}).
\begin{coro}
  There exists a fixed neighbourhood $U$ in $\RM^2=\{(E_1,E_2)\}$ of
  any critical value of $F$ in which the joint spectrum of $\H_1$ and
  $\H_2$ is distributed (modulo $\oh$) on the union of
  $(|N^+|,|N^-|)$-curves $L_n(h)$ ($n\in\ZM$) defined as the image by
  $F$ of the level sets $\e_1(h)=hn$. The principal part of these
  curves is thus given by the level sets of $H_p$.
\end{coro}

\begin{demo}
  Fix $k^\pm\in N^\pm$. Then Proposition \ref{prop:reeb} says that the
  restriction $(\e_1(h))_{\restr \k^\pm}$ is equal to $f_h\circ
  F_{\restr \k^\pm}$ for some smooth function
  $f_h\sim\sum_{\ell=0}^\infty f_\ell$ admitting an asymptotic
  expansion in $\Cinf(\RM^2,0)\formel{h}$. Since $f_0\circ F_{\restr
    \k^\pm}=(H_p)_{\restr \k^\pm}$, the hypothesis
  (\ref{equ:transv-Hp}) of Theorem \ref{theo:periodicham} implies that
  the foliation $f_h=\textrm{const}$ is transversal to $\CC_c$, and we
  can define the projection $\varpi_h:\RM^2\fleche\CC_c$ such that
  $\varpi_h(f=\textrm{const})$ is a point. The pre-image of
  $\{f_h=\textrm{const}\}\cap D^\pm$ by $F_{\restr \k^\pm}$ is a leaf
  of the foliation $\{\e_1(h)=\textrm{const}\}$ in $\k^\pm$. The value
  of $H_p$ on $\CC_c$ can be taken as a coordinate on $\CC_c$, and via
  this identification, it is natural to view $\e_1(h)$ as a function
  with values in $\CC_c$. By Corollary \ref{coro:cond1E}, any \mi\ 
  eigenfunction of $(\H_1,\H_2)$ \mi ised in $\k^\pm$ defines a joint
  eigenvalue that belongs to $\varpi_h^{-1}(h\ZM + \oh)$.
%
  \begin{equation}
    \label{equ:projection}
    \xymatrix{
      {\k^\pm} \ar[r]^{F}  \ar @/_/ [drr]^{H_p} & {U}
      \ar[r]^{\varpi_0} & {\CC_c} \ar[d]^\simeq \\
      & & \RM }
  \end{equation}
\end{demo}

\subsubsection{Regularization of $\lambda_h(E)$}
\label{sec:regularisation}

This section contains some of the most central results of this article
(Theorems \ref{theo:local-holonomie} and \ref{theo:holonomie-graphe}).
They give a new interpretation of the holonomy of section
\ref{sec:abstract} and provide for the proofs of the various
formul\ae\ claimed before.

To each $k^\pm\in N^\pm$ and $E_c\in\CC_c$ we associate the subset 
$T_{k^\pm}(E_c)=\k^\pm\cap\Lambda_{E_c}$ which is the ``limit'' of the 
torus $T_{k^\pm}(E)$ as $E\fleche E_c$.  Suppose we are given a 
continuous family of piecewise differentiable loops $(\delta_E)_{E\in 
D^\pm}$ on $T_{k^\pm}(E)$ that are everywhere transversal to the 
$S^1$-orbits (such a family can be constructed using for instance the 
normal form of Theorem \ref{theo:fn-classique}), and assume that they 
are oriented by the flow of $H_q$.  For non-singular values of $E$, 
$\delta_E$ together with an $S^1$-orbit form a basis of 
$\Hun(T_{k^\pm}(E))$, and $E\fleche \int_{\delta_E}\lambda_h$ is a 
smooth function.  To complete the result of Theorem 
\ref{theo:globalQNE}, it is natural to investigate here the behaviour 
of that function as $E$ approaches a critical value.

To each $E_c\in\CC_c$ corresponds a real number $x$ via the
diffeomorphism (\ref{equ:projection}); $T_{k^\pm}(E_c)$ is an
$S^1$-invariant subset of $H_p^{-1}(x)$, and hence can be reduced to a
cycle $G_{k^\pm}(E_c)$ of the graph $G(E_c):=G(\Lambda_{E_c})$ in the
reduced manifold $W(E_c)$.  The second goal of this section is to show
how the asymptotic behaviour of the function $E\fleche
\int_{\delta_E}\lambda_h$ is related to the \emph{holonomy} $\hol$ of
section \ref{sec:abstract}.  Actually we shall prove that
$\int_{\delta_E}\lambda_h$ diverges as $E$ approaches a critical
value; but there is a universal way of \emph{regularizing} this
divergence. The regularized value is precisely $\hol(\bar{\delta}_E)$
-- modulo some Maslov corrections in presence of vertices of degree 2
--, where $\bar{\delta}_E$ is the projection of $\delta_E$ onto the
reduced manifold and is actually equal to $G_{k^\pm}(E_c)$.

\begin{figure}[htbp]
  \begin{center}
    \input{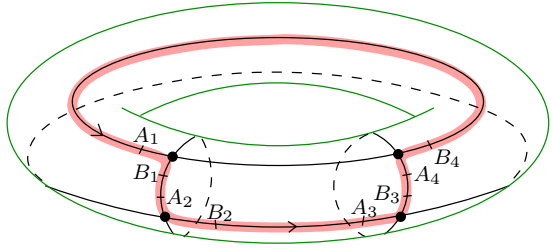}
    \caption{A graph on a torus. The cycle in gray cannot be obtained
      as a combination of boundary faces of the graph.}
    \label{fig:graphe-tore}
  \end{center}
\end{figure}

Unfortunately, the set of cycles $G_{k^\pm}(E_c)$ doesn't necessarily
generate the group $\Hun(G(E_c),\ZM)$ -- see eg.  Fig.
\ref{fig:graphe-tore} -- (but it does indeed if $G(E_c)$ is planar).
So there's a little bit more to it than just taking the limit of
regular cycles.

To compute $\hol$ of all possible cycles of $G(E_c)$, we have to
replace the ``natural'' object $\delta_E$ with a local path near a
critical circle.  In order to give some sense to the expression
$\int_{\delta_E}\lambda_h$, where $\delta_E$ is not closed, we could
abstractly choose a smooth family of closed 1-forms on $\Lambda_E$,
$E\in U^\pm$, whose cohomology class is $[\lambda_h(E)]$ -- which is
always possible as $\HUN(\Lambda_E)\fleche U^\pm$ is a trivial bundle
if $U$ is small enough. However, this does not allow us to have a
local control of the divergence of the holonomy.  Instead, we follow
\cite{san-focus} and interpret $\lambda_h$ as the phase of the
multiplicative \v Cech holonomy of the sheaf $(\L,\Lambda_E)$, as
follows:
\begin{defi}
  \label{defi:local-holonomy}
  Near each $\Gamma_j$, we let $A_j(E)$ and $B_j(E)$ in
  $\Cinf(D^\pm,\bar{k}^\pm\setminus\Gamma_j)$ be families of points
  such that for a critical value $E_c$, $A(E_c)$ and $B(E_c)$ lie
  respectively in the local stable or unstable manifold.  We endow a
  small neighbourhood of $A_j$ (or $B_j$) with a standard WKB \mi\ 
  solution $u_{A_j}$ (resp.  $u_{B_j}$) whose phase admits an
  asymptotic expansion in $h^{-1}\Cinf(D^\pm)\formel{h}$.
  
  Then the integrals $\int_{A_j}^{B_{j'}}\lambda_h$ are defined as the 
  phase of the \v Cech holonomy of $(\L,\Lambda_E)$ for paths joining 
  $A_j$ and $B_{j'}$ with the sections $u_{A_j}$ and $u_{B_{j'}}$ 
  fixed.
  
  In other words, if the path $\delta$ between $A_j$ and $B_{j'}$ is
  covered by open sets $\Omega_0,\dots,\Omega_\ell$, each of which
  being endowed with a \mi\ solution $u_\alpha$ with $u_0=u_{A_j}$ and
  $u_\ell=u_{B_{j'}}$, then
  \[ 
  \int_\delta \lambda_h = -i\ln\left(c_{0,1}c_{1,2}\cdots 
  c_{\ell-1,\ell}\right),
  \] 
  where $c_{i,j}$ is the transition constant $u_i=c_{i,j}u_j$ on
  $\Omega_i\cap\Omega_j$.

\end{defi}
Note that if another admissible choice for the local sections
$u_{A/B}$ is made, then the holonomy is modified by an additive term
admitting an asymptotic expansion in $h^{-1}\Cinf(D^\pm)\formel{h}$.
Therefore the singular behaviour of the holonomy at a critical value
is fully preserved. Note also that this additive term is necessary a
\v Cech \emph{coboundary}, and hence has no influence on the value of
the holonomy along a \emph{closed} loop.

\begin{defi}
  \label{defi:deltaE}
  In what follow, $(\delta_E)_{E\in D^\pm}$ designates a continuous 
  family of paths in $T_{k^\pm}(E)$ such that:
  \begin{itemize}
  \item for each $E\in U^\pm$, $\delta_E$ is smooth;
   
  \item either for all $E_c\in\CC_c$, $\delta_{E_c}$ does not meet the
    critical set $\Gamma$ (then $(\delta_E)$ is called
    \textbf{regular}) or for each $E_c$, $\delta_{E_c}$ meets uniquely
    a unique critical circle $\gamma_j$ (in which case $(\delta_E)$ is
    called \textbf{local});
    
  \item the end points $A_j(E)$ and $B_{j'}(E)$ are one of those
    defined in Definition \ref{defi:local-holonomy}, and we will write
    $\delta_E=[A_j(E),B_{j'}(E)]$;
  
  \item $\delta_E$ is always transversal to the $S^1$-orbits.
    
  \end{itemize}
\end{defi}
Here again, the normal form of Theorem \ref{theo:fn-classique} proves 
the existence of such a family of paths near any critical circle 
$\gamma_j$. 

The goal of this section is finally to investigate the behaviour of
the function $E\fleche\int_{\delta_E}\lambda_h$ as $E$ tends to a
critical value, and to relate it to the holonomy of the sheaf
$(\bar{\L},G(E_c))$.  The previous case where $\delta_E$ was
a loop can always be recovered by composing paths of the type of
Definition \ref{defi:deltaE}. Moreover, the regular theory implies
that the so-called ``local'' paths can indeed be restricted to paths
that are local in small neighbourhoods of the critical circles, since
the following proposition holds:
\begin{prop}
  \label{prop:regular-holonomy}
  If $\delta_E=[B_j(E),A_{j'}(E)]$ is a \emph{regular} family of paths
  (in the sense of Definition \ref{defi:deltaE}), then $E\fleche
  \int_{\delta_E}\lambda_h$ is smooth in $D^\pm$ and admits an
  asymptotic expansion in $h^{-1}\Cinf(D^\pm)\formel{h}$. The first
  terms of this expansion are the following:
  \begin{equation}
    \label{equ:transport}
    \int_{\delta_E}\lambda_h = \Phi_{B_j}(B_j)-\Phi_{A_{j'}}(A_{j'}) +
    \frac{1}{h}\int_{\delta_E} \alpha + \int_{\delta_E} \kappa +
    \mu(\delta_E)\frac{\pi}{2} + O(h),
  \end{equation}
  where $\Phi_{B_j}$ (resp. $\Phi_{A_{j'}}$) is the phase of the
  principal symbol (viewed as a section of the Keller-Maslov bundle
  over the Lagrangian manifold $\Lambda_E\setminus\Gamma_j$ -- see eg.
  \cite{duistermaat-oscillatory} or \cite{weinstein-bates}) of the
  fixed solution $u_{B_j}$ (resp.  $u_{A_{j'}}$).
\end{prop}


To study the neighbourhood of a critical circle, we shall use Theorem 
\ref{theo:fn-local} to generalise the \semicla\ invariant 
$\epsilon_{2,j}(h)$ of equation (\ref{equ:rho}) in a better way than 
Proposition \ref{prop:ansatzE} would do.  That theorem still holds if 
$H_1$ and $H_2$ are replaced by $H_1-E_1$ and $H_2-E_2$, for a 
parameter $E=(E_1,E_2)$ varying near $(a_o,b_o)$.  If we fix a 
critical circle $\gamma_j$ and $m\in\gamma_j$, the theorem yields a 
canonical change of coordinates $(x,y,\xi,\eta)$ near $m$, depending 
smoothly on $E$, and a function $\Phi_E\in\Cinf(\RM^2,0)$ depending 
smoothly on $E$, such that
\begin{equation}
\label{equ:fn-localeE}
H_1-E_1 = \xi, \quad H_2-E_2=\Phi_E(\xi,y\eta).
\end{equation}
(We have still $\partial_2\Phi_E(0,0)>0$.)  This leads to yet another
\semicla\ normal form~:
\begin{prop}
  \label{prop:fn-locale-semicla}
  Let $\gamma_j$ be a critical circle and $m\in\gamma_j$.  There
  exists a \mi ly unitary \fio\ $U(h)$ associated to the canonical
  coordinates $(x,y,\xi,\eta)$, elliptic \pdo s $M_1(h)$, $M_2(h)$
  commuting (modulo $\oh$) with $Q_1$ and $Q_2$, and a real-valued
  function of $h$~ (independent on $m$, $U$, $M_i$):
  $\e_2=\e_{2,j}(E;h)\sim\sum_{\ell=0}^\infty \e_{2,j}^{(\ell)}(E)
  h^\ell$, such that, microlocally near $m$~:
  \begin{eqnarray*}
    U^{-1}(\H_1-E_1)U & =  & Q_1 + \oh \\
    U^{-1}(\H_2-E_2)U & = &  M_1Q_1 +  M_2.(Q_2-\e_2) + \oh.
  \end{eqnarray*}
  $M_i$, $U$ and $\e_2$ depend smoothly on $E$.
 \begin{itemize}
 \item $\e_{2,j}^{(0)}(E)$ is equal to the value of $y\eta$ on
   $\Lambda_E$. In particular, $\e_{2,j}^{(0)}>0$ if $E\in U^+$ and
   $\e_{2,j}^{(0)}<0$ if $E\in U^-$;
 
 \item If $E=E_c\in\CC_c$,
   \[ 
   \e_{2,j}^{(1)}(E) = \left(\frac{\lambda r_1-r_2}
     {|\mathcal{H}_{\Sigma}(H_2)|^{1/2}}\right)_{\restr
     \Gamma_j\cap\Lambda_E},
   \]
   where $r_i$ is the sub-principal symbol of $\H_i$,
   $\lambda=\lambda(E_c)$ is the unique real number such that
   $H_2-\lambda H_1$ is critical on $\Gamma_j\cap\Lambda_E$ (see Lemma
   \ref{lemm:critic}), and $|\mathcal{H}_{\Sigma}(H_2)|$ is the
   absolute value of the determinant of the transversal Hessian of
   $H_2$.  Note also that this denominator is equal to
   $\partial_2\Phi_E(0,0)$, and $\lambda=\partial_1\Phi_E(0,0)$.
 \end{itemize}
\end{prop}
\begin{rema}
  (See Remark \ref{rema:epsilon1}).  If $E$ is restricted to a domain
  of the form $E=(a_o,b_o)+h(e_1,e_2)$, where $(a_o,b_o)\in\CC_c$,
  then Theorem \ref{theo:fn-semicla} applies with $e_1$ and $e_2$ as
  parameters, and yields an invariant $\epsilon_{2,j}$, which can be
  recovered from $\e_{2,j}$ by the following formula~:
   \[
   \epsilon_{2,j}(e_1,e_2) = \frac{1}{h}\e_{2,j}((a_o,b_o)+h(e_1,e_2))
   + \oh,
   \]
   or merely by viewing $-(e_1,e_2)$ as a correction of the
   subprincipal symbols and applying the formul\ae\ of the
   Proposition. This proves the claim (\ref{equ:epsilon2}) of Theorem
   \ref{theo:fn-semicla}.
\end{rema}
Using this proposition, we let $\beta=\frac{1}{2}+i\e_2/h$ and
$\zeta^\pm_j$ be the $h$-dependent functions in $\Cinf(U)$ defined by
\begin{eqnarray}
  \zeta^+_j & :=  & \frac{1}{\sqrt{2\pi h}}
  \Gamma(\beta)\ex^{\beta\ln h}  \ex^{-i\beta\frac{\pi}{2}} =
  \ex^{-i\frac{\pi}{4}} \mathcal{E}_j(\e_{2,j}/h); \\ 
  \zeta^-_j & :=  & \frac{1}{\sqrt{2\pi h}}
    \Gamma(\beta)\ex^{\beta\ln h}  \ex^{i\beta\frac{\pi}{2}} =
  \ex^{i\frac{\pi}{4}} \ex^{-\e_{2,j}\pi/h}\mathcal{E}_j(\e_{2,j}/h).
  \label{equ:coeff-transition}
\end{eqnarray}
($\mathcal{E}_j$ was defined in (\ref{equ:constE}).) Next lemma, which
directly follows from Stirling's formula, will be crucial for our
analysis.
\begin{lemm}
  \label{lemm:stirling}
  For any $E\in U^\pm$,
  \[
  \frac{1}{i}\ln\zeta^\pm_j = \frac{1}{h}(\e_{2,j}^{(0)}\ln
  |\e_{2,j}^{(0)}| - \e_{2,j}^{(0)}) + \e_{2,j}^{(1)}\ln
  |\e_{2,j}^{(0)}| \mp \frac{\pi}{4} + O_E(h).
  \]
\end{lemm}
\begin{theo}
  \label{theo:local-holonomie}
  Fix a component $k^\pm$, and let $\delta_E=[A_j(E),B_j(E)]$ be a
  local path near the critical component $\Gamma_j$ (see Definition
  \ref{defi:deltaE}).  Assume moreover that $\delta_E$ is oriented
  according to the flow of $H_q$ (otherwise just take the opposite of
  the holonomy !).  Then there exists an $h$-dependent
  $\RM/2\pi\ZM$-valued function $g_{\delta}(h) : E\fleche
  g_{\delta_E}(h) \in\Cinf(D^\pm)$ admitting a uniform asymptotic
  expansion of the form
    \[
    g_{\delta}(E;h) \sim \sum_{\ell =-1}^\infty
    g_{\delta}^{(\ell)}(E)h^\ell, \quad
    g_{\delta}^{(\ell)}\in\Cinf(D^\pm),
    \]
    such that
    \begin{equation}
      \label{equ:regularisation}
      \forall E\in U^\pm, \quad g_{\delta}(E;h) = \int_{\delta_E}\lambda_h -
      i\ln(\zeta_j^\pm(E)) \quad(\textup{mod } 2\pi\ZM).
    \end{equation}
    The principal terms of $g_{\delta}(h)$ are given by the following
    formul\ae, for $E\in U^\pm$:
    \begin{equation}
      \label{equ:holonomieE0}
      g_{\delta}^{(-1)}(E)  =  \int_{\delta_E}\alpha + \left(
        \e_{2,j}^{(0)}\ln |\e_{2,j}^{(0)}| -\e_{2,j}^{(0)}\right) +
      \Phi_{A_j}^{(-1)}(A_j) - \Phi_{B_j}^{(-1)}(B_j);
    \end{equation}
    \begin{equation}
      \label{equ:holonomieE1}
      g_{\delta}^{(0)}(E)  =  \int_{\delta_E}\kappa_E + 
      \mu(\delta_E)\frac{\pi}{2} + \left( \mp\frac{\pi}{4} + 
      \e_{2,j}^{(1)}\ln |\e_{2,j}^{(0)}|\right) +
        \Phi_{A_j}^{(0)}(A_j) - \Phi_{B_j}^{(0)}(B_j),
    \end{equation}
    where $\Phi_{A_j}=\frac{1}{h}\Phi_{A_j}^{(-1)}+\Phi_{A_j}^{(0)}$
    is the phase of the principal symbol of the fixed solution
    $u_{A_j}$ (and similarly for $u_{B_j}$).
\end{theo}
\begin{theo}
  \label{theo:holonomie-graphe}
  Let $E_c=(a_o,b_o)\in\CC_c$.  Let $\tilde{\delta}_{E_c}$ be a loop
  in $\Lambda_{E_c}$ oriented according to the flow of $H_q$ and of
  the form
  \[
  \tilde{\delta}_{E_c} = \delta^{\textup{loc}}_1\cdot
  \delta^{\textup{reg}}_1 \cdot\delta^{\textup{loc}}_2\cdot
  \delta^{\textup{reg}}_2\cdots \delta^{\textup{loc}}_q\cdot
  \delta^{\textup{reg}}_q,
  \]
  where $\delta^{\textup{loc}}_\ell$ and $\delta^{\textup{reg}}_\ell$
  are respectively ``local'' and ``regular'' paths in the sense of
  Definition \ref{defi:deltaE}.  (The components $k^\pm$ used for
  these paths may vary. -- see eg. Fig \ref{fig:graphe-tore}). Let
  \[
  g(E_c;h) \sim \sum_{\ell =-1}^\infty g^{(\ell)}(E_c)h^\ell, \quad
  g^{(\ell)}\in\Cinf(\CC_c)
  \]
  be defined as the sum
  \[
  g(E_c;h) := \left(g_{\delta^{\textup{loc}}_1} + g_{\delta^{\textup{reg}}_1} +
  \cdots + g_{\delta^{\textup{loc}}_q} +
  g_{\delta^{\textup{reg}}_q}\right)_{\restr E=E_c},
  \]
  where $g_{\delta^{\textup{loc}}_k}$ is given by Theorem
  \ref{theo:local-holonomie} and $g_{\delta^{\textup{reg}}_k}\egdef
  \int_{\delta^{\textup{reg}}_k}\lambda_h$ (see Proposition
  \ref{prop:regular-holonomy}).  Then, under the hypothesis of section
  \ref{sec:abstract},
  \begin{equation}
    \label{equ:holE}
    \hol(\bar{\delta}_{E_c}) := \ex^{ig(E_c;h)}
    \ex^{i\frac{\pi}{2}n(N_2^--N_2^+)} + \oh,
  \end{equation}
  where $\bar{\delta}_{E_c}$ is the projection of
  $\tilde{\delta}_{E_c}$ onto the graph $G(E_c)$ in the reduced
  orbifold $W(E_c)$, $n$ is the ``global quantum number'' of Theorem
  \ref{theo:globalQN}, and $N_2^\pm$ is the number of local paths
  through a vertex of degree 2 and defined by a component $k\in
  N^\pm$.
\end{theo}
\begin{demo}[of Theorem \ref{theo:local-holonomie}]
  Fix a critical value $E_c\in\CC_c$.  Let $m$ be the intersection of
  the cycle $\delta_{E_c}$ with $\gamma_j$, and $\Omega_m$ an open set
  in which Proposition \ref{prop:fn-locale-semicla} applies.  We can
  assume that the paths $\delta_E$, $E\in D^\pm$ all entirely lie in
  $\Omega_m$ (using Proposition \ref{prop:regular-holonomy}, this will
  only modify $\int_{\delta_E}\lambda_h$ by an additive term entering
  in $g_{\delta}(E;h)$).  As before, label the local un/stable
  manifolds with cyclic order $I_1$, $I_3$, $I_2$ and $I_4$.
  $\delta_{E_c}$ enters $\Omega_m$ on a local stable manifold $I_s$,
  $s=3,4$, and leaves it on a local unstable manifold $I_u$, $u=1,2$
  ($u$ is the index \emph{preceding} $s$ in the cycle $(1,3,2,4)$ if
  $k^\pm=k^+\in N^+$ and the index \emph{following} $s$ if
  $k^\pm=k^-\in N^-$).  As before, we endow a neighbourhood of each
  $I_\alpha$ with the distribution $e_\alpha$:
  \begin{eqnarray*}
    e_1 & := & u_+^{\e_{2,j}/h}:= 1_{\pm
      y>0}\frac{1}{\sqrt{|y|}}\ex^{i\e_{2,j}\ln |y|/h}, \\
    e_2 & := &  u_-^{\e_{2,j}/h}:= 1_{\pm
      y<0}\frac{1}{\sqrt{|y|}}\ex^{i\e_{2,j}\ln |y|/h}, \\
    e_3 & := & \fourier^{-1}(u_+^{-\e_{2,j}/h}), \\
    e_4 & := & \fourier^{-1}(u_-^{-\e_{2,j}/h}).
  \end{eqnarray*}
  These distributions are classical Lagrangian distributions whose
  phases admit an asymptotic expansion in
  $\frac{1}{h}\Cinf(D^\pm)\formel{h}$. Moreover, they are \mi\ 
  solutions of (\ref{equ:systemE}) in $\Omega_m$ uniformly for $E\in
  D^\pm$, and hence constitute an admissible choice in view of
  Definition \ref{defi:local-holonomy}. Note that this choice possibly
  implies another additive term entering in $g_{\delta_E}$.
  \begin{figure}[htbp]
    \begin{center}
      \input{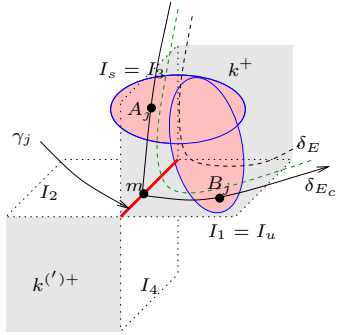}
      \caption{The local holonomy}
      \label{fig:local-holo}
    \end{center}
  \end{figure}
  In a small ball around $A_j(E_c)$, the space of solutions has
  dimension 1 and we must have a constant $C^\pm(E,h)$ such that
  $e_u\sim C^\pm(E,h)e_s$. Now Definition \ref{defi:local-holonomy}
  says that, with respect to the fixed solutions $e_u$ and $e_s$,
  \[
  C^\pm(E,h) = \ex^{-i\int_{\delta_E}\lambda_h} \qquad (\textrm{mod
    }\oh).
  \]
  The expression of the Fourier transform involved in $e_s$ shows --
  as in \cite{colin-p} -- that
  \begin{equation}
    \label{equ:facteur-local}
    C^\pm(E,h)= \zeta^\pm_j \qquad (\textrm{mod } \oh).
  \end{equation}
 
  This proves the existence of the claimed $g_{\delta}$ satisfying
  Eq.(\ref{equ:regularisation}).  The formul\ae\ 
  (\ref{equ:holonomieE0}) and (\ref{equ:holonomieE1}) now follow from
  Stirling's formula (Lemma \ref{lemm:stirling}) and from the fact
  that the asymptotic expansion for $g_{\delta}$ is uniform in
  $D^\pm$, and hence is given for any fixed $E\in U^\pm$ by the
  difference of the (non-uniform) asymptotic expansions of
  $\int_{\delta_E}\lambda_h$ and $i\ln(\zeta_j^\pm(E))$. We have also
  used the formula (\ref{equ:transport}) which holds only for
  non-critical $E$, and which comes from the definition of the bundle
  of principal symbols.
  
  Note that we don't see here any difference between critical circles
  of direct/reverse type, for the statement is purely local near a
  point of $\gamma_j$.
\end{demo}
\begin{demo}[of Theorem \ref{theo:holonomie-graphe}]
  To each $\delta_\ell^{\textrm{loc}}$ is associated a unique
  component $k^{s_\ell}$, $s_\ell=\pm$, a unique critical circle
  $\gamma_{j_\ell}$, and a (non-unique) local family of paths
  $\delta_\ell^{\textrm{loc}}(E)$ (in the sense of Definition
  \ref{defi:deltaE}) such that $\delta_\ell^{\textrm{loc}} =
  \delta_\ell^{\textrm{loc}}(E_c)$. Let $\mathcal{I}_h^{s_\ell}$ be
  the set of $E\in D^{s_\ell}$ such that $(\e_1)_{\restr
    k^{s_\ell}}\in h\ZM + \oh$ (ie. $\mathcal{I}_h^{s_\ell} =
  \varpi_h^{-1}[h\ZM + \oh]$ in the notation of the proof of Corollary
  \ref{coro:cond1E}). Then the (assumed) hypothesis of section
  \ref{sec:abstract} ``$\epsilon_1=n$'' says that $E_c\in
  \mathcal{I}_h^{s_\ell}$. Using Proposition \ref{prop:ansatzE} one
  can construct smooth families $(u_E)_{E\in
    \mathcal{I}_h^{s_\ell}}^{1,2,3,4}$ of solutions on a neighbourhood
  of $\gamma_{j_\ell}$ in $\overline{k^{s_\ell}}$ such that 
  \begin{equation}
    \label{equ:base-std-4}
    (u_E^1,u_E^2,u_E^3,u_E^4)
  \end{equation}
  -- in the direct case -- or
  \begin{equation}
    \left\{
      \label{equ:base-std-2}
      \begin{array}{ccc}
        u_E^1 & + & (-1)^n u_E^2 \\
        \ex^{i\frac{\pi}{2}n}u_E^3 & + & \ex^{-i\frac{\pi}{2}n}u_E^4
      \end{array}
    \right.
  \end{equation}
  -- in the reverse case -- form at $E=E_c$ a \emph{standard basis}
  for the graph $G(E_c)$ at the vertex $g_{j_\ell}$.
    
  Since these solutions are smooth WKB solutions and hence admissible
  in the sense of Definition \ref{defi:local-holonomy}, we shall use
  them to define the local holonomies
  $\int_{\delta_\ell^{\textrm{loc/reg}}(E)} \lambda_h$; since they are
  standard basis at $E=E_c$, we shall in the same way use them to
  define the local ``reduced'' holonomies
  $\hol(\overline{\delta}_\ell^{\textrm{loc/reg}})$. But then by
  definition of the sheaf $(\bar{\L},G(E_c))$ we have
  \begin{equation}
    \label{equ:hol-loc-reg} 
    \left\{
      \begin{array}{l}
        \hol(\overline{\delta}_\ell^{\textrm{loc}}) = 1 \qquad
        \textrm{ and}\\
        \hol(\overline{\delta}_\ell^{\textrm{reg}}) = \exp \left(i
          \int_{\delta_\ell^{\textrm{reg}}(E_c)} \lambda_h\right) =
        \exp \left(i g_{\delta_\ell^{\textrm{reg}}}(E_c)\right).
      \end{array}
    \right.
  \end{equation}
  
  On the other hand, we know from the proof of Theorem
  \ref{theo:local-holonomie} that for such a choice of \mi\ solutions,
  we have (modulo $2\pi$)
  \begin{equation}
    \label{equ:g-loc-reg}
    \left\{
      \begin{array}{ccc}
        g_{\delta_\ell^{\textrm{loc}}}(E) = 0 & \textrm{if
          $\gamma_{j_\ell}$ is of degree 4}; \\
        g_{\delta_\ell^{\textrm{loc}}}(E) = s_\ell\frac{\pi}{2}n & \textrm{if
          $\gamma_{j_\ell}$ is of degree 2}.
      \end{array}
    \right.
  \end{equation}
  Therefore, if we decompose
  \[
  \hol(\overline{\delta}_{E_c}) = \prod_{\ell/\deg\gamma_{j_\ell}=2}
  \hol(\overline{\delta}_\ell^{\textrm{loc}}) \times
  \prod_{\ell/\deg\gamma_{j_\ell}=4}
  \hol(\overline{\delta}_\ell^{\textrm{loc}}) \times \prod_\ell
  \hol(\overline{\delta}_\ell^{\textrm{reg}}),
  \]
  we obtain by (\ref{equ:hol-loc-reg}) and (\ref{equ:g-loc-reg}):
  \begin{eqnarray*}
    \lefteqn{\hol(\overline{\delta}_{E_c}) = \exp
    \left(i\sum_{\ell/\deg\gamma_{j_\ell}=2}
      \!\!\left(g_{\delta_\ell^{\textrm{loc}}}(E_c) -
        s_\ell\frac{\pi}{2}n\right) +{}\right. }\\ 
      & & {}+ \left.\sum_{\ell/\deg\gamma_{j_\ell}=4}
      \!\!g_{\delta_\ell^{\textrm{loc}}}(E_c) + \sum_\ell
      g_{\delta_\ell^{\textrm{reg}}}(E_c)\right),  
  \end{eqnarray*}
  which proves the theorem.
\end{demo}
\begin{coro}
  \label{coro:coeff}
  Theorem \ref{theo:holonomie-graphe} together with formul\ae\ 
  (\ref{equ:holonomieE0}) and (\ref{equ:holonomieE1}) finally prove
  the second point of Proposition \ref{prop:holo}.
\end{coro}
\begin{demo}
  \paragraph{1. The principal action --} Since
  $\e_{2,j}^{(0)}(E_c)=0$, it is clear from (\ref{equ:holonomieE0})
  that for any $\ell=1,\dots,q$,
  \[
  g_{\delta_\ell^{\textrm{loc}}}^{(-1)}(E_c) =
  \int_{\delta_\ell^{\textrm{loc}}} \alpha +
  \Phi_{A_{j_\ell}}^{(-1)}(A_{j_\ell}(E_c)) -
  \Phi_{B_{j_\ell}}^{(-1)}(B_{j_\ell}(E_c))
  \]
  and from (\ref{equ:transport}) that for any $\ell=1,\dots,q$
  (identifying $\ell=q+1$ with $\ell=1$),
  \[
  g_{\delta_\ell^{\textrm{reg}}}^{(-1)}(E_c) =
  \int_{\delta_\ell^{\textrm{reg}}} \alpha +
  \Phi_{B_{j_\ell}}^{(-1)}(B_{j_\ell}(E_c)) -
  \Phi_{A_{j_{\ell+1}}}^{(-1)}(B_{j_{\ell+1}}(E_c)).
  \]
  Therefore,
  \[
  g^{(-1)}(E_c) = \int_{\tilde{\delta}_{E_c}} \alpha.
  \]
  \paragraph{2. The sub-principal action and the Maslov index --} Fix
  $\ell=1,\dots,q$ and let $\gamma_{j}=\gamma_{j_\ell}$. Let
  $m=m(E_c)$ be the point where $\delta_\ell^{\textrm{loc}}$ meets
  $\gamma_j$. As in Proposition \ref{prop:fn-locale-semicla}, we shall
  use the local canonical coordinates at $m$ given by Theorem
  \ref{theo:fn-local}.
  
  Recalling the notation of Proposition \ref{prop:fn-locale-semicla},
  (\ref{equ:fn-localeE}) implies that
  \[
  (\ham{1},\ham{2}) = \left(
    \begin{array}{cc}
      1 & 0 \\ \partial_1\Phi_E & \partial_2\Phi_E
    \end{array}
  \right)\cdot(\ham{\xi},\ham{y\eta}).
  \]
  Since $\partial_2\Phi_E\neq 0$, there exist a smooth function
  $\rho_2=\frac{-\partial_1\Phi_Er_1+r_2}{\partial_2\Phi_E}$
  (depending also smoothly on $E$) such that the sub-principal form
  $\kappa_E$ is given by
  \[
  \kappa_E.(\ham{\xi},\ham{y\eta}) = -(r_1,\rho_2).
  \]
  Note that for a critical value $E=E_c$, $(\rho_2)_{\restr
    y=\eta=0}=-\e^{(1)}_{2,j}$.  The closedness of $\kappa_E$ on each
  $\Lambda_E$ implies that
  \[
  \{r_1,y\eta\} = \{\rho_2,\xi\}.
  \]
  Using a local analogue of Lemma \ref{lemm:poincare}, we can
  decompose $(r_1,\rho_2)$ in the following way:
  \begin{equation}
    \label{equ:decompose-kappa}
    (r_1,\rho_2) = (0,K) -
    (\{\ham{\xi},\tilde{f}\},\{\ham{y\eta},\tilde{f}\})
  \end{equation}
  for some smooth functions $\tilde{f}$, $K$ where $K$ commutes with
  $y\eta$ and $\xi$. Therefore the function
  \[
  f:= \tilde{f} - K\ln |y| \qquad \textrm{ (or } \tilde{f} + K\ln
  |\eta| \textrm{ where } y=0 \textrm{)},
  \]
  restricted to $\Lambda_E$, satisfies
  $d_{\Lambda_E}\tilde{f}=\kappa_E$. We can now compute
  \[
  \int_{\delta_\ell^{\textrm{loc}}} \kappa_E = f(B_j)-f(A_j) =
  \tilde{f}(B_j) - \tilde{f}(A_j) - K\ln |y_{B_j}\eta_{A_j}| +
  K\ln|y_{A_j}\eta_{A_j}|
  \]
  Since $A_j(E)$ and $B_j(E)$ are in $\Lambda_E$,
  $y_{A_j}\eta_{A_j}=\e_{2,j}^{(0)}$ and we have
  \[
  \int_{\delta_\ell^{\textrm{loc}}} \kappa_E + \e_{2,j}^{(1)}\ln
  |\e_{2,j}^{(0)}| = \tilde{f}(B_j) - \tilde{f}(A_j) - K\ln
  |y_{B_j}\eta_{A_j}| + (\e_{2,j}^{(1)}+K)\ln |\e_{2,j}^{(0)}|.
  \]
  Because of (\ref{equ:decompose-kappa}), $\rho_2-K$ vanish at
  $y=\eta=0$, hence $K=-\e_{2,j}^{(1)} + O(y\eta)$ and
  $(\e_{2,j}^{(1)}+K)\ln |\e_{2,j}^{(0)}|$ tends to zero as $E$ tends
  to a critical value $E_c$. Since
  $g_{\delta_\ell^{\textrm{loc}}}^{(0)}$ is smooth at $E_c$, the
  formula (\ref{equ:holonomieE1}) implies that
  $\mu(\delta_\ell^{\textrm{loc}}(E))$ is continuous at $E_c$ and
  hence constant; let us denote it by
  $\mu(\delta_\ell^{\textrm{loc}}(E_c))$.
  
  Suppose now that $E=E_c$ and let $a$, $b$ be points on
  $\delta_\ell^{\textrm{loc}}(E_c)$ located respectively in $[A_j,m]$
  and $[m,B_j]$. Then
  \[
  \tilde{f}(B_j) - \tilde{f}(A_j) - K\ln |y_{B_j}\eta_{A_j}| = 
  \]
  \[ {} \qquad =  \lim_{a,b\fleche m} \left(\tilde{f}(a) - \tilde{f}(A_j) +
    \tilde{f}(B_j) - \tilde{f}(b) - K\ln |y_{B_j}\eta_{A_j}|\right).
  \]
  The term in the limit is equal to
  \[
  \int_{[A_j,a]}\kappa_{E_c} + \int_{[b,B_j]}\kappa_{E_c} +
  \e_{2,j}^{(1)}\ln |y_{b}\eta_{a}|
  \]
  Therefore, by Definition \ref{defi:invariants}, the limit is equal
  to $\int_{\delta_\ell^{\textrm{loc}}}\tilde{\kappa}_{E_c}$, and
  (\ref{equ:holonomieE1}) yields:
  \[
  g_{\delta_\ell^{\textrm{loc}}}^{(0)}(E_c) =
  \int_{\delta_\ell^{\textrm{loc}}}\tilde{\kappa}_{E_c} +
  \mu(\delta_\ell^{\textrm{loc}}(E_c))\frac{\pi}{2}
  -s_\ell\frac{\pi}{4} \quad + \]
  \[
  {} \qquad + \quad \Phi_{A_j}^{(0)}(A_j(E_c)) -
  \Phi_{B_j}^{(0)}(B_j(E_c)).
  \]
  Then as before, if we sum up all the contributions from regular and
  local paths, we obtain
  \[
  g^{(0)}(E_c) = \int_{\tilde{\delta}_{E_c}} \tilde{\kappa}_{E_c} +
  \mu(\tilde{\delta}_{E_c})\frac{\pi}{2} + \sum_\ell
  -s_\ell\frac{\pi}{4}.
  \]
  Using the definition
  \ref{defi:invariants} of the regularized Maslov cocycle, we finally
  obtain
  \[
  g^{(0)}(E_c) + \frac{\pi}{2}(N_2^- - N_2^+) =
  \int_{\tilde{\delta}_{E_c}} \tilde{\kappa}_{E_c} +
  \tilde{\mu}(\tilde{\delta}_{E_c})\frac{\pi}{2},
  \]
  and equation (\ref{equ:holE}) concludes the proof.
\end{demo}

\section{Examples}
We propose in this section several examples for which our theory
applies. Many other could probably be found; the ones presented
here are interesting by their simplicity and yet by their rich
structure and behaviour.

\subsection*{The recipe}
Let us recall here briefly the recipe for obtaining the \semicla\ 
quantisation rules. The first thing to do is to locate the critical
value of transversal hyperbolic type in the image of the momentum map
$F=(H_1,H_2)$. Then choose one of these points $o$ and describe the
singular level set $\Lambda_o=F^{-1}(o)$, in order to have: a) the
graph $G$, b) a formula for the vector fields $\ham{1}$ and $\ham{2}$
on $\Lambda_o$ -- and if it is not one of these, for the periodic
vector field $\ham{p}$.

Compute the \semicla\ invariants (action integral, Maslov index) for a
periodic cycle -- this implies only regular tools -- in order to
derive the first quantisation condition of Theorem \ref{theo:globalQN}
up to $O(h)$, and fix the quantum number $n$.

From the graph, apply Theorem \ref{theo:BSformel} to obtain the second
quantisation rule in the form of a determinantal equation. It remains
to compute the holonomy $\hol$ up to $O(h)$, which involves the
singular \semicla\ invariants of Definition~\ref{defi:invariants}. The
fulfilment of these quantisations rules determine the spectrum up to
an error of order $O(h^2)$ in a window of size $O(h)$ (in fact, it is
easy to determine the smooth dependence of the \semicla\ invariants in
$o$, and the spectral window can be extended to a rectangular domain
of size $O(1)$ along the curve of critical values, and of size $O(h)$
in the transversal direction).

\paragraph{Notation--} The reader must be warned that the symbol $e$
(italic) is used as a subprincipal spectral parameter (as in
``$E=he$''), while the exponential is denoted by $\ex^a=\exp(a)$.

\subsection{Laplacians on Ellipsoids}

Let us consider the ellipsoid in the Euclidian space $\R^3$ defined
by~:
\[
E=\{ \frac{x_1^2}{a_1^2}+ \frac{x_2^2}{a_2^2}+ \frac{x_3^2}{a_3^2}= 1
\} \] with $0<a_1<a_2<a_3 $.  The geodesic flow on $E$ has been
discovered to be integrable by Jacobi in 1838 using Abelian integrals.
For a recent presentation, one can read \cite{moser-quadrics},
\cite{audin-qua}, \cite{zung-al} or \cite{zung-torus}.

\subsubsection{Classics}
Let us denote by $P,Q,P'=-P,Q'=-Q$ the four umbilics of $E$ which are
located on the ellipse $\{ x_2=0 \}$.  If
$$X_1=\sqrt{a_1(a_2-a_1)(a_3-a_1)}{\rm ~and~}X_3=
\sqrt{a_3(a_3-a_2)(a_3-a_1)}~,$$ we have
$$P=(X_1,0,X_3),~Q=(-X_1,0,X_3)~.$$ We will consider the (unique up to
global dilatation) conformal representation $\Phi $ of $E_+=E\cap
\{x_2>0\} $ on a rectangle $R=]0,T_1[\times ]0, T_2[$ such that the
four umbilics are going on the four vertices of $R$, according to
figure \ref{fig:fundamental}.  Using such coordinates $(x,y)\in R$, we
get (see \cite{darboux} (vol. 2 p. 308 and vol. 3 p.13) or
\cite{klingenberg}) the following expression for the metric of $E$~:
\begin{equation}
  ds^2=\big(a^2(x)+b^2(y)\big)(dx^2+dy^2)
\end{equation}
where $a,b$ are given in terms of hyperelliptic integrals and extends
to smooth functions on $\R$ which satisfy: $a$ is $>0$ on $]0,T_1[$,
vanishes exactly at the points $kT_1,~k\in \Z$, and is odd with
respect to $T_1\Z$ and $b$ satisfies the same properties with respect
to $T_2$.  Moreover $a'(0)=b'(0)>0$.  Let us denote by $\Gamma $ the
lattice $T_1\Z \oplus T_2 \Z$.  Then $ds^2$ extends into a smooth
metric on $\R^2 \setminus \Gamma $, which is $\Gamma $-periodic.  Let
us consider the torus $T=\R^2 /2\Gamma $.  Then the map $\gs :T\ra T$
defined by $T(z)=-z$ defines an isometric involution of $T$ with four
fixed points and we get a natural identification of $E$ with $T/\gs $
as a 2-sheeted branched covering $\Pi $ of $T$ over $E$ with
automorphism $\gs$.  The metric $ds^2$ admits conical singular points
of total angle $4\pi$ at the umbilics which makes the metric on $E$
smooth.  More precisely, it follows from the formulae of
\cite{klingenberg} that there exists an analytic function $G$ defined
near $0$, with $G(0)=0,G'(0)>0$ such that near $(0,0)$ we have $
ds^2=(G(x^2)-G(-y^2) )(dx^2+dy^2)$.  It is rather easy to check that
there exists an analytic function $A(u,v)$ with $A(0,0)>0$ such that
$G(x^2)-G(-y^2)=(x^2+y^2) A(x^2-y^2,2xy) $ and, if locally
$Z=\Pi(z)=z^2$, we have~: $ ds^2 =\Pi ^\star \big(4A(X,Y)
(dX^2+dY^2)\big) $.

We will use the fundamental domain $D=[0,2T_1]\times [0, T_2]$.  $E$
can be recovered from this rectangle by gluing edges as indicated on
the figure \ref{fig:fundamental}.

If $a,b$ were non vanishing, $ds^2$ would be called a {\it Liouville
  metric} on $T$. Our case corresponds to a degenerate Liouville
metric on the sphere.  It is well known that Liouville metrics are
integrable.  Let us denote by
$$H_1=\frac{\xi^2+\eta ^2}{a^2(x)+b^2(y)}$$ the geodesic flow and by
$$H_2= \frac{b^2(y)\xi^2-a^2(x)\eta^2}{a^2(x)+b^2(y)}~.$$ The manifold
$L_{E,F}=\{ H_1=E,~H_2=F\}$ is given by:
$$ L_{E,F}=\{ \xi^2=F+a^2(x)E,~\eta^2=b^2(y)E-F\}$$ which is obviously
Lagrangian.

We are interested in the singular value $o=(E=1,F=0)$ of the moment
map $(H_1,H_2)$ and the corresponding $\Lambda_o$.  Geodesics passing
through $P$ (resp. $Q$) contain also $P'$ (resp. $Q'$) and vice-versa.
$\Lambda_o$ is the set of unit covectors corresponding to geodesics
passing through $P$ or $Q$.  We have $\Lambda_o=\cup L_{\pm, \pm}$
where $L_{\pm, \pm}=\{ \xi=\pm a(x),~\eta =\pm b(y) \}$.  In
particular $L_P=L_{+,+}\cup L_{-,-}$ and $L_Q=L_{+,-}\cup L_{-,+}$ are
smooth Lagrangian tori whose intersection is $\gg_+\cup \gg_-$. Here
$\gg_+$ (resp. $\gg_-$) is the lift of the ellipse $x_2=0$ with
orientation $(P,Q,P',Q')$ (resp. opposite).  $L_P$ (resp. $L_Q$) is
the set of unit covectors of geodesics of $E$ containing $P$ and $P'$
(resp. $Q$ and $Q'$).  $L_P$ (resp. $L_Q$) is the stable manifold of
$\gg_-$ (resp. $\gg_+$) and the unstable manifold of $\gg_+$ (resp.
$\gg_-$).
\begin{figure}[hbtp]
  \begin{center}
    \leavevmode \input{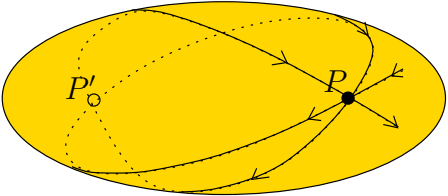}
    \caption{a geodesic passing through $P$ and $P'$}
    \label{fig:geodesiques}
  \end{center}
\end{figure}

The associated graph $G$ is the union of 2 circles corresponding to $L_P$
and $L_Q$ intersecting at 2 points corresponding to $\gg_\pm $.

\begin{figure}[hbtp]
  \begin{center}
    \leavevmode \input{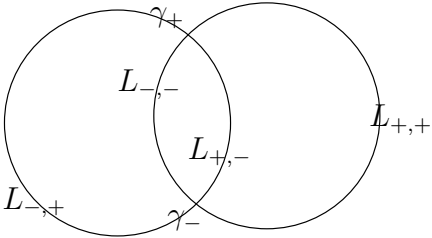}
    \caption{associated graph for $E$}
    \label{fig:graphell}
  \end{center}
\end{figure}

\begin{figure}[hbtp]
  \begin{center}
    \leavevmode \input{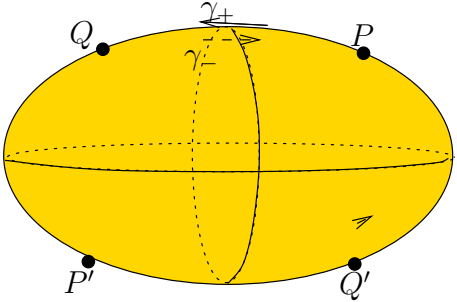}
    \caption{the ellipsoid}
    \label{fig:ellipsoid}
  \end{center}
\end{figure}

\begin{figure}[hbtp]
  \begin{center}
    \leavevmode \input{fundamental.pstex_t}
    \caption{fundamental domain }
    \label{fig:fundamental}
  \end{center}
\end{figure}

\subsubsection{Quantum}

We now introduce the quantum Hamiltonian $\hat{H}_1=h^2 \Delta _E$
which is given in the coordinates $(x,y)$ by:
$$\hat{H}_1=-\frac{h^2}{a^2+b^2}({\pa _x^2+\pa _y^2})$$ and the
operator
$$\hat{H}_2=-\frac{h^2}{a^2+b^2}({b^2\pa _x^2-a^2\pa _y^2})~.$$ It is
possible to check directly that $[\hat{H}_1,\hat{H}_2]=0$ outside of
the lattice $\Gamma $ which is the set of singular points.  We need a
stronger form of commutation, namely if
\begin{equation} \label{equ:eigenlaplace}
  \hat{H}_1\gf =\gl \gf ~,
\end{equation}
we need to prove that $\hat{H}_2 \gf $ is in the domain of
$\hat{H}_1$ and $\hat{H}_1(\hat{H}_2 \gf) =\gl \hat{H}_2 \gf $ ~, so
that $\hat{H}_1$ and $\hat{H}_2$ have a commun eigenbasis.

For that, we need a caracterisation of the functions $\tilde{\gf}=\gf
\circ \Pi$ where $\gf $ is an eigenfunction of the Riemannian
Laplacian $\Delta _E $ on $E$.
\begin{lemm}
  $\tilde{\gf} :T\ra \CM$ is of the form $\tilde{\gf}=\gf \circ \Pi$
  with $\Delta _E \gf =\gl \gf $ if and only if $\tilde{\gf} $ is $\gs$
  invariant and satifies
  \begin{equation}\label{equ:deltaT}
    -\tilde{\gf} ''_{x^2}-\tilde{\gf}''_{y^2}=\gl (a^2+b^2)\tilde{\gf}~.
  \end{equation}

\end{lemm}
\begin{demo} Starting from $\gf $ an eigenfunction of
  $\Delta _E$ and using smoothness of $\Pi$ gives the trivial direction.
  For the other, using the fact that $\tilde{\gf}\circ \gs
  =\tilde{\gf}$, we get a bounded $\gf $ with $\tilde{\gf}=\gf \circ \Pi
  $ and $(\Delta _E -\gl )\gf =T$ where $T$ is supported insides the
  (finite) set of umbilics.  We deduce that $T$ is $0$ or that $\gf $
  (and hence $\tilde{\gf}$) is unbounded.
\end{demo}

It is easy to check that $\hat{H}_2 \tilde{\gf} $ is smooth using
equation (\ref{equ:eigenlaplace}).  Then $\hat{H}_2 \tilde{\gf} $
satisfies equation (\ref{equ:deltaT}) outside $\Gamma $ and hence
everywhere on $T$ and we get a commun eigenbasis for $\hat{H}_1$ and
$\hat{H}_2$.

We prefer to rewrite the eigenvectors equations
$$\hat{H}_1\tilde{\gf} =\gl \tilde{\gf} ,~ \hat{H}_2 \tilde{\gf} =\mu
\tilde{\gf}$$ in the following simpler way:
$$ \hat{P}\tilde{\gf} := h^2\frac{\pa ^2 \tilde{\gf} }{\pa
  x^2}+(a^2(x) \gl +\gm)\tilde{\gf} =0$$
$$ \hat{Q}\tilde{\gf} := h^2\frac{\pa ^2 \tilde{\gf} }{\pa
  y^2}+(b^2(y) \gl -\gm)\tilde{\gf} =0~.$$ We are interested in
solutions of this system which are $\gs $ invariant. If we denote by
$\gs_1(x,y)=(-x,y)$ and $\gs _2(x,y)=(x,-y)$, we get $\gs =\gs _1
\circ \gs _2$ and because $\gs _j$ commutes with $\hat{P}$ and
$\hat{Q}$ we are reduced to find solutions of the form $\tilde{\gf}
(x,y)=f(x)g(y)$ with $f$ a $2T_1$-periodic solution of $\hat{P}f=0$
and $g$ a $2T_2$-periodic solution of $\hat{Q}f=0$.  We ask moreover
that $f$ and $g$ are both even or both odd.  We assume $\gl =1$ which
corresponds to quantize $h$ and $\gm=\ge h $. The associated fiber of
the momentum map is then $\Lambda_o$.

This way we are reduced to 2 one dimensional problems and because
$\hat{P} $ and $\hat{Q}$ are semi-classical stationnary Schrödinger
operators with potentials $-a^2$ and $-b^2$, we are reduced to the
computations of \cite[p. 489-490]{colin-p3} for periodic double wells.


\newpage

\subsection{$1:2$-resonance}

\subsubsection{Birkhoff normal forms}

Consider a Hamiltonian $H:T^*\RM^2 \ra \RM$ with a non-degenerate
minimum at the origin. We can assume using a symplectic linear change
that $H(z_1,z_2)=K_2(z)+O(|z|^3)$ with
\[
 K_2(z)=\go _1|z_1|^2+\go _2 |z_2|^2~,
\]
and $\go_j>0$.  Here $z_j=x_j+i\xi_j$, where $(x_1,x_2,\xi_1,\xi_2)$
are canonical coordinates for $T^*\RM^2$. We will say that the
quadratic part is resonant if $\go_1/\go_2$ is a rational number.  It
is possible to derive a Birkhoff normal form $H$ of the following form
\[
H=K_2+R +O(|z|)^\infty
\] 
with $R=O(|z|^3)$ and $\{ K_2, R\} =0$.  The same result is true on
the quantum level (see chapter 5 of \cite{san-these}) with commuting
operators $\hat{K}_2$ and $ \hat{R}$ .  If we are able to analyse the
joint spectrum of the operators $\hat{K}_2, \hat{R}$ we can deduce
some sharp results for eigenstates in the energy domain $E=O(h^\ga )$
with $\ga >0$.  In the case of the $1:1$ resonance -- ie.
$\omega_1=\omega_2$ -- the flow of $K_2$ induces a free circle action
on the energy hypersurface $K_2=\textrm{const}$ and the reduced space
is smooth. Then, via the use of Toeplitz operators, the problem is
fully reduced to a 1-dimensional one. This is no longer the case for
the $1:2$ resonance, where the reduced phase space has a conical
singularity. For this simple example, we will show that our analysis
applies. Another application would be the \emph{near} $1:2$ resonance
with $R=\varepsilon|z_1|^2+R'$.

\subsubsection{$1:2$ resonance}

\begin{figure}[htbp]
  \begin{center}
    \input{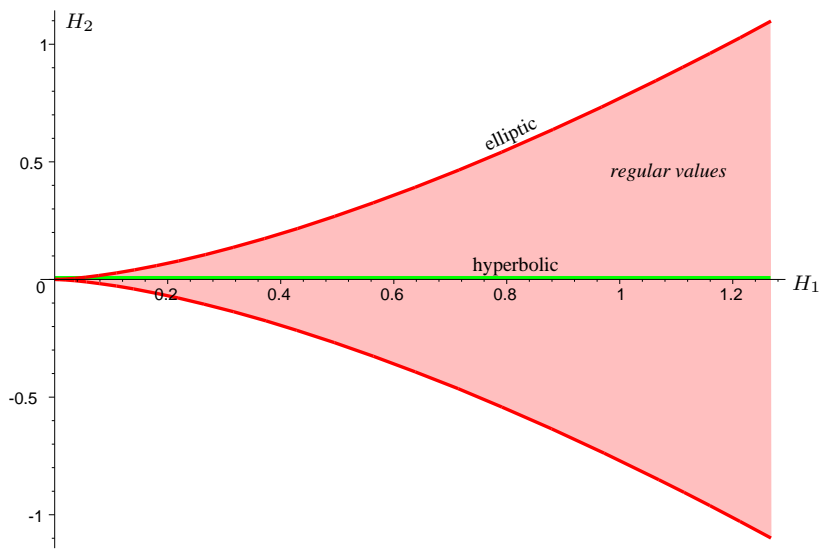}
    \caption{Image of the momentum map for the $1:2$ resonance}
    \label{fig:moment-reso12}
  \end{center}
\end{figure}

We consider the following Poisson commuting Hamiltonians
on $T^\star \R^2$:
\begin{equation}\label{equ:hami-reso12}
H_1=\ha |z_1|^2+|z_2|^2, \quad
H_2=(x_1^2-\xi_1^2)x_2+2x_1\xi_1\xi_2=\Re (z_1^2 \bar{z}_2)
\end{equation}
with $z_j=x_j+i\xi_j,~j=1,2$. The image of the momentum map
$F=(H_1,H_2)$ is
\[
F(T^*\RM^2) = \{(X,Y) | \quad 16X^3 \geq 27Y^2\}.
\]
The singular values consist of the boundary (which corresponds to
transversal elliptic points, except for the origin which is
degenerate) and the half line $\CC_c=\{(X,0)$, $X>0\}$, whose points
are transversally hyperbolic (see Fig.\ref{fig:moment-reso12}).
Here we shall be interested in the critical values on $\CC_c$.
Because of the homogeneity of $H_j$, it is sufficient to consider the
point $o=(1,0)$.

The corresponding commuting quantum Hamiltonians are:
\begin{equation}
\label{equ:qhami1-reso12}
\hat{H}_1=\ha (-h^2\frac{\pa ^2}{\pa x_1^2}+x_1^2)+
(-h^2\frac{\pa ^2}{\pa x_2^2}+x_2^2),
\end{equation}
\begin{equation}
\label{equ:qhami2-reso12}
\hat{H}_2=x_2(h^2\frac{\pa ^2}{\pa x_1^2}+x_1^2)-h^2\frac{\pa}{\pa x_2}
(2x_1\frac{\pa}{\pa x_1}+1).
\end{equation}

\subsubsection{Classical description}
Here we are interested in the singular Lagrangian leaf $\Lambda_o$
defined by $H_1=1, H_2=0$. The singular part of $\Lambda_o$ is the
closed trajectory $\gamma_0=\{z_1=0\} \cap \{ |z_2|=1 \}$.  From its
defining equations, it is easy to find a parameterisation that shows
that $\Lambda_o$ is a Lagrangian immersion of a Klein bottle
$\mathbb{K}$ with $\gamma_0$ as a double loop:
\begin{equation} \label{equ:klein}
\Phi: \mathbb{K}\ni(\theta ,\varphi) \mapsto (\sqrt{2}\ex^{i\theta }\sin
\varphi , -i\ex^{2i\theta }\cos \varphi ) \in \Lambda_o,
\end{equation}
where $\mathbb{K}$ is the quotient of
$\T^2_{(\theta,\varphi)}=\RM^2/\ZM^2$ by the equivalence relation
\[
(\theta +\pi,-\varphi) \sim (\theta,\varphi)
\]

\begin{figure}[hbtp]
  \begin{center}
    \leavevmode
    \input{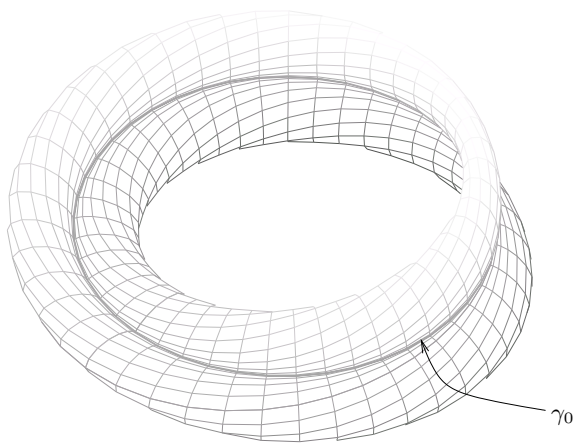}
    \caption{The manifold $\Lambda_o$}
    \label{fig:klein }
  \end{center}
\end{figure}

A fundamental domain $D$ is given by $D=\{ (\theta ,\varphi )\,|\quad
0\leq \theta \leq \pi ,~-\pi \leq \varphi \leq \pi \}$.
\begin{figure}[hbtp]
  \begin{center}
    \leavevmode
    \input{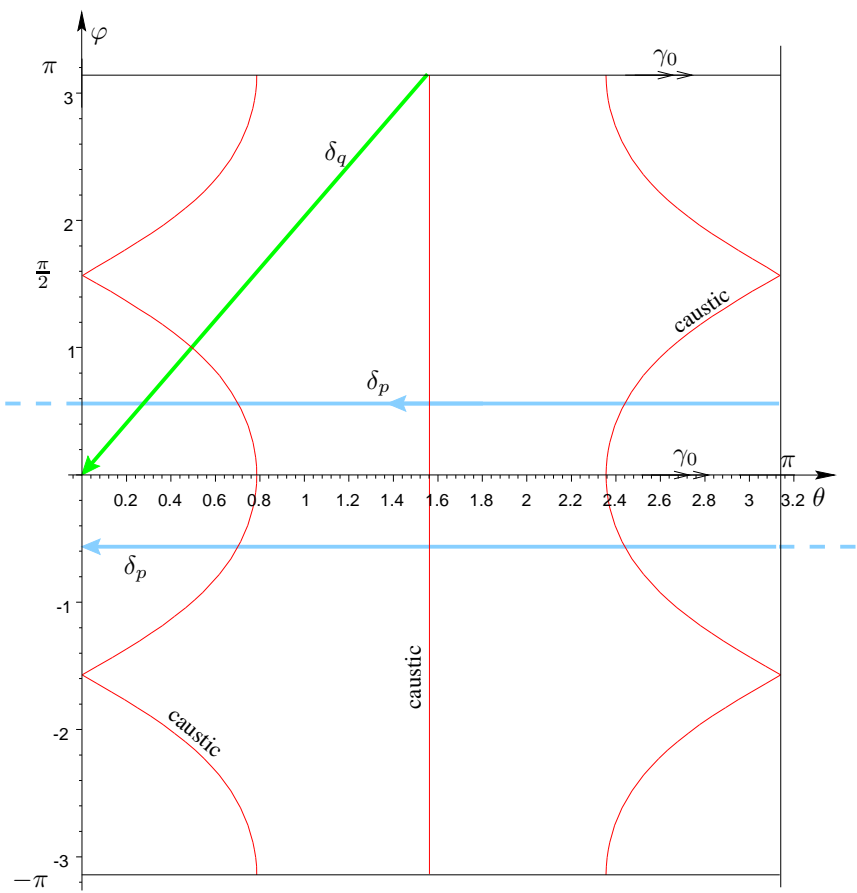}
    \caption{Parameters set for $\Lambda_o$}
    \label{fig:kleinpara}
  \end{center}
\end{figure}
The singular line $\gg_0$ corresponds to $\{ \varphi =0 \}\cup
\{\varphi =\pm \pi \}$ and we have there the identifications $\Phi
(\theta , 0)=\Phi (\theta +\pi /2,\pm \pi )$.  The graph
$G=\epsfig{file=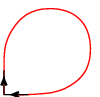}$ corresponding to $\Lambda_o$ has
just one vertex $\gamma_0$ and one edge.

\begin{rema}
  Although we don't really need it, it can be helpful to have a
  representation of the reduced phase space $W=H_1^{-1}(1)/S^1$ --
  where the $S^1$-action is the flow of the harmonic oscillator $H_1$.
  Using a priori argument, one can show that $W$ is a 2-sphere with a
  conical singular point; however, one can find an explicit equation
  for $W$. The algebra of $S^1$-invariant polynomials -- that is,
  those that commute with $H_1$ -- is generated by
  \[
  \pi_1=|z_1|^2,\quad \pi_2=|z_2|^2, \quad \pi_3=\re(z_1^2\bar{z}_2),
  \quad \pi_4=\im(z_1^2\bar{z}_2),
  \]
  which are subject to the relation $\pi_3^2+\pi_4^2=\pi_1^2\pi_2$.
  One can show (see the book \cite{cushman-book}) that this relation
  together with $\pi_1\geq 0$ and $\pi_2\geq 0$ define the orbit space
  $T^*\RM^2/S^1$ in terms of the variables $\pi_j$. The energy level
  set is the section $\{\pi_1\geq 0\}\cap\{\pi_2\geq 0\}$ of the
  3-dimensional hyperplane $\pi_1+2\pi_2=1$. Therefore, $W$ is defined
  in the space $\RM^3=(\pi_1,\pi_4,\pi_3)$ by the equation
  \[
  \pi_3^2 = \pi_1^2(1-\pi_1)/2 - \pi_4^2, \quad \textrm{ with }
  \pi_1\in[0,1].
  \]
  $W$ is a surface of revolution around the $\pi_1$-axis, homeomorphic
  to a 2-sphere, with a conical singularity at the origin.  Note that
  $H_2=\pi_3$ so that the restriction $(H_2)_{\restr W}$ is just the
  height function (see Fig.~\ref{fig:poire}) and is a Morse function
  on $W\setminus\{0\}$. The manifold $\Lambda_o$
  reduces to the singular equator $\pi_3=0$.
  \begin{figure}[htbp]
    \begin{center}
      \input{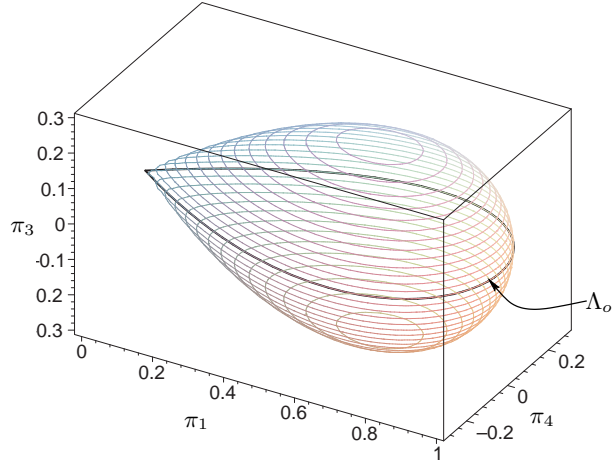}
      \caption{The singular reduced phase space $W$. Here the function
        $(H_2)_{\restr W}$ is equal to the $\pi_3$ coordinate.}
      \label{fig:poire}
    \end{center}
  \end{figure}
\end{rema}

\subsubsection{Semi-classical computations}
We consider the solutions of the system
\begin{equation} \label{equ:systq-reso12}
  (\H_1-1)u=0, \quad (\H_2-eh)u=0,
\end{equation}
for bounded $e$. The microsupport of the solutions is $\Lambda_o$.

If we denote by ${\cal X}_j$ the Hamiltonian
vector fields of $H_j$ we get on $\Lambda_o$:
\begin{equation}
{\cal X}_1=-\frac{\pa}{\pa \theta },~
{\cal X}_2=-2\sin \varphi \frac{\pa}{\pa \varphi  }~,
\end{equation}
which leads to the following sub-principal form $\kappa$ for
$\F=(\H_1,\H_2-eh)$:
\begin{equation}
\kappa =\frac{-e}{2\sin \varphi }d\varphi~.
\end{equation}
Note also that, since the flow of $\ham{1}$ is $2\pi$-periodic outside
$\gamma_0$, we have $\ham{p}=\ham{1}$.  The canonical $1$-form $\alpha
=\xi_1dx_1+\xi_2dx_2$ is given by:
$$\ga= -2(\sin ^2 \theta \sin ^2 \varphi +\cos ^2 2\theta
\cos ^2 \varphi)d\theta +\ha (\sin 2 \theta (1+\cos 2\theta)
\sin 2 \varphi )d\varphi ~.$$
Finally, the caustic set $C$ of $\Lambda_o$ is given by
\[
 C=\{ \cos \theta =0\} \cup \{ \tan^2 \theta = \cos^2 \varphi \} .
\]

In order to compute the quantisation rules, let us
introduce the following loops on $\Lambda_o$:
\begin{eqnarray}
\delta_p(\pi-s)=\Phi(-s,\pm \varphi_0), & s\in[0,\pi] & \quad \textrm{ for
  some } \varphi_0\neq 0 (\pi), \\
\delta_q(\frac{\pi}{2}-s)=\Phi (-s,-2s), & s\in [0,\pi/2] & .
\label{equ:loops-reso12}
\end{eqnarray}
$\delta_p$ is an oriented $S^1$-orbit, and $\delta_q$ is a loop which
is everywhere transversal to the $S^1$-action and oriented according
to the flow of $\ham{2}$.

The first quantisation condition:
\[
\frac{1}{2\pi}\int_{\delta_p}\alpha + h\mu(\delta_p)/4 \in h\ZM + O(h^2)
\]
is actually exact since $\H_1$ is a harmonic oscillator, and therefore
reads:
\begin{equation}
  \label{equ:cond1-reso12}
  1-\frac{6}{4}h=hn \quad \textrm{ or }\  \quad h=\frac{1}{n+3/2}.
\end{equation}
Because of the homogeneity property, we choose here to view this
condition as a discrete quantisation of $h$ (see remark
\ref{rema:quantif-h}).

Assuming now that (\ref{equ:cond1-reso12}) holds, we can compute the
\semicla\ invariants associated to $\delta_q$:

The action integral is easily computed to be
\[
\int_{\delta_q}\alpha = \frac{\pi}{2},
\] 
and the invariant $\varepsilon$ given by (\ref{equ:epsilon2}) is equal
to $e/2$. The sub-principal action $I_{\delta_q}$ of Definition
\ref{defi:invariants} is given by
\[
I_{\delta_q} := \int_{\delta_q}\tilde{\kappa} = 3\varepsilon\ln 2.
\]
Finally, we can show by slightly shifting $\delta_q$ to the right (in
the $\theta$ direction) that its Maslov index is $-2$. Moreover
$\delta_q$ turns around $\gamma_0$ in the direct sense, hence the
regularised Maslov index is $-2+(\frac{1}{2}+n)$.

We can now write down the second quantisation condition:
$\mathcal{C}_0=\hol(\delta_q)$, which reads:
\begin{eqnarray}
  \label{equ:C0-reso12}
  \ex^{-i\frac{\pi}{4}-in\frac{\pi}{2}}
  \left(1+i(-1)^n\ex^{-\varepsilon\pi}\right)
  \Gamma\left(\frac{1}{2}+i\varepsilon\right)
  \ex^{\varepsilon\left(\frac{\pi}{2}+i\ln h\right)} = \\
  \label{equ:hol-reso12}
  = \ex^{i\left(\frac{\pi}{2h} +
      3\varepsilon\ln 2 + \frac{\pi}{2}(-2+(\frac{1}{2}+n)) +
      O(h)\right)}.
\end{eqnarray}
Using (\ref{equ:cond1-reso12}) and $\varepsilon=e/2$, we obtain the
equation in $e$ and $n$:
\begin{equation}
  \label{equ:cond2-reso12}
  \left(1+i(-1)^n\ex^{-\frac{e\pi}{2}}\right)
\Gamma\left(\frac{1}{2}+\frac{ie}{2}\right)
\ex^{\frac{e}{2}\left(\frac{\pi}{2}-i\ln (n+3/2)\right)} =
\ex^{i\left(\frac{3}{2}e\ln 2 -\frac{\pi}{2}n + \frac{\pi}{4} +
  O(h)\right)}.
\end{equation}

\begin{rema}
  The \semicla\ invariants were computed explicitly; this is related
  to the fact that $\Lambda_o$ admits a parameterisation as a
  ``rational'' variety. Somewhat paradoxically, it would be more
  technical to compute the WKB invariants attached to \emph{regular}
  tori, since no rational parameterisation of these tori exists. On
  the other hand, the regular invariants can be asymptotically
  recovered from the singular ones using Stirling's
  formula~\eqref{equ:regularisation}.
\end{rema}
\begin{rema}
  The obtained formula~\eqref{equ:cond2-reso12} yields easily (exactly
  as in \cite{san-focus}) the fact that the level spacings for the
  eigenvalues in a region of bounded $e$ are of order $O(1/n\ln n)$ =
  $O(h/\ln h)$ -- while they are of order $O(h)$ in a regular region.
  Moreover, the precise shape of the spacing function is readily
  derived and involves the log-derivative of the Gamma function.
  \label{rema:ecarts}
\end{rema}

\subsubsection{Matrix form for $\hat{H}_2$}
The goal here is to study the restriction of $\H_2$ to the eigenspace
$\mathfrak{E}_n$ of $\H_1$ corresponding to the quantum number $n$
(ie. to the eigenvalue $h(n+3/2)$).

In analogy with formul\ae\ (\ref{equ:systq-reso12}), the operators
(\ref{equ:qhami1-reso12}) and (\ref{equ:qhami2-reso12}) can be written
\begin{equation}
  \label{equ:qhami1bis-reso12}
  \hat{H}_1=h\left(a_1(h)b_1(h) + 2a_2(h)b_2(h) -\frac{3}{2}\right),
\end{equation}
\begin{equation}
  \label{equ:qhami2bis-reso12}
  \hat{H}_2=\sqrt{2}h^{\frac{3}{2}}\left(a_2(h)b_1(h)^2 +
  a_1(h)^2b_2(h)\right),
\end{equation}
with $a_1(h)=(2h)^{-1/2}(h\deriv{}{x}+x)$ and
$b_1(h)=a_1(h)^*=(2h)^{-1/2}(-h\deriv{}{x}+x)$ (and similarly for
$a_2(h)$ and $b_2(h)$ with the variable $y$).  Using the unitary
transform in $L^2(\RM^2)$ : $f(x)\fleche \sqrt{h}f(\sqrt{h}x)$, the
operators $\H_1$ and $\H_2$ are transformed into those given by the
equations (\ref{equ:qhami1bis-reso12}) and
(\ref{equ:qhami2bis-reso12}) with $a_j(h)$ and $b_j(h)$ replaced by
$a_j:=a_j(1)$ and $b_j:=b_j(1)$. Note that this shows that the
homogeneity argument used for the classical analysis has an analogue
in the quantum setting: if we know the spectrum for some value of
$h>0$, then the spectrum for any other value of $h$ immediately
follows.

Now, using the Bargmann representation, we identify $a_j$
(respectively $b_j$) with the operator $\deriv{}{z_j}$ (resp. $z_j$),
and let them act on the monomials $\frac{z_1^kz_2^\ell}{k!\ell!}$
which form a Hilbert basis of eigenvectors of $\H_1$ (corresponding to
the eigenvalues $E_1=h(k+2\ell+3/2)$). Then it is easy to find the
matrix of $\H_2$ in this basis of $\mathfrak{E}_n$ ($n=k+2\ell$)~:
\begin{equation}
\label{equ:matrix-reso12}
\mbox{$\H_2$}_{\restr \mathfrak{E}_n} = \sqrt{2}h^{\frac{3}{2}}\left(
  \vcenter{\xymatrix@!C@!R@=2pt@M=2pt{
      0 & A_{n,1} & & & &  \\
      A_{n,1}  & 0 \ar@{--}[ddddrrrr]& A_{n,2}\ar@{.}[dddrrr] & & 0 &  \\
      & A_{n,2}\ar@{.}[dddrrr] & & & &\\
      & & & & & \\
      & 0 & & & & \\
      & & & & & 0 }} \right),
\end{equation}
\[
\textrm{with } A_{n,\ell}=\sqrt{\ell(n-2\ell+1)(n-2\ell+2)}, \qquad
\ell=1,2,\dots E[\frac{n}{2}].
\]

\subsubsection{Numerical computations}
Since $h$ is of order $1/n$, one sees that the coefficients of
(\ref{equ:matrix-reso12}) are bounded as $n\fleche\infty$. Moreover,
since no coefficient $A_{n,\ell}$ vanishes, the spectrum is simple.
For these reasons, it is reasonable to expect a good accuracy of
numerically computed eigenvalues. The resulting spectrum will be
called the ``quantum'' spectrum.

On the other hand, numerically solving equation
(\ref{equ:cond2-reso12}) in the variable $e$ -- assuming that $he$
remains in the bounds of the momentum map -- yields the so-called
``semi-classical'' spectrum for $E_1=1$ and $h=(n+3/2)^{-1}$. If we
wish now to fix $h$ and compute the rest of the joint spectrum, the
same formul\ae\ (\ref{equ:C0-reso12}) and~\eqref{equ:hol-reso12} can
be used if one lets $\varepsilon=\frac{e}{2\sqrt{E_1}}$ and replaces
$h$ by $\tilde{h}=h/E_1$.

The results are displayed in the following figures.  In
Fig.~\ref{fig:reso12-numerics} we have superposed the quantum and the
semi-classical joint spectra. The differences are hardly noticeable
(they are theoretically of order $h^{1}$ -- which means $h^{2}$ for
the unscaled spectrum -- and experimentally much better at the
critical value $H_2=0$), even for a very large $h$ and very small
$E_1$'s -- both of these conditions are supposed to reach the
limitations of our analysis. 
\begin{figure}[htbp]
  \begin{center}
    \input{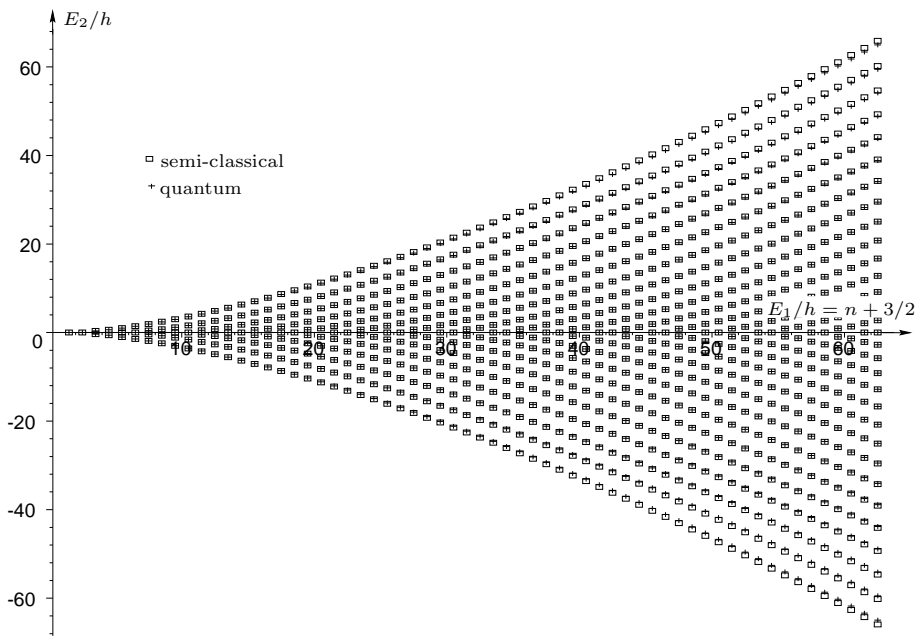}
    \caption{A comparison between semi-classical and quantum
      results. Here $h=2/63$ (so that $E_1=1$ corresponds to $n=30$)}
    \label{fig:reso12-numerics}
  \end{center}
\end{figure}

In the other figures we focus on one spectrum (here at $E_1=1$) around
the critical value $E_2=0$ -- which is the most interesting feature.
\begin{figure}[htbp]
  \begin{center}
    \input{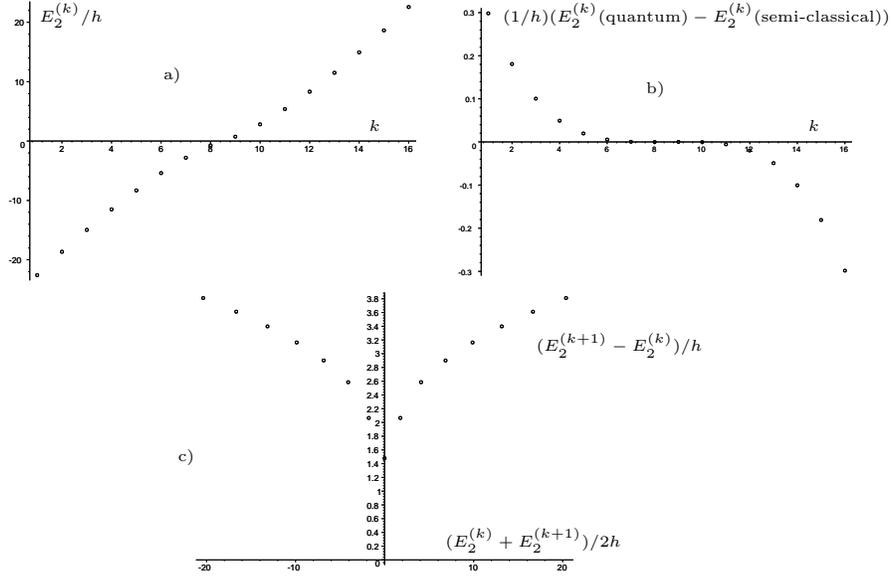}
    \caption{The spectrum at $E_1=1$, $n=30$ (quantum and \semicla\ are
      indistinguishable). a) The spectrum sorted in increasing order
      and displayed versus the eigenvalue number. b) The difference
      ``quantum-\semicla''. c) The energy spacings.}
    \label{fig:reso12-ecarts} 
  \end{center} 
\end{figure}

\begin{figure}[htbp]
  \begin{center} 
    \input{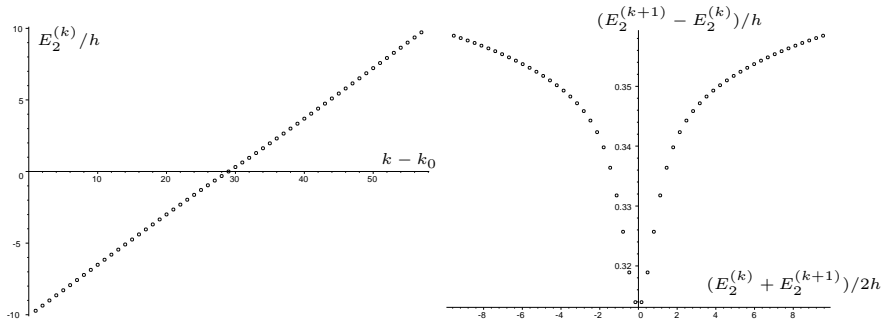}
    \caption{Contrary to the numerical diagonalisation of the
      matrix~\eqref{equ:matrix-reso12}, the \semicla\ formula allows
      very small values of $h$ -- and the results are supposedly even
      more accurate. Here is displayed a window of size $[-10h,10h]$
      for the spectrum and the corresponding eigenvalue spacings,
      where we have let $E_1=1$ and $n=10^{15}$ ($h\simeq 10^{-15}$ !!
      of course we haven't tried the matrix diagonalisation !)}
    \label{fig:fun} 
  \end{center}
\end{figure} 


\newpage

\newcommand{\harm}{\mathcal{H}}
\newcommand{\poly}{\mathcal{P}}
\newcommand{\Bas}{\mathfrak{B}}
\subsection{Schrödinger Operators  on $S^2$}

\subsubsection{Setting of the problem}
We consider now the operator $\hat{H}=\Delta + V$ where $\Delta $ is
the canonical Laplacian on $S^2=\{ (x,y,z)\in \R^3| x^2+y^2+z^2=1 \} $
whose spectral theory is given by the spherical harmonics and
$V:S^2\ra \R$ is a smooth potential.  We introduce the
pseudo-differential operator $\hat{H}_2$ on $S^2$ which is obtained by
averaging $V$ using the $2\pi$-periodic quantum unitary flow
$U(t)={\rm exp}{(it \sqrt{\Delta +1/4})}$~:
\begin{defi}
  $$\hat{H}_2=\frac{1}{2\pi}\int_0^{2\pi}U(-t)VU(t)dt ~.$$
\end{defi}
The following results have been obtained by Weinstein and Guillemin
(see \cite{weinstein-cluster}, \cite{guillemin-sphere},
\cite{guillemin-band}, \cite{guillemin-rank1} and also
\cite{colin-bica} \cite{colinII}):
\begin{theo}
  \begin{itemize}
  \item $\hat{H}_2$ commute with $\Delta$.
  \item $\hat{H}_2$ is a PDO of order $0$ whose principal symbol is the
    Radon transform of $V$:
    $$ H_2(z)=\frac{1}{2\pi}\int_0^{2\pi}V(\gf _t (z))dt~,$$ where $\gf _t
    $ is the geodesic flow with unit speed.  The sub-principal symbol of
    $\hat{H}_2$ vanishes.
  \item There exists an unitary FIO $\Omega$ such that
    $$ \Omega^{-1}\hat{H}\Omega =\Delta +\hat{H}_2 +R$$ where $R$ is a PDO
    of order $-2$.
  \item
    $$ \hat{H}_2=\oplus _{l=0}^\infty \Pi_l V \Pi _l $$ where the $\Pi
    _l$'s are the orthogonal projections on the spaces $\harm_l $ of
    spherical harmonics of degree $l$.
  \end{itemize}
\end{theo}
Proofs can be found in \cite{guillemin-sphere} and
\cite{guillemin-band}.

In such a way, we get a quantum integrable system
$\hat{H}_1=h^2\Delta,~\hat{H}_2$.  The spectrum of $\hat{H}$ is
related to the joint spectrum
$$(h^2l(l+1),\mu_{l,m}),~l=0,\cdots,\infty,~-l \leq m \leq l$$ of
$(\hat{H}_1, \hat{H}_2)$ by $\gl_{l,m}=l(l+1)+\gm_{l,m}+O(l^{-2})$ and
high energy asymptotics ($l\ra \infty $) corresponds in the usual way
to semi-classical asymptotics $h^2.l(l+1)=1,~h\ra 0$.  We will study
the system~:
\[ \hat{H}_1\gf =h^2l(l+1)\gf=\gf ,~\hat{H}_2 \gf =eh \gf ~,\]
assuming that $0$ is a critical value of saddle type of $H_2$.

The Radon transform $H_2$ of $V$ is a function on the manifold $Geod $
of oriented closed geodesics of $S^2$. $Geod $ is a global Poincaré
section for $H_1$ and can be identified with $S^2\subset \R^3_{X,Y,Z}$
by associating to the circle $t \ra \gg(t)=u\cos t +v\sin t $ the unit
vector $u\wedge v $.  Then reversing the orientation of $\gg$
corresponds to antipodal symmetry $\gs $ on $Geod=S^2$ and $H_2$ is
even with respect to that symmetry.  We can then interpret $H_2$ as a
function on the projective plane.  This fact implies that, if $H_2$ is
a Morse-Bott function, it cannot have only local maxima and minima: it
has always saddle points for which our analysis is needed.

We will from now assume that we are in the simplest situation where
$H_2:Geod \ra \R$ has only 2 maxima, 2 minima and 2 nondegenerate
saddle points.  The singular manifold $\Lambda_0$ is then the union of
2 tori which intersect along 2 circles.  The projection of $\Lambda
_0$ on $Geod =S^2$, i.e. the reduction of $\Lambda _0$, is the union
of 2 circles which are invariant by $\gs$ and which intersect at 2
antipodal points.

\begin{figure}[hbtp]
  \begin{center}
    \leavevmode \input{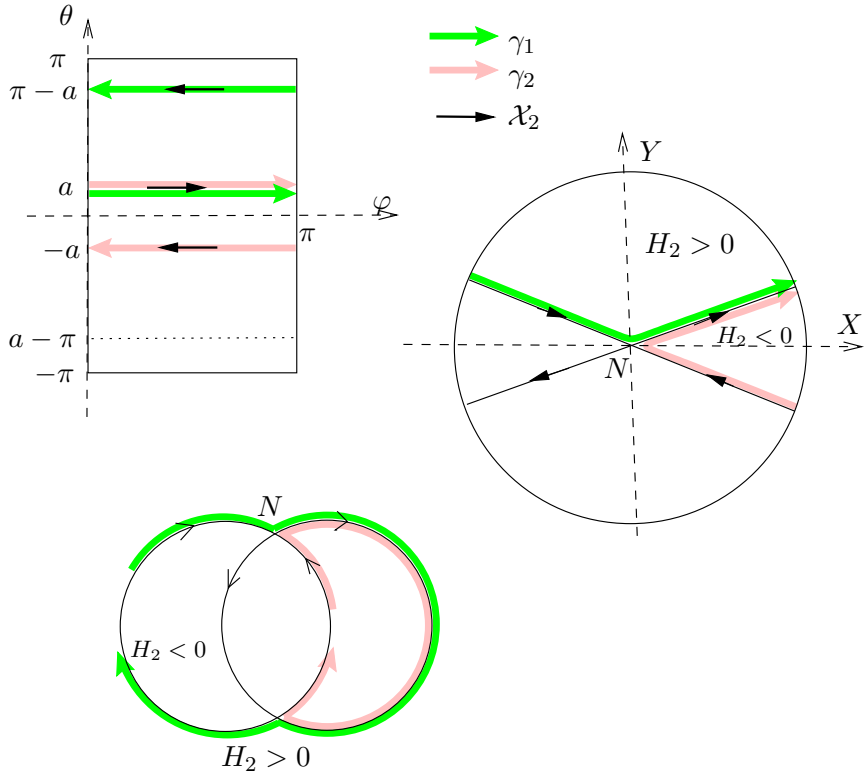}
    \caption{$G$ in $Geod$}
    \label{fig:geod}
  \end{center}
\end{figure}

\subsubsection{Semi-classical computations}

We will complete the computations in the simplest case where $V$
itself is a generic harmonic polynomial of degree 2.  Using $SO(3)$
invariance, we have only to consider the 2-parameter family given by
$$V_{a,b,c}(x,y,z)=2(ax^2+by^2+cz^2)~,$$ with $a+b+c=0$ and $a<c<b$.

Because the Radon transform commutes with the $SO(3)$ action, by Schur's
lemma, $H_2$ itself is an harmonic polynomial of degree 2 on $Geod $
which is given by
$$H_2=aX^2+bY^2+cZ^2~. $$ The critical values of $H_2$ are $a<c<b$
and $H_2-c=B^2 Y^2-A^2 X^2$ with $A=\sqrt{c-a}$ and $B=\sqrt{b-c}$.
It is clear that the singular value $(1,c)$ of $(H_1, H_2)$ is of
hyperbolic type, so that we can apply the previous tools.

Let us denote by $G$ the projection of $\Lambda _0$ on $Geod $ by the
map $\pi$ which associates to a point $m \in U^\star S^2$ the geodesic
to which $m$ belongs; $G$ is the graph introduced in the general
situation and is the union of the two circles $ C_\tau = \{ AX= \tau
BY,~\tau=\pm 1 \}$.  We can compute the projections on $G$ of the
vector field $\ham{2}$ and of the sub-principal form $\kappa $ on $\Lambda
_0$ (because $\kappa ({\cal X}_1)=0$).  We will use spherical
coordinates $(\theta , \gf), 0\leq \theta \leq 2\pi, ~0\leq \gf \leq
\pi, $ on $Geod$:
\[ X= \sin \gf \cos \theta ,~Y=\sin \gf \sin \theta,~Z=\cos \gf ~.\]
The symplectic form on $Geod$ is $SO(3)$ invariant and of total area
$4\pi$. We will assume that:
$$ \go =\sin \gf~d\theta \wedge d\gf ~.$$ By direct computation we get
on $C_\tau $:
$$ \ham{2}=2\tau AB \sin \gf \frac{\pa }{\pa \gf}$$ and
$$\kappa = \frac{\tau e}{2AB \sin \gf }d\gf ~.$$ Let us denote by
$\gg_j,~j=1,\cdots , 4, $ the four cycles of $G$ oriented by $\ham{2}$,
consisting of the union of one arc of $C_1$ and one arc of $C_{-1}$
and bounding a topological disk in $Geod$. $\gg_1$ and $\gg_2$ are
defined on figure \ref{fig:geod} and $\gg_3=\gs (\gg_1),~\gg_4=\gs
(\gg_2)$.  It is then easy to check using the explicit formula for the
integral $\int d\gf /\sin \gf $ that
\begin{equation}\forall j=1,\cdots ,4,~
  I_{\gg_j}=\frac{e}{AB} \log \frac{8AB}{A^2+B^2}~.
\end{equation}
Moreover, we find easily, using formula~\eqref{equ:epsilon2}, that
\begin{equation}
  \ge_2=\frac{e}{2AB} ~.\end{equation}

We put $\ga =\atan \frac{A}{B}$ and $\gb =\atan \frac{B}{A}=\frac{\pi
  }{2} -\ga $.  Action integrals are $A _1 =A_3=- 4 \gb $ and $A_2=A_4
= 4 \ga $.

Let $\tau_j =\hol(\gg _j )$.  We observe first, using the fact that
$h(l+\ha)=1$, the relations $\tau_1 =\tau_3$ and $\tau_2=\tau_4$
(modulo $O(h)$).  We get
\begin{equation}
\tau_1=\ex^{i(  -(4l+2) \gb +I+\pi   +O(h))}~.
\end{equation}
and
\begin{equation}
\tau_2=\ex^{i(  (4l+2)  \ga +I  +O(h))}~.
\end{equation}
It follows that $H:=\tau_1= \tau_2=\tau_3=\tau_4$.  It would be nice
to prove that this relation holds mod $O(h^\infty )$.

Using the computations of \cite{colin-p3} p. 493 and putting
$T=T(\ge)$ with $\ge=\ge _2+ O(h)$ we get the following quantisation
rule:
\begin{equation}
  {\rm det}\bigg({\rm Id}-T
  \left( 
    \begin{array}{cc}
      0&1 \\
      1&0
    \end{array}
  \right)
  T
  \left( 
    \begin{array}{cc}
      0& H^{-1} \\
      H^{-1} &0
    \end{array}
  \right)\bigg)=0~. 
\end{equation}

Putting
\begin{equation*}
  T ={\cal E}
  \left( 
    \begin{array}{cc}
      1&\go \\
      \go &1
    \end{array}
  \right),
\end{equation*}
we get:
\[
H= (1\pm \go)^2 {\cal E}^2 
\] 
and the quantisation rule:
\begin{equation}
  \ex^{i\big( (4l+2)\ga +I +O(h)\big)}= \frac{1}{2\pi}
  \Gamma (\ha +i \ge)^2 \ex^{\ge \big(\pi -2 i \log (l+\ha)\big)}
  (1\pm i \ex^{-\ge \pi})^2 ~,
\end{equation}
which has to be solved in $e$ where $e$ enters in $I$ and in $\ge$.

\begin{rema}
  Because of the $\pm$ sign in~\eqref{equ:cond2-reso12}, the spectrum
  can be separated in two spectra. For each of these, the spacings of
  eigenvalues are easily computed as in remark \ref{rema:ecarts}.
  Moreover, formula~\eqref{equ:cond2-reso12} shows that far from the
  critical value (ie. $e\fleche\pm\infty$), the ``$+$'' and ``$-$''
  eigenvalues associate in doublets, and that there is a universal
  transition happening when crossing the critical value ($e=0$), where
  a doublet ``$++$'' becomes a doublet ``$-+$''. The details for these
  formul\ae\ are similar to \cite{colin-p2} and left to the reader.
\end{rema}

\subsubsection{Matrix form for $\H_2$}
Since $V$ is a harmonic polynomial of degree $2$, the Toeplitz
operator $\H_2$ is given by
\[
\H_2 = \Pi_2V\Pi_2,
\]
where $V$ here is just the multiplication by $V$ and $\Pi_2$ is the
orthogonal projector on the space $\harm_2$ of spherical harmonics of
degree 2. We shall first determine the explicit formula for $\H_2$
with a generic $V\in \harm_2$ and then apply it to the specific form
$V=2ax^2+2by^2+2cz^2$.

The spaces $\harm_l$ are seen as the spaces of the irreducible
representation  of $\Lie{so}(3)$ on $L^2(S^2)$, acting via the
differential operators 
\[
\begin{array}{rcl}
L_x & = & y\deriv{}{z}-z\deriv{}{y}\\
L_y & = & z\deriv{}{x}-x\deriv{}{z}\\
L_z & = & z\deriv{}{x}-x\deriv{}{z},
\end{array}
\]
which are subject to the relation
\[
[L_x,L_y]=L_z \quad \textrm{ ( and cyclic permutations of } (x,y,z)
\textrm{ )}.
\] 
These operators commute with $\Delta$ and thus preserve $\harm_l$. As
usual, (see eg. \cite{weyl-qm}), we use the coordinates $\zeta=x+iy$
and $z$, and let
\[
\Omega_\pm = iL_x \pm L_y,
\]
so that $\Omega_+ = -\bar{\zeta}\deriv{}{z} + 2z\deriv{}{\zeta}$ and
$\Omega_- = \Omega_+^* = \zeta\deriv{}{z} - 2z\deriv{}{\bar{\zeta}}$.
A natural basis of $\harm_l$ is then the following:
\begin{equation}
  \label{equ:basis-harm}
  \Bas_l=(\zeta^l,\Omega_+\zeta^l,\Omega_+^2\zeta^l, \dots,
  \Omega_+^{m}\zeta^l, \dots, \Omega_+^{2l}\zeta^l).
\end{equation}

We shall use the convenient equivalent representation given by the
action of $\Lie{su}(2)$ on the spaces $\poly_{2l}$ of homogeneous
polynomials of degree $2l$ in $\CM^2=(\xi,\eta)$. Using the following
identification
\[
\begin{array}{rcl}
  L_x & = &
  \frac{1}{2i}\left(\eta\deriv{}{\xi}+\xi\deriv{}{\eta}\right)\\
  L_y & = &
  \frac{1}{2}\left(-\eta\deriv{}{\xi}+\xi\deriv{}{\eta}\right)\\
  L_z & = &
  \frac{1}{2i}\left(\xi\deriv{}{\xi}-\eta\deriv{}{\eta}\right),\\
\end{array}
\]
we get $\Omega_+=\xi\deriv{}{\eta}$, and a natural basis of
$\poly_{2l}$ is the following:
\[
(\eta^{2l},\Omega_+\eta^{2l},\Omega_+^2\eta^{2l},\dots,
\Omega_+^{m}\eta^{2l}, \dots, \Omega_+^{2l}\eta^{2l}).
\]
In the rest of the argument, this basis together with the basis
(\ref{equ:basis-harm}) will be used to identify $\harm_l$ and
$\poly_{2l}$. With this identification, $V$ assumes the form:
\[
V = \frac{a-b}{2}(\eta^4+\xi^4) - 3(a+b)\eta^2\xi^2.
\]

Up to a multiplicative constant (depending on $l$), there exists, for
each $l$, a unique equivariant morphism $\mbold{\Pi}:
\harm_2\otimes\harm_l\fleche \harm_l$. Hence $\mbold{\Pi}$ is a
multiple of $\mathcal{D}^2$, where
\[
\mathcal{D} = \deriv{}{\eta}\otimes\deriv{}{\xi} -
\deriv{}{\xi}\otimes\deriv{}{\eta}.
\]
(A subsequent composition by the multiplication $f\otimes g\fleche fg$
is always assumed.)

On the other hand, using the fact that every element of the form $fg$,
$f\in\harm_2$ and $g\in\harm_l$ splits into
\[
fg = f_{l+2}+r^2f_l + r^4f_{l-2},
\]
where $r$ is the radial distance and $f_j\in\harm_j$, one easily
computes
\[
f_l = \mbold{\Pi}(f\otimes g) = \frac{1}{6+4l}\left(\Delta_{\RM^3}(fg) -
  \frac{r^2}{2(2l-1)}\Delta_{\RM^3}^2(fg)\right).
\]
Testing this formula and $\mathcal{D}^2$ with, for instance,
$f=x^2-z^2$ and $g=(x+iy)^l$, one gets
\[
\mbold{\Pi} = -\frac{1}{6(2l-1)(3+2l)} \mathcal{D}^2
\]
\[
= -\frac{1}{6(2l-1)(3+2l)}\left(\deriv{^2}{\xi^2}\otimes
  \deriv{^2}{\eta^2} - 2\frac{\partial^2}{\partial \xi\partial
    \eta}\otimes \frac{\partial^2}{\partial \xi\partial \eta} +
  \deriv{^2}{\eta^2}\otimes\deriv{^2}{\xi^2}\right).
\]

Now it is easy to let this operator act on $V\harm_l$, and one finally
obtains the matrix representation of $\H_2$ in the basis $\Bas_l$
\begin{equation}
  \label{equ:matrix-s2}
  \mbox{$\H_2$}_{\restr \harm_l} = -\frac{1}{3(2l-1)(3+2l)}\left(
    \vcenter{\xymatrix@!R@=2pt@M=2pt{
        A_0\ar@{.}[dddddrrrrr] & 0 \ar@{--}[ddddrrrr] &
        B_0\ar@{.}[dddrrr] & & & \\ 
        0 \ar@{--}[ddddrrrr] &   & & & \textrm{\Large 0} &  \\
        B_{2l-2}\ar@{.}[dddrrr]  & & &  & &\\
        & & & & & B_{2l-2}\\
        & \textrm{\Large 0} &  & & & 0 \\
        & & &  B_{0} & 0 & A_{2l} }} \right),
\end{equation}
with
\[
\left\{\begin{array}{lcl}
    A_m & = &   6(a+b)(3m(2l-m)-l(2l-1))\\
    B_m & = & 3(a-b)(m+1)(m+2)
  \end{array}\right.
\]
Notice the symmetry of the matrix: $A_m=A_{2l-m}$.

\subsubsection{Numerical computations}
As it is, the matrix~\eqref{equ:matrix-reso12} is very badly prepared
for being numerically diagonalised. Indeed, the spectrum exhibits near
degeneracies -- as is expected from tunneling effects -- and usual
algorithms will rapidly fail as $l$ increases. Fortunately, there is
an easy way to cope with this, for two commuting transformations can
split the matrix: the projection onto the subspace spanned by the
vectors from the basis~\eqref{equ:basis-harm} having even index, and
the central symmetry ($m\fleche 2l-m$) of the matrix.  We arrive at
the following 4-blocks decomposition for the matrix
$-3(2l-1)(3+2l)\mbox{$\H_2$}_{\restr \harm_l}$, each block being a
tridiagonal matrix:
\begin{itemize}
\item if $l=2k+1$ is odd:
  \[
\left(
  \def\objectstyle{\scriptstyle}
  \vcenter{\xymatrix@!R@=0pt@M=1pt{
      A_0 & B_0 & & & &  \\
      B_{2l-2}  & A_2 \ar@{--}[dddrrr]& B_2\ar@{.}[dddrrr] & & 0 &  \\
      & B_{2l-4}\ar@{.}[dddrrr] & & & &\\
      & & & & & \\
      & 0 & & & A_{2k-2} & B_{2k-2}\\
      & & & & B_{2k+2} & (A_{2k}+B_{2k}) }} \right) \oplus 
\left(
  \def\objectstyle{\scriptstyle}
  \vcenter{\xymatrix@!R@=0pt@M=1pt{
      A_0 & B_0 & & & &  \\
      B_{2l-2}  & A_2 \ar@{--}[dddrrr]& B_2\ar@{.}[dddrrr] & & 0 &  \\
      & B_{2l-4}\ar@{.}[dddrrr] & & & &\\
      & & & & & \\
      & 0 & & & A_{2k-2} & B_{2k-2}\\
      & & & & B_{2k+2} & (A_{2k}-B_{2k}) }} \right)
\]
\[
\oplus \left(
  \def\objectstyle{\scriptstyle}
  \vcenter{\xymatrix@!R@=0pt@M=1pt{
      A_1 & B_1 & & & &  \\
      B_{2l-3}  & A_3 \ar@{--}[dddrrr]& B_3\ar@{.}[ddrr] & & 0 &  \\
      & B_{2l-5}\ar@{.}[dddrrr] & & & &\\
      & & & & B_{2k-3} & \\
      & 0 & & & A_{2k-1} & 2B_{2k-1}\\
      & & & & B_{2k+1} & A_{2k+1} }} \right) \oplus 
 \left(
  \def\objectstyle{\scriptstyle}
  \vcenter{\xymatrix@!R@=0pt@M=1pt{
      A_1 & B_1 & & & &  \\
      B_{2l-3}  & A_3 \ar@{--}[ddddrrrr]& B_3\ar@{.}[dddrrr] & & 0 &  \\
      & B_{2l-5}\ar@{.}[dddrrr] & & & &\\
      & & & & & \\
      & 0 & & &  & B_{2k-3}\\
      & & & & B_{2k+3} & A_{2k-1} }} \right).
\]
\item if $l=2k$ is even:
\[
\left(
  \def\objectstyle{\scriptstyle}
  \vcenter{\xymatrix@!R@=0pt@M=1pt{
      A_0 & B_0 & & & &  \\
      B_{2l-2}  & A_2 \ar@{--}[dddrrr]& B_2\ar@{.}[ddrr] & & 0 &  \\
      & B_{2l-4}\ar@{.}[dddrrr] & & & &\\
      & & & & B_{2k-4} & \\
      & 0 & & & A_{2k-2} & 2B_{2k-2}\\
      & & & & B_{2k} & A_{2k} }} \right) \oplus 
\left(
  \def\objectstyle{\scriptstyle}
  \vcenter{\xymatrix@!R@=0pt@M=1pt{
      A_0 & B_0 & & & &  \\
      B_{2l-2}  & A_2 \ar@{--}[ddddrrrr]& B_2\ar@{.}[dddrrr] & & 0 &  \\
      & B_{2l-4}\ar@{.}[dddrrr] & & & &\\
      & & & & & \\
      & 0 & & & & B_{2k-4}\\
      & & & & B_{2k+2} & A_{2k-2} }} \right)
\]
\[
\oplus \left(
  \def\objectstyle{\scriptstyle}
  \vcenter{\xymatrix@!R@=0pt@M=1pt{
      A_1 & B_1 & & & &  \\
      B_{2l-3}  & A_3 \ar@{--}[dddrrr]& B_3\ar@{.}[dddrrr] & & 0 &  \\
      & B_{2l-5}\ar@{.}[dddrrr] & & & &\\
      & & & &  & \\
      & 0 & & & A_{2k-3} & B_{2k-3}\\
      & & & & B_{2k+1} & (A_{2k-1}+B_{2k-1}) }} \right) \oplus 
 \left(
  \def\objectstyle{\scriptstyle}
  \vcenter{\xymatrix@!R@=0pt@M=1pt{
    A_1 & B_1 & & & &  \\
      B_{2l-3}  & A_3 \ar@{--}[dddrrr]& B_3\ar@{.}[dddrrr] & & 0 &  \\
      & B_{2l-5}\ar@{.}[dddrrr] & & & &\\
      & & & &  & \\
      & 0 & & & A_{2k-3} & B_{2k-3}\\
      & & & & B_{2k+1} & (A_{2k-1}-B_{2k-1}) }} \right).
\]
\end{itemize}

\begin{figure}[hbtp]
  \begin{center}
    \leavevmode \input{s2-compare.pstex_t}
    \caption{A comparison between \semicla\ and quantum computations. 
      Here is displayed the spectrum of $\mbox{$\H_2$}_{\restr
        \harm_l}$ in increasing order versus the eigenvalue number,
      for $l=40$ and the potential $V$ defined by $a=-1$ and $b=1$.
      The light crosses linked by line segments are the \semicla\ 
      computations while the quantum eigenvalues are the dark
      diamonds. We observe a very good accuracy (of mean order
      $O(h)=O(1/l)$ predicted by the theory, but much better near the
      critical value). Notice also how the eigenvalue doublets
      reassociate when passing through the critical value.}
    \label{fig:s2-compare}
  \end{center}
\end{figure}



\bibliographystyle{plain}
\bibliography{bibli}

\begin{thebibliography}{10}

\bibitem{ahlfors-sario}
L.V. Ahlfors and S.~Sario.
\newblock {\em Riemann Surfaces}.
\newblock Princeton University Press, 1960.

\bibitem{audin-qua}
M.~Audin.
\newblock Courbes alg{\'e}briques et syst{\`e}mes int{\'e}grables :
  g{\'e}od{\'e}siques des quadriques.
\newblock {\em Expositiones Math.}, 12:193--226, 1994.

\bibitem{cushman-book}
L.M. Bates and R.H. Cushman.
\newblock {\em Global aspects of classical integrable systems}.
\newblock Birha{\"u}ser, 1998.

\bibitem{weinstein-bates}
S.~Bates and A.~Weinstein.
\newblock {\em Lectures on the Geometry of Quantization}, volume~8 of {\em
  Berkeley Mathematics Lecture Notes}.
\newblock AMS, 1997.

\bibitem{child-book}
M.S. Child.
\newblock {\em Semiclassical Mechanics with Molecular Applications}.
\newblock Oxford University Press, 1991.

\bibitem{colin-bica}
Y.~Colin~de Verdi{\`e}re.
\newblock Sur le spectre des op{\'e}rateurs elliptiques {\`a}
  bicaract{\'e}ristiques toutes p{\'e}riodiques.
\newblock {\em Comment. Math. Helv.}, 54:508--522, 1979.

\bibitem{colinII}
Y.~Colin~de Verdi{\`e}re.
\newblock Spectre conjoint d'op{\'e}rateurs pseudo-diff{\'e}rentiels qui
  commutent {II}.
\newblock {\em Math. Z.}, 171:51--73, 1980.

\bibitem{colin-p}
Y.~Colin~de Verdi\`ere and B.~Parisse.
\newblock {\'E}quilibre instable en r{\'e}gime semi-classique {I}~:
  Concentration microlocale.
\newblock {\em Comm. Partial Differential Equations}, 19(9--10):1535--1563,
  1994.

\bibitem{colin-p2}
Y.~Colin~de Verdi\`ere and B.~Parisse.
\newblock {\'E}quilibre instable en r{\'e}gime semi-classique {II}~: Conditions
  de {B}ohr-{S}ommerfeld.
\newblock {\em Ann. Inst. H. Poincar{\'e}. Phys. Th{\'e}or.}, 61(3):347--367,
  1994.

\bibitem{colin-p3}
Y.~Colin~de Verdi\`ere and B.~Parisse.
\newblock Singular {B}ohr-{S}ommerfeld rules.
\newblock {\em Commun. Math. Phys.}, 205:459--500, 1999.

\bibitem{colin-vey}
Y.~Colin~de Verdi{\`e}re and J.~Vey.
\newblock Le lemme de {M}orse isochore.
\newblock {\em Topology}, 18:283--293, 1979.

\bibitem{darboux}
G.~Darboux.
\newblock {\em Th{\'e}orie g{\'e}n{\'e}rale des surfaces}.
\newblock Chelsea, 1972.

\bibitem{duistermaat-oscillatory}
J.J. Duistermaat.
\newblock Oscillatory integrals, {L}agrange immersions and unfoldings of
  singularities.
\newblock {\em Comm. Pure Appl. Math.}, 27:207--281, 1974.

\bibitem{fomenko}
A.T. Fomenko.
\newblock {\em Topological classification of integrable systems}, volume~6 of
  {\em Advances in soviet mathematics}.
\newblock AMS, 1991.

\bibitem{guillemin-sphere}
V.~Guillemin.
\newblock Some spectral results for the laplace operator with potential on the
  n-sphere.
\newblock {\em Adv. in Math.}, 27:273--286, 1978.

\bibitem{guillemin-rank1}
V.~Guillemin.
\newblock Some spectral results on rank one symmetric spaces.
\newblock {\em Adv. in Math.}, 28:129--137, 1978.

\bibitem{guillemin-band}
V.~Guillemin.
\newblock Band asymptotics in two dimensions.
\newblock {\em Adv. in Math.}, 42:248--282, 1981.

\bibitem{guillemin-schaeffer}
V.~Guillemin and D.~Schaeffer.
\newblock On a certain class of {F}uchsian partial differential equations.
\newblock {\em Duke Math. J.}, 44(1):157--199, 1977.

\bibitem{hirzebruch}
F.~Hirzebruch.
\newblock {\em Topological Methods in Algebraic Geometry}, volume 131 of {\em
  Grundlehren der math. {W}.}
\newblock Springer, New York, 1966.

\bibitem{klingenberg}
W.~Klingenberg.
\newblock {\em Riemannian Geometry}.
\newblock de Gruyter, 1982.

\bibitem{moser-quadrics}
J.~Moser.
\newblock Geometry of quadrics and spectral theory.
\newblock In {\em The Chern Symposium}, pages 147--188. Springer, 1980.

\bibitem{zung-torus}
Z.~Nguy{\^e}n~Ti{\^e}n.
\newblock Singularities of integrable geodesic flows on multidimensional torus
  and sphere.
\newblock {\em J. Geom. Phys.}, 18:147--162, 1996.

\bibitem{zung-I}
Z.~Nguy{\^e}n~Ti{\^e}n.
\newblock Symplectic topology of integrable hamiltonian systems, {I}:
  {A}rnold-{L}iouville with singularities.
\newblock {\em Compositio Math.}, 101:179--215, 1996.

\bibitem{zung-al}
Z.~Nguy{\^e}n~Ti{\^e}n, L.S. Polyakova, and E.N. Selianova.
\newblock Topological classification of integrable geodesic flows on orientable
  two-dimensional manifolds...
\newblock {\em Functional Anal. Appl.}, 27:186--196, 1993.

\bibitem{san-fn}
S.~V{\~u}~Ng{\d o}c.
\newblock Formes normales semi-classiques des syst{\`e}mes compl{\`e}tement
  int{\'e}grables au voisinage d'un point critique de l'application moment.
\newblock Preprint Institut Fourier 377, \texttt{http://www-\ellastiq
  fourier.ujf-\ellastiq grenoble.fr/\ellastiq}
  \verb+~+\texttt{svungoc/\ellastiq fn-tout.ps.gz}, 1997.

\bibitem{san-these}
S.~V{\~u}~Ng{\d o}c.
\newblock {\em Sur le spectre des syst{\`e}mes compl{\`e}tement int{\'e}grables
  semi-classiques avec singularit{\'e}s}.
\newblock PhD thesis, Universit{\'e} Grenoble 1, 1998.

\bibitem{san-focus}
S.~V{\~u}~Ng{\d o}c.
\newblock Bohr-{S}ommerfeld conditions for integrable systems with critical
  manifolds of focus-focus type.
\newblock {\em Comm. Pure Appl. Math.}, 53(2):143--217, 2000.

\bibitem{weinstein}
A.~Weinstein.
\newblock {\em Lectures on Symplectic Manifolds}.
\newblock Number~29 in Regional Conference Series in Mathematics. AMS, 1976.

\bibitem{weinstein-cluster}
A.~Weinstein.
\newblock Asymptotics of eigenvalue clusters for the laplacian plus a
  potential.
\newblock {\em Duke Math. J.}, 44(4):883--892, 1977.

\bibitem{weyl-qm}
H.~Weyl.
\newblock {\em The theory of groups and quantum mechanics}.
\newblock Dover, 1950.
\newblock Translated from the (second) german edition.

\end{thebibliography}

\end{document}